\documentclass[a4paper]{article}%
\usepackage{graphicx}
\usepackage{geometry}
\usepackage{amsmath}
\usepackage{amssymb}
\usepackage{amsfonts}
\usepackage{amsthm,amscd}
\usepackage[colorlinks=true]{hyperref}%
\providecommand{\U}[1]{\protect\rule{.1in}{.1in}}
\hypersetup{urlcolor=blue, citecolor=red}
\def\figurename{Figure}
\makeatletter
\renewcommand{\fnum@figure}[1]{\figurename~\thefigure.}
\makeatother
\def\tablename{Table}
\makeatletter
\renewcommand{\fnum@table}[1]{\tablename~\thetable.}
\makeatother
\def \bop {\noindent\textbf{Proof. }}
\def \eop {\hbox{}\nobreak\hfill
\vrule width 2mm height 2mm depth 0mm
\par \goodbreak \smallskip}
\newtheorem{theorem}{Theorem}[section]
\newtheorem{lemma}[theorem]{Lemma}

\newtheorem{proposition}[theorem]{Proposition}
\theoremstyle{definition}
\newtheorem{definition}[theorem]{Definition}

\theoremstyle{remark}
\newtheorem{remark}[theorem]{Remark}
\numberwithin{equation}{section}

\begin{document}

\title{A PDE approach for open-loop equilibriums in time-inconsistent stochastic
optimal control problems}
\author{\textbf{Ishak Alia}{ }\thanks{Department of Mathematics University of Bordj
Bou Arreridj 34000 Algeria. E-mail: izacalia@yahoo.com}}
\maketitle

\begin{abstract}
This paper studies open-loop equilibriums for a general class of
time-inconsistent stochastic control problems under jump-diffusion SDEs with
deterministic coefficients. Inspired by the idea of Four-Step-Scheme for
forward-backward stochastic differential equations with jumps (FBSDEJs, for
short), we derive two systems of integro-partial differential equations
(IPDEs, for short). Then, we rigorously prove a verification theorem which
provides a sufficient condition for open-loop equilibrium strategies. As an
illustration of the general theory, we discuss a mean-variance portfolio
selection problem under a jump-diffusion model.

\end{abstract}

\textbf{Keys words}: Time-inconsistency, stochastic control, \ open-loop
equilibrium control, Poisson point process, parabolic integro-partial
differential equations, mean-variance problem.


\section{Introduction\label{section1}}

Traditional optimal control problems are time-consistent in the sense that an
optimal control constructed for a given initial pair of time and state will
remain optimal thereafter. Time-consistency provides a powerful approach to
solving optimal control problems, namely, the method of \textit{dynamic
programming}. The basic idea of this method is to consider a family of
(time-consistent)\textit{ }optimal control problems with different initial
times and states, to establish relationships among these problems via the
so-called Hamilton-Jacobi-Bellman equation (HJB, for short). If the HJB
equation is solvable, then one can obtain an optimal control in closed-loop
form by taking the maximizer/minimizer of the Hamiltonian involved in the HJB
equation; see e.g. the textbook \cite{Yong} for a detailed discussion.

However, there are growing evidences which tend to support that the
time-consistency could be lost in real world situations (i.e., an optimal
control selected at a given moment might not remain optimal at later time
moments). Among many possible reasons causing the time-inconsistency, there
are two playing some essential roles: (i) People are impatient about choices
in the short term but are patient when choosing between long-term alternatives
(see e.g. \cite{Loewenstein}); (ii) and people's attitude regarding risks is
subjective rather than objective, meaning that, different (group of) people
will have different opinions on the risks that contained in some coming event
(see e.g. \cite{Yan}). Mathematically, the first situation can be described by
the general discounting, and the second situation can lead to a certain
nonlinear appearance of conditional expectations for the state process and/or
control process in the cost functional. Two well-known examples are Merton's
portfolio problem with non-exponential discounting \cite{EkePir} and
continuous-time mean-variance (MV, for short) portfolio selection model
\cite{ZhouLi}.

The theory of time-inconsistent optimal control problems was pioneered by R.H.
Strotz in his work \cite{Strotz} on a general discounting Ramsay problem. It
became now very popular, and is an important field of research due to its
application in mathematical finance. To deal with the time-inconsistency,
Strotz \cite{Strotz} suggested the so-called \textit{consistent planning
}approach\textit{.} The basic idea of this method is that an action the
controller selects at every instant of time is considered as a non-cooperative
game against all the actions the controller is going to make in the future. A
Nash equilibrium strategy is therefore a decision such that any deviation from
it at any time instant will be worse off. Further work along this line in
continuous and discrete time can be found in \cite{Goldman}, \cite{Phelps},
\cite{Krusell}, \cite{Pollak} and \cite{Kydland}. However, at that time, the
results in continuous-time setting were essentially obtained under intuitive
descriptions of the Nash equilibrium concept. The first precise definition of
the \textit{equilibrium strategy} in continuous-time setting was introduced in
the papers of Ekeland and Lazrak \cite{EkeLaz} and Ekeland and Pirvu
\cite{EkePir}, where the source of inconsistency is general discounting. In
the series of works carried out by Bj\"{o}rk and Murgoci \cite{Bjork1} and
Bj\"{o}rk et al. \cite{Bjork4}, the authors investigated equilibrium
strategies within the class of closed-loop controls for a general class of
time inconsistent Markovian problems. Inspired by the discrete case, they
first derived in an heuristic way the so-called \textit{extended HJB
equations} (that is a system of three-coupled fully nonlinear parabolic PDEs),
and then rigorously proved a verification theorem. Recently, Lindensj\"{o}
\cite{Lindensjo1} derived rigourously the extended HJB equations without using
arguments from the discrete-time case. He and Jiang \cite{He} and
Hern\'{a}ndez and Possama\"{\i} \cite{Possamai} generalized \cite{Bjork4} by
refining the definition of the closed-loop equilibrium concept. For
equilibrium stopping times in time-inconsistent Markovian Problems, see e.g.
\cite{Lindensjo2}, \cite{Lindensjo3} and references therein.

In the series of works carried out by Yong [\cite{Yong2011}, \cite{Yong2012b},
\cite{Yong2012}, \cite{Yong2013}], the author studied a general class of time
inconsistent optimal control problems. He constructed closed-loop equilibrium
strategies in a multi-person differential game framework with a hierarchical
structure and introduced the so-called \textit{equilibrium HJB equation}; a
new type of\textit{ forward-backward Riccati-Volterra equations} were also
introduced. Yong's\ approach has inspired many studies over the last five
years. Let us just mention a few: see \cite{Wang2} for controlled systems with
random coefficients, see \cite{Wei2}, \cite{Mei} for stochastic systems with
Markov regime switching and see \cite{Wei1} for time-inconsistent recursive
utilities. More recently, Wang and Yong \cite{WangYong} considered the case
where the cost functional is determined by a backward stochastic Volterra
integral equation which covers most of the existence literature on the general
discounting situation. Dou and L\"{u} \cite{Dou} dealt with closed-loop
equilibrium controls in the framework of linear-quadratic optimal control
problems with time-inconsistent cost functionals for stochastic evolution
equations in a Hilbert space.

The extended HJB-method and Yong's\ differential game approach are extensions
of the classical dynamic programming approach for the determination of
closed-loop equilibrium strategies. In contrast to the above mentioned
literature, Hu et al. (\cite{Huetal2011}, \cite{Huetal2017}) introduced the
concept of \textit{open-loop equilibrium control }by using a spike variation
formulation, which is different from the closed-loop equilibrium concepts. The
authors considered a stochastic linear-quadratic model (SLQ, for short) in a
non-Markovian framework, where the time-inconsistency arises from the presence
of a quadratic term of the expected state as well as a state-dependent term in
the objective functional. They used a variational method in the spirit of
Peng's stochastic maximum principle to characterize the existence and
uniqueness of the equilibrium control in terms of the solvability of a
\textit{"flow"} of FBSDEs; we refer the readers to the work of Hamaguchi
\cite{Hamaguchi1} for a detailed discussion about this new type of FBSDEs.
Some recent studies devoted to the open-loop equilibrium concept can be found
in \cite{Djehiche}, \cite{Alia1}, \cite{Alia3}, \cite{Hamaguchi2},
\cite{Wang}, \cite{Wang11}, \cite{Wang12}, \cite{Huetal2019}, \cite{Sun} and
\cite{Zhang}. Specially, Wang \cite{Wang} discussed open-loop equilibrium
controls and their particular closed-loop representations for a class of SLQ
problems under mean-field SDEs. Alia \cite{Alia1} investigated open-loop
equilibrium controls in general discounting situations. Sun and Guo \cite{Sun}
generalized \cite{Huetal2017} to the case of SLQ models with random jumps.
Wang et al. \cite{Wang11} found the open-loop equilibrium strategy to a
mean-variance problem under a non-markovian regime-switching model. Hamaguchi
\cite{Hamaguchi2} performed a sophisticated FBSDEs-approach and characterized
the unique open-loop equilibrium solution to a general discounting Merton
portfolio problem in an incomplete market with random parameters. Finally, in
the work of Hamaguchi \cite{Hamaguchi3}, the author investigated open-loop
equilibrium controls in the framework of time-inconsistent stochastic
recursive control problems, where the cost functional is defined by the
solution to a backward stochastic Volterra integral equation. He managed, by a
number of very clever ideas, to derive a necessary and sufficient condition
for an open-loop equilibrium control via variational methods.

In this paper, we suggest a new PDEs-approach to characterize open-loop
equilibrium controls for a general class of time-inconsistent stochastic
control problems with deterministic coefficients. We combine the ideas from
the generalized HJB-approaches [\cite{Bjork4}, \cite{Yong2012}] and
forward-backward stochastic differential equations with jumps \cite{MY}. In
the following, we provide a brief outline of our approach:

\begin{enumerate}
\item Due to the Markovian structure of our problem\footnote{All the
coefficients in our model are deterministic functions.}, we introduce the
concept of \textit{open-loop equilibrium strategy }(see Definition \ref{Def2}
below) which can be seen as a closed-loop representation of an open-loop
equilibrium control; see e.g. \cite{Wang} and \cite{Wang12}.

\item Combining the spike perturbation of the equilibrium strategy with the
idea of Four-Step-Scheme for FBSDEs with jumps introduced in \cite{MY}, we
derive two systems of parabolic integro-partial differential equations; we
point out that this type of IPDEs appears for the first time in the literature.

\item We then rigorously prove a verification theorem by which one can
construct an open-loop equilibrium strategy by solving the above-mentioned IPDEs.
\end{enumerate}

Recently, Alia \cite{Alia2} investigated open-loop equilibriums for a similar
class of time-inconsistent stochastic control problems by using a backward
stochastic partial differential equations (BSPDEs, for short) approach. He
derived a verification theorem that provides a general sufficient condition
for equilibriums via a coupled stochastic system of a single forward SDE and a
flow of BSPDEs. Note that the author introduced BSPDEs in his approach while
all the involved coefficients are deterministic. In contrast to the above
mentioned paper, here we derive a verification theorem by using deterministic
IPDEs instead of BSPDEs. Moreover, it is worth mentioning that compared with
the existing literature, the previous proposed approach demonstrates several
new advantages on the treatment of open-loop equilibrium controls; see Remark
\ref{Rem1} below for more details.

The plan of the paper is as follows, in the second section, we give necessary
notations and some preliminaries on FBSDEs with jumps. In Section
\ref{section3}, we formulate our time-inconsistent stochastic optimal control
problem. Section \ref{section4} is devoted to the presentation of the
verification theorem. In Section \ref{section5}, we discuss the connections
between the PDEs-approach of the present article and the variational approach
of Hu et al. (\cite{Huetal2011}, \cite{Huetal2017}). Finally, in Section
\ref{section6}, we discuss a mean-variance portfolio selection model.

\section{Preliminaries\label{section2}}

\subsection{Notations}

Throughout this paper $(\Omega,\mathcal{F},\left(  \mathcal{F}_{t}\right)
_{t\in\left[  0,T\right]  },\mathbb{P})$ is a filtered probability space such
that $\mathcal{F}_{0}$ contains all $\mathbb{P}$-null sets, $\mathcal{F}%
_{T}=\mathcal{F}$ for an arbitrarily fixed finite time horizon $T>0$, and
$\left(  \mathcal{F}_{t}\right)  _{t\in\left[  0,T\right]  }$ satisfies the
usual conditions. We assume that $\left(  \mathcal{F}_{t}\right)
_{t\in\left[  0,T\right]  }$ is generated by a d-dimensional standard Brownian
motion $\left(  W\left(  t\right)  \right)  _{t\in\left[  0,T\right]  }$ and
an independent Poisson measure $N$ on $\left[  0,T\right]  \times E$ where $E$
$\subseteq$ $%
\mathbb{R}
\backslash\left\{  0\right\}  $. We assume that the compensator of $N$ has the
form $\mu\left(  ds,de\right)  =\vartheta\left(  de\right)  ds$ for some
positive and finite L\'{e}vy measure on $E$, endowed with its Borel $\sigma
$-field $\mathcal{B}\left(  E\right)  $. We write $\tilde{N}\left(
ds,de\right)  =N\left(  ds,de\right)  -\vartheta\left(  de\right)  ds$ for the
compensated jump martingale random measure of $N$. Obviously, we have
$\mathcal{F}_{t}=\sigma\left[
{\textstyle\int}
{\textstyle\int_{\left(  0,s\right]  \times A}}
N\left(  dr,de\right)  ;s\leq t,\text{ }A\in\mathcal{B}\left(  E\right)
\right]  \vee\sigma\left[  W\left(  s\right)  ;s\leq t\right]  \vee
\mathcal{N}$, \ where $\mathcal{N}$ denotes the totality of $\vartheta$-null
sets, and $\sigma_{1}\vee\sigma_{2}$\ denotes the $\sigma$-field generated by
$\sigma_{1}\cup\sigma_{2}.$

We use $H^{\top}$ to denote the transpose of any vector or matrix $H$ and
$\chi_{A}$ to denote the indicator function of the set $A$. For a function $f$
we denote by $\nabla f\left(  x,z\right)  $ (resp. $\nabla^{2}f\left(
x,y\right)  $) the gradient (resp. the Hessian) of $f$ with respect to the
variables $\left(  x,z\right)  \in\mathbb{R}^{n}\times\mathbb{R}^{n}$.
Particularly, $f_{y}\left(  x,z\right)  $ (resp. $f_{yy}\left(  x,z\right)  $)
denotes the first derivative (resp. the second derivative) of $f$ with respect
to the variable $y=x$,$z$. When $\left(  X\left(  s\right)  \right)
_{s\in\left[  t,T\right]  }$ is a c\`{a}dl\`{a}g processes, we define the
process $\left(  X_{-}\left(  s\right)  \right)  _{s\in\left[  t,T\right]  }$
by%
\[
X_{-}\left(  t\right)  =X\left(  t\right)  \text{ and }X_{-}\left(  s\right)
=\lim_{r\uparrow s}X\left(  r\right)  \text{, for }s\left(  t,T\right]
\text{.}%
\]
In addition, we use the following notations for several sets and spaces of
processes on the filtered probability space, which will be used later:

\begin{enumerate}
\item[$\bullet$] $D\left[  0,T\right]  =\left\{  \left(  t,s\right)
\in\left[  0,T\right]  \times\left[  0,T\right]  \text{, such that }s\geq
t\right\}  $.

\item[$\bullet$] $\mathbb{S}^{n}:$ the set of $\left(  n\times n\right)  $
symmetric matrices.

\item[$\bullet$] $\mathbb{L}^{q}\left(  \Omega,\mathcal{F}_{t},\mathbb{P};%
\mathbb{R}
^{n}\right)  :$the set of $\mathbb{R}^{n}$-valued, $\mathcal{F}_{t}%
$-measurable random variables $\zeta$, with%
\[
\left\Vert \zeta\right\Vert _{\mathbb{L}^{q}\left(  \Omega,\mathcal{F}%
_{t},\mathbb{P};%
\mathbb{R}
^{n}\right)  }^{q}=\mathbb{E}\left[  \left\vert \zeta\right\vert ^{q}\right]
<\infty.
\]

\item[$\bullet$] $\mathbb{L}^{q}\left(  E,\mathcal{B}\left(  E\right)
,\vartheta;%
\mathbb{R}
^{n}\right)  :$the space of functions $r:E\rightarrow%
\mathbb{R}
^{n}$ such that $r\left(  \cdot\right)  $ is $\mathcal{B}\left(  E\right)
-$measurable, with%
\[
\left\Vert r\left(  \cdot\right)  \right\Vert _{\mathbb{L}^{q}\left(
E,\mathcal{B}\left(  E\right)  ,\vartheta;%
\mathbb{R}
^{n}\right)  }^{q}=%
{\displaystyle\int_{E}}
\left\vert r\left(  e\right)  \right\vert ^{q}\vartheta\left(  de\right)
<\infty.
\]

\item[$\bullet$] $\mathcal{S}_{%
\mathcal{F}%
}^{q}\left(  t,T;%
\mathbb{R}
^{n}\right)  :$ the space of$\
\mathbb{R}
^{n}-$valued, $\left(  \mathcal{F}_{s}\right)  _{s\in\left[  t,T\right]  }%
$-adapted c\`{a}dl\`{a}g processes $X\left(  \cdot\right)  $, with%
\[
\left\Vert X\left(  \cdot\right)  \right\Vert _{\mathcal{S}_{%
\mathcal{F}%
}^{q}\left(  t,T;%
\mathbb{R}
^{n}\right)  }^{q}=\mathbb{E}\left[  \sup\limits_{s\in\left[  t,T\right]
}\left\vert X\left(  s\right)  \right\vert ^{q}\right]  <\infty.
\]

\item[$\bullet$] $\mathcal{L}_{%
\mathcal{F}%
}^{2}\left(  t,T;%
\mathbb{R}
^{n\times d}\right)  :$ the space of $\
\mathbb{R}
^{n\times d}-$valued, $\left(  \mathcal{F}_{s}\right)  _{s\in\left[
t,T\right]  }-$adapted processes $Q\left(  \cdot\right)  ,$ with%
\[
\left\Vert Q\left(  \cdot\right)  \right\Vert _{\mathcal{L}_{%
\mathcal{F}%
}^{2}\left(  t,T;%
\mathbb{R}
^{n\times d}\right)  }^{2}=\mathbb{E}\left[
{\displaystyle\int\limits_{t}^{T}}
\text{\textbf{tr}}\left[  Q\left(  s\right)  ^{\top}Q\left(  s\right)
\right]  ds\right]  <\infty.
\]

\item[$\bullet$] $\mathcal{L}_{%
\mathcal{F}%
,p}^{q}\left(  t,T;%
\mathbb{R}
^{l}\right)  :$ the space of $\
\mathbb{R}
^{l}-$valued, $\left(  \mathcal{F}_{s}\right)  _{s\in\left[  t,T\right]  }%
-$predictable processes $u\left(  \cdot\right)  ,$ with%
\[
\left\Vert u\left(  \cdot\right)  \right\Vert _{\mathcal{L}_{%
\mathcal{F}%
,p}^{q}\left(  t,T;%
\mathbb{R}
^{l}\right)  }^{q}=\mathbb{E}\left[
{\displaystyle\int\limits_{t}^{T}}
\left\vert u\left(  s\right)  \right\vert ^{q}ds\right]  <\infty.
\]

\item[$\bullet$] $\mathcal{L}_{%
\mathcal{F}%
,p}^{\vartheta,q}\left(  \left[  t,T\right]  \times E;%
\mathbb{R}
^{n}\right)  $ the space of $%
\mathbb{R}
^{n}$-valued, $\left(  \mathcal{F}_{s}\right)  _{s\in\left[  t,T\right]  }%
$-predictable processes $R\left(  \cdot,\cdot\right)  $,$\ $with%
\[
\left\Vert R\left(  \cdot,\cdot\right)  \right\Vert _{\mathcal{L}_{%
\mathcal{F}%
,p}^{\vartheta,q}\left(  \left[  t,T\right]  \times E;%
\mathbb{R}
^{n}\right)  }^{q}=\mathbb{E}\left[
{\displaystyle\int_{t}^{T}}
{\displaystyle\int_{E}}
\left\vert R\left(  s,e\right)  \right\vert ^{q}\vartheta\left(  de\right)
ds\right]  <\infty\text{.}%
\]

\item[$\bullet$] $\mathcal{C}^{1,2}\left(  \left[  t,T\right]  \times%
\mathbb{R}
^{n}\times%
\mathbb{R}
^{n};%
\mathbb{R}
\right)  :$ the space of functions $V:\left[  t,T\right]  \times%
\mathbb{R}
^{n}\times%
\mathbb{R}
^{n}\rightarrow%
\mathbb{R}
$, such that $V$, $V_{s}$, $\nabla V$ and $\nabla^{2}V$ are continuous.
\end{enumerate}

\subsection{FBSDEs with jumps and IPDEs}

As preparations, in this paragraph we investigate the connection between a
specific class of FBSDEs with jumps and IPDEs. Let $\bar{\mu}:\left[
0,T\right]  \times\mathbb{R}^{n}\times\mathbb{R}^{n}\rightarrow\mathbb{R}^{n}%
$, $\bar{\sigma}:\left[  0,T\right]  \times\mathbb{R}^{n}\times\mathbb{R}%
^{n}\rightarrow\mathbb{R}^{n\times d}$, $\bar{c}:\left[  0,T\right]
\times\mathbb{R}^{n}\times\mathbb{R}^{n}\times E\rightarrow\mathbb{R}^{n}$,
$\bar{f}:\left[  0,T\right]  \times\mathbb{R}^{n}\times\mathbb{R}%
^{n}\rightarrow\mathbb{R}$ and $\bar{F}:\mathbb{R}^{n}\rightarrow\mathbb{R}$
be deterministic measurable functions.

For some fixed $\left(  t_{0},\xi,\eta\right)  \in\left[  0,T\right]
\times\mathbb{L}^{4}\left(  \Omega,\mathcal{F}_{t_{0}},\mathbb{P};%
\mathbb{R}
^{n}\right)  \times\mathbb{L}^{4}\left(  \Omega,\mathcal{F}_{t_{0}}%
,\mathbb{P};%
\mathbb{R}
^{n}\right)  $, consider in time interval $\left[  t_{0},T\right]  $ the
following coupled system of FBSDEs,%
\begin{equation}
\left\{
\begin{array}
[c]{l}%
dZ\left(  s\right)  =\bar{\mu}\left(  s,Z\left(  s\right)  ,Z\left(  s\right)
\right)  ds+\bar{\sigma}\left(  s,Z\left(  s\right)  ,Z\left(  s\right)
\right)  dW\left(  s\right) \\
\text{ \ \ \ \ \ \ \ \ \ \ \ \ }+%
{\displaystyle\int_{E}}
\bar{c}\left(  s,Z_{-}\left(  s\right)  ,Z_{-}\left(  s\right)  ,e\right)
\tilde{N}\left(  ds,de\right)  \text{,}\\
dX\left(  s\right)  =\bar{\mu}\left(  s,X\left(  s\right)  ,Z\left(  s\right)
\right)  ds+\bar{\sigma}\left(  s,X\left(  s\right)  ,Z\left(  s\right)
\right)  dW\left(  s\right) \\
\text{ \ \ \ \ \ \ \ \ \ \ \ }+%
{\displaystyle\int_{E}}
\bar{c}\left(  s,X_{-}\left(  s\right)  ,Z_{-}\left(  s\right)  ,e\right)
\tilde{N}\left(  ds,de\right)  \text{,}\\
dY\left(  s\right)  =-\bar{f}\left(  s,X\left(  s\right)  ,Z\left(  s\right)
\right)  ds+Q\left(  s\right)  ^{\top}dW\left(  s\right) \\
\text{ \ \ \ \ \ \ \ \ \ \ \ \ \ \ }+%
{\displaystyle\int_{E}}
R\left(  s,e\right)  \tilde{N}\left(  ds,de\right)  \text{, for }s\in\left[
t_{0},T\right]  ,\\
Z\left(  t_{0}\right)  =\xi\text{, }X\left(  t_{0}\right)  =\eta\text{,
}Y\left(  T\right)  =\bar{F}\left(  X\left(  T\right)  \right)  \text{.}%
\end{array}
\right.  \label{Eq01}%
\end{equation}

We introduce the following assumptions.

\begin{enumerate}
\item[\textbf{(A1)}] The maps $\bar{\mu}$, $\bar{\sigma}$ and $\bar{c}$ are
continuous and there exists a constant $C>0$ such that for all $s\in\left[
0,T\right]  $, $e\in E$, $\left(  x,x^{\prime},z,z^{\prime}\right)  \in\left(
%
\mathbb{R}
^{n}\right)  ^{4}$, $q\geq2$%
\[%
\begin{array}
[c]{rr}%
\left\vert \bar{\mu}\left(  s,x,z\right)  -\bar{\mu}\left(  s,x^{\prime
},z^{\prime}\right)  \right\vert +\left\vert \bar{\sigma}\left(  s,x,z\right)
-\bar{\sigma}\left(  s,x^{\prime},z^{\prime}\right)  \right\vert  & \\
+\left\vert \bar{c}\left(  s,x,z,e\right)  -\bar{c}\left(  s,x^{\prime
},z^{\prime},e\right)  \right\vert  & \leq C\left(  \left\vert x-x^{\prime
}\right\vert +\left\vert z-z^{\prime}\right\vert \right)
\end{array}
\]
and%
\[
\left\vert \bar{\mu}\left(  s,x,z\right)  \right\vert +\left\vert \bar{\sigma
}\left(  s,x,z\right)  \right\vert +\left(
{\displaystyle\int_{E}}
\left\vert \bar{c}\left(  s,x,z,e\right)  \right\vert ^{q}\vartheta\left(
de\right)  \right)  ^{\frac{1}{q}}\leq C\left(  1+\left\vert x\right\vert
+\left\vert z\right\vert \right)  \text{. }%
\]

\item[\textbf{(A2)}] The maps $\bar{f}$ and $\bar{F}$ are continuous and
quadratic growth on $\left(  x,z\right)  $ uniformly in time, i.e. there
exists a constant $C>0$ such that for all $\left(  s,x,z\right)  \in\left[
0,T\right]  \times%
\mathbb{R}
^{n}\times%
\mathbb{R}
^{n}$,%
\[
\left\vert \bar{f}\left(  s,x,z\right)  \right\vert \leq C\left(  1+\left\vert
x\right\vert ^{2}+\left\vert z\right\vert ^{2}\right)
\]
and%
\[
\left\vert \bar{F}\left(  x\right)  \right\vert \leq C\left(  1+\left\vert
x\right\vert ^{2}\right)  \text{.}%
\]

\end{enumerate}

Under Assumptions \textbf{(A1)-(A2)}, it is not difficult to verify that the
forward SDEs in (\ref{Eq01}) are uniquely solvable in $\left(  Z\left(
\cdot\right)  ,X\left(  \cdot\right)  \right)  \in\mathcal{S}_{%
\mathcal{F}%
}^{4}\left(  t_{0},T;%
\mathbb{R}
^{n}\right)  ^{2}$, and the BSDEs in (\ref{Eq01}) admits a unique adapted
solution $\left(  Y\left(  \cdot\right)  ,Q\left(  \cdot\right)  ,R\left(
\cdot,\cdot\right)  \right)  \in\mathcal{S}_{%
\mathcal{F}%
}^{2}\left(  t_{0},T;%
\mathbb{R}
\right)  \times\mathcal{L}_{%
\mathcal{F}%
}^{2}\left(  t_{0},T;%
\mathbb{R}
^{d}\right)  \times\mathcal{L}_{%
\mathcal{F}%
,p}^{\vartheta,q}\left(  \left[  t_{0},T\right]  \times E;%
\mathbb{R}
\right)  $; see e.g [\cite{Romuald}, pp. 125-126].

Inspired by the idea of Four Step Scheme introduced in \cite{MY} for FBSDEs
with jumps, we are going to show that the family of FBSDEs (\ref{Eq01}) is
closely linked to the following linear integro-partial differential equation:%
\begin{equation}
\left\{
\begin{array}
[c]{l}%
0=\Theta_{s}\left(  s,x,z\right)  +\left\langle \Theta_{x}\left(
s,x,z\right)  ,\bar{\mu}\left(  s,x,z\right)  \right\rangle +\left\langle
\Theta_{z}\left(  s,x,z\right)  ,\bar{\mu}\left(  s,z,z\right)  \right\rangle
\\
\text{ \ \ }+\frac{1}{2}\text{\textbf{tr}}\left[  \bar{\sigma}\bar{\sigma
}^{\top}\left(  s,x,z\right)  \Theta_{xx}\left(  s,x,z\right)  \right]
+\frac{1}{2}\text{\textbf{tr}}\left[  \bar{\sigma}\bar{\sigma}^{\top}\left(
s,z,z\right)  \Theta_{zz}\left(  s,x,z\right)  \right] \\
\text{ \ }+\text{\textbf{tr}}\left[  \bar{\sigma}\left(  s,x,z\right)  ^{\top
}\Theta_{xz}\left(  s,x,z\right)  \bar{\sigma}\left(  s,z,z\right)  \right]
+\text{\ }\bar{f}\left(  s,x,z\right) \\
\text{ \ \ }+%
{\displaystyle\int_{E}}
\left\{  \Theta\left(  s,x+\bar{c}\left(  s,x,z,e\right)  ,z+\bar{c}\left(
s,z,z,e\right)  \right)  -\Theta\left(  s,x,z\right)  \right. \\
\text{\ }%
\begin{array}
[c]{r}%
\left.  -\left\langle \Theta_{x}\left(  s,x,z\right)  ,\bar{c}\left(
s,x,z,e\right)  \right\rangle -\left\langle \Theta_{z}\left(  s,x,z\right)
,\bar{c}\left(  s,z,z,e\right)  \right\rangle \right\}  \vartheta\left(
de\right)  \text{, }\\
\text{for }\left(  s,x,z\right)  \in\left[  t_{0},T\right]  \times
\mathbb{R}^{n}\times\mathbb{R}^{n}\text{,}%
\end{array}
\\
\Theta\left(  T,x,z\right)  =\bar{F}\left(  x\right)  \text{, for }\left(
x,z\right)  \in\mathbb{R}^{n}\times\mathbb{R}^{n}\text{.}%
\end{array}
\right.  \label{Eq02}%
\end{equation}

The following result will be used later.

\begin{theorem}
\label{result00}Assume that \textbf{(A1)-(A2)} are satisfied. Let $\left(
X\left(  \cdot\right)  ,Z\left(  \cdot\right)  ,Y\left(  \cdot\right)
,Q\left(  \cdot\right)  ,R\left(  \cdot,\cdot\right)  \right)  $ be the unique
adapted solution of FBSDEs (\ref{Eq01}) and suppose that $\Theta\left(
\cdot,\cdot,\cdot\right)  \in\mathcal{C}^{1,2,2}\left(  \left[  t_{0}%
,T\right]  \times%
\mathbb{R}
^{n}\times%
\mathbb{R}
^{n};%
\mathbb{R}
\right)  $ is a classical solution to the IPDE (\ref{Eq02}) such that the
following hold%
\begin{align}
\int_{t_{0}}^{T}\mathbb{E}\left[  \left\vert \bar{\sigma}\left(  s,Z\left(
s\right)  ,Z\left(  s\right)  \right)  ^{\top}\Theta_{z}\left(  s,X\left(
s\right)  ,Z\left(  s\right)  \right)  \right\vert ^{2}\right]  ds  &
<\infty\text{,}\label{Eq06}\\
\int_{t_{0}}^{T}\mathbb{E}\left[  \left\vert \bar{\sigma}\left(  s,X\left(
s\right)  ,Z\left(  s\right)  \right)  ^{\top}\Theta_{x}\left(  s,X\left(
s\right)  ,Z\left(  s\right)  \right)  \right\vert ^{2}\right]  ds  &
<\infty\text{,} \label{Eq07}%
\end{align}
and%
\begin{equation}
\mathbb{E}\left[  \int\limits_{t_{0}}^{T}%
{\displaystyle\int\limits_{E}}
\left\vert \Theta\left(  s,X\left(  s\right)  +\bar{c}_{1}\left(  s,e\right)
,Z\left(  s\right)  +\bar{c}_{2}\left(  s,e\right)  \right)  -\Theta\left(
s,X\left(  s\right)  ,Z\left(  s\right)  \right)  \right\vert ^{2}%
\vartheta\left(  de\right)  ds\right]  <\infty\text{.} \label{Eq08}%
\end{equation}
where%
\begin{align}
\bar{c}_{1}\left(  s,e\right)   &  :=\bar{c}\left(  s,X_{-}\left(  s\right)
,Z_{-}\left(  s\right)  ,e\right)  \text{,}\label{Eq04}\\
\bar{c}_{2}\left(  s,e\right)   &  :=\bar{c}\left(  s,Z_{-}\left(  s\right)
,Z_{-}\left(  s\right)  ,e\right)  \text{.} \label{Eq05}%
\end{align}
Then for a.e. $\left(  s,e\right)  \in\left[  t_{0},T\right]  \times E$, a.s.,%
\begin{align*}
Y\left(  s\right)   &  =\Theta\left(  s,X\left(  s\right)  ,Z\left(  s\right)
\right)  \text{,}\\
Q\left(  s\right)   &  =\bar{\sigma}\left(  s,X\left(  s\right)  ,Z\left(
s\right)  \right)  ^{\top}\Theta_{x}\left(  s,X\left(  s\right)  ,Z\left(
s\right)  \right) \\
&  +\bar{\sigma}\left(  s,Z\left(  s\right)  ,Z\left(  s\right)  \right)
^{\top}\Theta_{z}\left(  s,X\left(  s\right)  ,Z\left(  s\right)  \right)
\text{,}\\
R\left(  s,e\right)   &  =\Theta\left(  s,X_{-}\left(  s\right)  +\bar{c}%
_{1}\left(  s,e\right)  ,Z_{-}\left(  s\right)  +\bar{c}_{2}\left(
s,e\right)  \right)  -\Theta\left(  s,X\left(  s\right)  ,Z\left(  s\right)
\right)  \text{.}%
\end{align*}

\end{theorem}

Before presenting a proof of the above theorem we first recall a version of
It\^{o}'s formula to jump-diffusion processes; see e.g. [\cite{Oksendal},
Theorem 1.16].

\begin{lemma}
\label{result001}Let $X_{1}\left(  \cdot\right)  $, $X_{2}\left(
\cdot\right)  \in\mathcal{S}_{\mathcal{F}}^{2}\left(  t_{0},T;\mathbb{R}%
^{n}\right)  $ be two processes of the form%
\[
dX_{i}\left(  s\right)  =b_{i}\left(  s\right)  ds+\sigma_{i}\left(  s\right)
dW\left(  s\right)  +%
{\displaystyle\int_{E}}
c_{i}\left(  s,e\right)  \tilde{N}\left(  ds,de\right)  ,\text{ for }i=1,2,
\]
where $b_{i}\left(  \cdot\right)  \in\mathcal{L}_{\mathcal{F}}^{1}\left(
t_{0},T;\mathbb{R}^{n}\right)  ,$ $\sigma_{i}\left(  \cdot\right)
\in\mathcal{L}_{\mathcal{F}}^{2}\left(  t_{0},T;\mathbb{R}^{n\times d}\right)
$ and $c_{i}\left(  \cdot,\cdot\right)  \in\mathcal{L}_{%
\mathcal{F}%
,p}^{\vartheta,q}\left(  \left[  t_{0},T\right]  \times E;%
\mathbb{R}
^{n}\right)  $. Suppose that $V\left(  \cdot,\cdot,\cdot\right)
\in\mathcal{C}^{1,2,2}\left(  \left[  t_{0},T\right]  \times%
\mathbb{R}
^{n}\times%
\mathbb{R}
^{n};%
\mathbb{R}
\right)  $. Then for all $s\in\left[  t_{0},T\right]  $, a.s.,
\begin{align*}
&  dV\left(  s,X_{1}\left(  s\right)  ,X_{2}\left(  s\right)  \right) \\
&  =\left\{  V_{s}\left(  s,X_{1}\left(  s\right)  ,X_{2}\left(  s\right)
\right)  +\left\langle V_{x}\left(  s,X_{1}\left(  s\right)  ,X_{2}\left(
s\right)  \right)  ,b_{1}\left(  s\right)  \right\rangle \right.  \ \\
&  +\left\langle V_{z}\left(  s,X_{1}\left(  s\right)  ,X_{2}\left(  s\right)
\right)  ,b_{2}\left(  s\right)  \right\rangle +\frac{1}{2}\text{\textbf{tr}%
}\left[  \sigma_{1}\left(  s\right)  \sigma_{1}\left(  s\right)  ^{\top}%
V_{xx}\left(  s,X_{1}\left(  s\right)  ,X_{2}\left(  s\right)  \right)
\right] \\
&  +\frac{1}{2}\text{\textbf{tr}}\left[  \sigma_{2}\left(  s\right)
\sigma_{2}\left(  s\right)  ^{\top}V_{zz}\left(  s,X_{1}\left(  s\right)
,X_{2}\left(  s\right)  \right)  \right]  +\text{\textbf{tr}}\left[
\sigma_{1}\left(  s\right)  ^{\top}V_{xz}\left(  s,X_{1}\left(  s\right)
,X_{2}\left(  s\right)  \right)  \sigma_{2}\left(  s\right)  \right]
\end{align*}%
\begin{align*}
&  +%
{\displaystyle\int_{E}}
\left\{  V\left(  s,X_{1}\left(  s\right)  +c_{1}\left(  s,e\right)
,X_{2}\left(  s\right)  +c_{2}\left(  s,e\right)  \right)  -V\left(
s,X_{1}\left(  s\right)  ,X_{2}\left(  s\right)  \right)  \right. \\
&  \left.  \left.  -\left\langle V_{x}\left(  s,X_{1}\left(  s\right)
,X_{2}\left(  s\right)  \right)  ,c_{1}\left(  s,e\right)  \right\rangle
-\left\langle V_{z}\left(  s,X_{1}\left(  s\right)  ,X_{2}\left(  s\right)
\right)  ,c_{2}\left(  s,e\right)  \right\rangle \right\}  \vartheta\left(
de\right)  \right\}  ds\\
&  +\left\{  V_{x}\left(  s,X_{1}\left(  s\right)  ,X_{2}\left(  s\right)
\right)  ^{\top}\sigma_{1}\left(  s\right)  +V_{z}\left(  s,X_{1}\left(
s\right)  ,X_{2}\left(  s\right)  \right)  ^{\top}\sigma_{2}\left(  s\right)
\right\}  dW\left(  s\right) \\
&  +%
{\displaystyle\int\limits_{E}}
\left\{  V\left(  s,X_{1}\left(  s\right)  +c_{1}\left(  s,e\right)
,X_{2}\left(  s\right)  +c_{2}\left(  s,e\right)  \right)  -V\left(
s,X_{1}\left(  s\right)  ,X_{2}\left(  s\right)  \right)  \right\}  \tilde
{N}\left(  ds,de\right)  \text{.}%
\end{align*}

\end{lemma}

\textbf{Proof of Theorem \ref{result00}}. First define for $s\in\left[
t_{0},T\right]  $,%
\begin{align*}
\bar{Y}\left(  s\right)   &  :=\Theta\left(  s,X\left(  s\right)  ,Z\left(
s\right)  \right)  \text{,}\\
\bar{Q}\left(  s\right)   &  :=\bar{\sigma}\left(  s,X\left(  s\right)
,Z\left(  s\right)  \right)  ^{\top}\Theta_{x}\left(  s,X\left(  s\right)
,Z\left(  s\right)  \right) \\
&  +\bar{\sigma}\left(  s,Z\left(  s\right)  ,Z\left(  s\right)  \right)
^{\top}\Theta_{z}\left(  s,X\left(  s\right)  ,Z\left(  s\right)  \right)
\text{,}%
\end{align*}
and%
\[
\bar{R}\left(  s,e\right)  :=\Theta\left(  s,X_{-}\left(  s\right)  +\bar
{c}_{1}\left(  s,e\right)  ,Z_{-}\left(  s\right)  +\bar{c}_{2}\left(
s,e\right)  \right)  -\Theta\left(  s,X\left(  s\right)  ,Z\left(  s\right)
\right)  \text{,}%
\]
where $X\left(  \cdot\right)  $, $Z\left(  \cdot\right)  $ are the
corresponding solutions to the forward equations in (\ref{Eq01}) and
$c_{1}\left(  s,e\right)  $, $c_{2}\left(  s,e\right)  $ are as introduced in
(\ref{Eq04})-(\ref{Eq05}). Applying It\^{o}'s formula to $\Theta\left(
\cdot,X\left(  \cdot\right)  ,Z\left(  \cdot\right)  \right)  $, we get%
\begin{align}
&  d\bar{Y}\left(  s\right) \nonumber\\
&  =\left\{  \Theta_{s}\left(  s,X\left(  s\right)  ,Z\left(  s\right)
\right)  +\left\langle \Theta_{x}\left(  s,X\left(  s\right)  ,Z\left(
s\right)  \right)  ,\bar{\mu}\left(  s,X\left(  s\right)  ,Z\left(  s\right)
\right)  \right\rangle \right.  \ \nonumber\\
&  +\left\langle \Theta_{z}\left(  s,X\left(  s\right)  ,Z\left(  s\right)
\right)  ,\bar{\mu}\left(  s,Z\left(  s\right)  ,Z\left(  s\right)  \right)
\right\rangle +\frac{1}{2}\text{\textbf{tr}}\left[  \bar{\sigma}\bar{\sigma
}^{\top}\left(  s,X\left(  s\right)  ,Z\left(  s\right)  \right)  \Theta
_{xx}\left(  s,X\left(  s\right)  ,Z\left(  s\right)  \right)  \right]
\nonumber\\
&  +\frac{1}{2}\text{\textbf{tr}}\left[  \bar{\sigma}\bar{\sigma}^{\top
}\left(  s,Z\left(  s\right)  ,Z\left(  s\right)  \right)  \Theta_{zz}\left(
s,X\left(  s\right)  ,Z\left(  s\right)  \right)  \right] \nonumber\\
&  +\text{\textbf{tr}}\left[  \bar{\sigma}\left(  s,X\left(  s\right)
,Z\left(  s\right)  \right)  ^{\top}\Theta_{xz}\left(  s,X\left(  s\right)
,Z\left(  s\right)  \right)  \bar{\sigma}\left(  s,Z\left(  s\right)
,Z\left(  s\right)  \right)  \right] \nonumber\\
&  +%
{\displaystyle\int\limits_{E}}
\left\{  \Theta\left(  s,X\left(  s\right)  +\bar{c}_{1}\left(  s,e\right)
,Z\left(  s\right)  +\bar{c}_{2}\left(  s,e\right)  \right)  -\Theta\left(
s,X\left(  s\right)  ,Z\left(  s\right)  \right)  \right. \nonumber\\
&  \left.  \left.  -\left\langle \Theta_{x}\left(  s,X\left(  s\right)
,Z\left(  s\right)  \right)  ,\bar{c}_{1}\left(  s,e\right)  \right\rangle
-\left\langle \Theta_{z}\left(  s,X\left(  s\right)  ,Z\left(  s\right)
\right)  ,\bar{c}_{2}\left(  s,e\right)  \right\rangle \right\}
\vartheta\left(  de\right)  \right\}  ds\nonumber\\
&  +\left\{  \Theta_{x}\left(  s,X\left(  s\right)  ,Z\left(  s\right)
\right)  ^{\top}\bar{\sigma}\left(  s,X\left(  s\right)  ,Z\left(  s\right)
\right)  +\Theta_{z}\left(  s,X\left(  s\right)  ,Z\left(  s\right)  \right)
^{\top}\bar{\sigma}\left(  s,Z\left(  s\right)  ,Z\left(  s\right)  \right)
\right\}  dW\left(  s\right) \nonumber\\
&  +%
{\displaystyle\int\limits_{E}}
\left\{  \Theta\left(  s,X\left(  s\right)  +\bar{c}_{1}\left(  s,e\right)
,Z\left(  s\right)  +\bar{c}_{2}\left(  s,e\right)  \right)  -\Theta\left(
s,X\left(  s\right)  ,Z\left(  s\right)  \right)  \right\}  \tilde{N}\left(
ds,de\right)  \text{.} \label{Eq03}%
\end{align}

On the other hand, it follows from the IPDE (\ref{Eq02}) that%
\begin{align*}
&  \Theta_{s}\left(  s,X\left(  s\right)  ,Z\left(  s\right)  \right) \\
&  =-\text{\ }\bar{f}\left(  s,X\left(  s\right)  ,Z\left(  s\right)  \right)
-\left\langle \Theta_{x}\left(  s,X\left(  s\right)  ,Z\left(  s\right)
\right)  ,\bar{\mu}\left(  s,X\left(  s\right)  ,Z\left(  s\right)  \right)
\right\rangle \\
&  -\left\langle \Theta_{z}\left(  s,X\left(  s\right)  ,Z\left(  s\right)
\right)  ,\bar{\mu}\left(  s,Z\left(  s\right)  ,Z\left(  s\right)  \right)
\right\rangle -\frac{1}{2}\text{\textbf{tr}}\left[  \bar{\sigma}\bar{\sigma
}^{\top}\left(  s,X\left(  s\right)  ,Z\left(  s\right)  \right)  \Theta
_{xx}\left(  s,X\left(  s\right)  ,Z\left(  s\right)  \right)  \right] \\
&  -\frac{1}{2}\text{\textbf{tr}}\left[  \bar{\sigma}\bar{\sigma}^{\top
}\left(  s,Z\left(  s\right)  ,Z\left(  s\right)  \right)  \Theta_{zz}\left(
s,X\left(  s\right)  ,Z\left(  s\right)  \right)  \right] \\
&  -\text{\textbf{tr}}\left[  \bar{\sigma}\left(  s,X\left(  s\right)
,Z\left(  s\right)  \right)  ^{\top}\Theta_{xz}\left(  s,X\left(  s\right)
,Z\left(  s\right)  \right)  \bar{\sigma}\left(  s,Z\left(  s\right)
,Z\left(  s\right)  \right)  \right] \\
&  -%
{\displaystyle\int\limits_{E}}
\left\{  \Theta\left(  s,X\left(  s\right)  +\bar{c}_{1}\left(  s,e\right)
,Z\left(  s\right)  +\bar{c}_{2}\left(  s,e\right)  \right)  -\Theta\left(
s,X\left(  s\right)  ,Z\left(  s\right)  \right)  \right. \\
&  \left.  -\left\langle \Theta_{x}\left(  s,X\left(  s\right)  ,Z\left(
s\right)  \right)  ,\bar{c}_{1}\left(  s,e\right)  \right\rangle -\left\langle
\Theta_{z}\left(  s,X\left(  s\right)  ,Z\left(  s\right)  \right)  ,\bar
{c}_{2}\left(  s,e\right)  \right\rangle \right\}  \vartheta\left(  de\right)
.
\end{align*}

Invoking this into (\ref{Eq03}), we obtain that $\left(  \bar{Y}\left(
\cdot\right)  ,\bar{Q}\left(  \cdot\right)  ,\bar{R}\left(  \cdot
,\cdot\right)  \right)  $ satisfies the following BSDE%
\[
\left\{
\begin{array}
[c]{r}%
d\bar{Y}\left(  s\right)  =-\bar{f}\left(  s,X\left(  s\right)  ,Z\left(
s\right)  \right)  ds+\bar{Q}\left(  s\right)  ^{\top}dW\left(  s\right)
\text{\ }+%
{\displaystyle\int_{E}}
\bar{R}\left(  s,e\right)  \tilde{N}\left(  ds,de\right)  \text{,}\\
\text{for }s\in\left[  t_{0},T\right]  \text{,}\\
\multicolumn{1}{l}{\bar{Y}\left(  T\right)  =\bar{F}\left(  X\left(  T\right)
\right)  \text{.}}%
\end{array}
\right.
\]
Hence, by the uniqueness of the solutions to the BSDE in (\ref{Eq01}), we
obtain that%
\[
\left(  Y\left(  s\right)  ,Q\left(  s\right)  ,R\left(  s,e\right)  \right)
=\left(  \bar{Y}\left(  s\right)  ,\bar{Q}\left(  s\right)  ,\bar{R}\left(
s,e\right)  \right)  ,\text{ a.s., a.e. }\left(  s,e\right)  \in\left[
t_{0},T\right]  \times E\text{.}%
\]
\eop

\section{Formulation of the problem\label{section3}}

Given a subset $U\subset%
\mathbb{R}
^{l}$, let $\mu:\left[  0,T\right]  \times\mathbb{R}^{n}\times U\rightarrow
\mathbb{R}^{n}$, $\sigma:\left[  0,T\right]  \times\mathbb{R}^{n}\times
U\rightarrow\mathbb{R}^{n\times d}$ and $c:\left[  0,T\right]  \times%
\mathbb{R}
^{n}\times U\times E\rightarrow%
\mathbb{R}
^{n}$ be three deterministic measurable functions. Consider on the time
interval $\left[  0,T\right]  $ the following controlled stochastic
differential equation with jumps%
\begin{equation}
\left\{
\begin{array}
[c]{l}%
dX^{x_{0},u\left(  \cdot\right)  }\left(  s\right)  =\mu\left(  s,X^{x_{0}%
,u\left(  \cdot\right)  }\left(  s\right)  ,u\left(  s\right)  \right)
ds+\sigma\left(  s,X^{x_{0},u\left(  \cdot\right)  }\left(  s\right)
,u\left(  s\right)  \right)  dW\left(  s\right) \\
\text{ \ \ \ \ \ \ \ \ \ \ \ \ }+%
{\displaystyle\int_{E}}
c\left(  s,X_{-}^{x_{0},u\left(  \cdot\right)  }\left(  s\right)  ,u\left(
s\right)  ,e\right)  \tilde{N}\left(  ds,de\right)  \text{,\ }s\in\left[
0,T\right]  \text{,}\\
X^{x_{0},u\left(  \cdot\right)  }\left(  0\right)  =x_{0}\text{,}%
\end{array}
\right.  \label{Eq1}%
\end{equation}
where $u:\left[  0,T\right]  \times\Omega\rightarrow U$ represents the control
process, $X^{x_{0},u\left(  \cdot\right)  }\left(  \cdot\right)  $ is the
controlled state process and $x_{0}\in%
\mathbb{R}
^{n}$ is regarded as the initial state.

As time evolves, we need to consider the following controlled stochastic
differential equation starting from the situation $\left(  t,y\right)
\in\left[  0,T\right]  \times%
\mathbb{R}
^{n}$\emph{,}%
\begin{equation}
\left\{
\begin{array}
[c]{l}%
dX\left(  s\right)  =\mu\left(  s,X\left(  s\right)  ,u\left(  s\right)
\right)  ds+\sigma\left(  s,X\left(  s\right)  ,u\left(  s\right)  \right)
dW\left(  s\right) \\
\text{ \ \ \ \ \ \ \ \ \ \ \ \ }+%
{\displaystyle\int_{E}}
c\left(  s,X_{-}\left(  s\right)  ,u\left(  s\right)  ,e\right)  \tilde
{N}\left(  ds,de\right)  \text{,\ }s\in\left[  t,T\right]  \text{,}\\
X\left(  t\right)  =y\text{.}%
\end{array}
\right.  \label{Eq2}%
\end{equation}
where $X\left(  \cdot\right)  =X^{t,y,u\left(  \cdot\right)  }\left(
\cdot\right)  $ denotes its solution. For any initial state $\left(
t,y\right)  \in\left[  0,T\right]  \times%
\mathbb{R}
^{n}$, in order to evaluate the performance of a control process $u\left(
\cdot\right)  $, we introduce the cost functional%

\begin{equation}
\mathbf{J}\left(  t,y;u\left(  \cdot\right)  \right)  :=\mathbb{E}_{t}\left[
\int_{t}^{T}f\left(  t,s,X\left(  s\right)  ,u\left(  s\right)  \right)
ds+F\left(  t,X\left(  T\right)  \right)  \right]  +G\left(  t,y,\mathbb{E}%
_{t}\left[  \Psi\left(  X\left(  T\right)  \right)  \right]  \right)  \text{,}
\label{Eq3}%
\end{equation}
where $\mathbb{E}_{t}\left[  \cdot\right]  =\mathbb{E}_{t}\left[
\cdot|\mathcal{F}_{t}\right]  $; $f:\left[  0,T\right]  \times\left[
0,T\right]  \times%
\mathbb{R}
^{n}\times U\rightarrow%
\mathbb{R}
$, $F:\left[  0,T\right]  \times%
\mathbb{R}
^{n}\rightarrow%
\mathbb{R}
$, $G:\left[  0,T\right]  \times%
\mathbb{R}
^{n}\times%
\mathbb{R}
^{m}\rightarrow%
\mathbb{R}
$ and $\Psi=\left(  \Psi_{1},...,\Psi_{m}\right)  ^{\top}:%
\mathbb{R}
^{n}\rightarrow%
\mathbb{R}
^{m}$ are four deterministic measurable functions.

Of course the appearance of $t$, $y$ and $\mathbb{E}_{t}\left[  \Psi\left(
X\left(  T\right)  \right)  \right]  $ in the coefficients $f$, $F$ and $G$ is
not just means an extension from mathematical point of view, but also be of
great importance in applications. More specifically:

\begin{enumerate}
\item[(i)] The dependence of $f\left(  t,s,X\left(  s\right)  ,u\left(
s\right)  \right)  $ and $F\left(  t,X\left(  T\right)  \right)  $ on the
initial time $t$ is motivated by general discounting situations; see e.g.
\cite{Strotz}, \cite{EkePir}, \cite{Yong2012}, \cite{WangYong} and \cite{Wei1}.

\item[(iii)] In the term $G\left(  t,y,\mathbb{E}_{t}\left[  \Psi\left(
X\left(  T\right)  \right)  \right]  \right)  $ we have a non linear function
$G$ acting on $\left(  y,\mathbb{E}_{t}\left[  \Psi\left(  X\left(  T\right)
\right)  \right]  \right)  $, which can be motivated by the mean-variance
criterion with state dependent risk aversion \cite{BMZ}.
\end{enumerate}

\begin{remark}
(i) Let $Q:\left[  0,T\right]  \rightarrow\mathbb{S}^{n}$, $R:\left[
0,T\right]  \rightarrow\mathbb{S}^{l}$ be two deterministic measurable maps
and $h,$ $\bar{h}\in\mathbb{S}^{n}$, $\mu_{1}\in%
\mathbb{R}
^{n\times n}$, $\mu_{2}\in%
\mathbb{R}
^{n}$ be four matrices. If $f$, $F$, $\Psi$ and $G$ are of the following form%
\begin{align*}
f\left(  t,s,x,u\right)   &  :=\frac{1}{2}\left(  \left\langle Q\left(
s\right)  x,x\right\rangle +\left\langle R\left(  s\right)  u,u\right\rangle
\right)  \text{,}\\
F\left(  t,x\right)   &  :=\frac{1}{2}\left\langle hx,x\right\rangle
,\Psi\left(  x\right)  :=x\text{ and}\\
G\left(  t,y,\bar{x}\right)   &  :=-\left\langle \frac{1}{2}\bar{h}\bar{x}%
+\mu_{1}y+\mu_{2},\bar{x}\right\rangle \text{,}%
\end{align*}
then the cost functional (\ref{Eq3}) becomes the same as (2.3) in
\cite{Huetal2011} (in the case when the coefficients are
deterministic).\newline(ii) In the case when $G\left(  t,y,\bar{x}\right)
\equiv0$, the cost functional (\ref{Eq3}) reduces to the same as (3.2) in
\cite{Yong2012}.\newline
\end{remark}

We introduce the following assumptions.

\begin{enumerate}
\item[\textbf{(H1)}] The maps $\mu$, $\sigma$ and $c$ are continuous and there
exists a constant $C>0$ such that for all $\left(  x,y\right)  \in%
\mathbb{R}
^{n}\times%
\mathbb{R}
^{n}$, $\left(  t,u,e\right)  \in\left[  0,T\right]  \times U\times E$ and
$q\geq2$%
\[%
\begin{array}
[c]{rr}%
\left\vert \mu\left(  s,x,u\right)  -\mu\left(  s,y,u\right)  \right\vert
+\left\vert \sigma\left(  s,x,u\right)  -\sigma\left(  s,y,u\right)
\right\vert  & \\
+\left\vert c\left(  s,x,u,e\right)  -c\left(  s,y,u,e\right)  \right\vert  &
\leq C\left\vert x-y\right\vert
\end{array}
\]
and%
\[
\left\vert \mu\left(  s,0,u\right)  \right\vert +\left\vert \sigma\left(
s,0,u\right)  \right\vert +\left(
{\displaystyle\int_{E}}
\left\vert c\left(  s,0,u,e\right)  \right\vert ^{q}\vartheta\left(
de\right)  \right)  ^{\frac{1}{q}}\leq C\left(  1+\left\vert u\right\vert
\right)  \text{. }%
\]

\item[\textbf{(H2)}] \textbf{(i)} The maps $f$, $h$ are continuous and
quadratic growth on $\left(  x,u\right)  $ uniformly in time, i.e. there
exists a constant $C>0$ such that for all $x\in%
\mathbb{R}
^{n}$ and $\left(  t,s,u\right)  \in\left[  0,T\right]  \times\left[
0,T\right]  \times U$,%
\begin{align*}
\left\vert f\left(  t,s,x,u\right)  \right\vert  &  \leq C\left(  1+\left\vert
x\right\vert ^{2}+\left\vert u\right\vert ^{2}\right)  \text{,}\\
\left\vert F\left(  t,x\right)  \right\vert  &  \leq C\left(  1+\left\vert
x\right\vert ^{2}\right)  \text{.}%
\end{align*}
\textbf{(ii)} The map $G$ is continuously differentiable with respect to
$\bar{x}$ and there exists a constant $C>0$ such that for all $\left(
t,y,\bar{x}\right)  \in\left[  0,T\right]  \times%
\mathbb{R}
^{n}\times%
\mathbb{R}
^{m}$,%
\begin{align*}
\left\vert G\left(  t,y,\bar{x}\right)  \right\vert  &  \leq C\left(
1+\left\vert y\right\vert ^{2}+\left\vert \bar{x}\right\vert ^{2}\right)
\text{,}\\
\left\vert G_{\bar{x}}\left(  t,y,\bar{x}\right)  \right\vert  &  \leq
C\left(  1+\left\vert y\right\vert +\left\vert \bar{x}\right\vert \right)
\text{.}%
\end{align*}
\textbf{(iii)} The map $\Psi$ is continuous and linear growth on $x$, i.e.
there exists a constant $C>0$ such that for all $x\in%
\mathbb{R}
^{n}$,%
\[
\left\vert \Psi\left(  x\right)  \right\vert \leq C\left(  1+\left\vert
x\right\vert \right)  \text{.}%
\]

\end{enumerate}

Under Assumptions \textbf{(H1)-(H2)}, for any initial pair $\left(
t,y\right)  $ $\in$ $\left[  0,T\right]  \times\mathbb{R}^{n}$ and a control
$u\left(  \cdot\right)  \in\mathcal{L}_{%
\mathcal{F}%
,p}^{q}\left(  t,T;%
\mathbb{R}
^{l}\right)  $, with $q\geq2$, the state equation (\ref{Eq1}) admits a unique
solution $X\left(  \cdot\right)  =X^{t,y,u(\text{.})}\left(  \cdot\right)
\in\mathcal{S}_{\mathcal{F}}^{q}\left(  t,T;\mathbb{R}^{n}\right)  $ (see e.g.
\cite{Romuald}, pp. 125-126) and the cost functional (\ref{Eq3}) is
well-defined. Moreover, there exists a constant $C>0$ such that%
\[
\mathbb{E}\left[  \sup_{t\leq s\leq T}\left\vert X\left(  s\right)
\right\vert ^{q}\right]  \leq C\left(  1+\left\vert y\right\vert ^{q}\right)
\text{.}%
\]

We now introduce the class of admissible controls.

\begin{definition}
[Admissible control]An admissible control $u\left(  \cdot\right)  $ over
$\left[  t,T\right]  $ is a $U$-valued $\left(  \mathcal{F}_{s}\right)
_{s\in\left[  t,T\right]  }$-predictable process such that%
\[
\mathbb{E}\left[  \int_{s}^{T}\left\vert u\left(  s\right)  \right\vert
^{4}ds\right]  <\infty\text{.}%
\]
The class of admissible controls over $\left[  t,T\right]  $ is denoted by
$\mathcal{U}\left[  t,T\right]  $.
\end{definition}

We also introduce the class of admissible closed-loop controls.

\begin{definition}
[Admissible Strategy]\label{Def1}A map $\varphi:\left[  0,T\right]  \times%
\mathbb{R}
^{n}\rightarrow U$ is called an admissible strategy if for every $\left(
t,y\right)  \in\left[  0,T\right]  \times%
\mathbb{R}
^{n}$ the following SDE%
\begin{equation}
\left\{
\begin{array}
[c]{l}%
dX\left(  s\right)  =\mu\left(  s,X\left(  s\right)  ,\varphi\left(
s,X\left(  s\right)  \right)  \right)  ds+\sigma\left(  s,X\left(  s\right)
,\varphi\left(  s,X\left(  s\right)  \right)  \right)  dW\left(  s\right) \\
\text{ \ \ \ \ \ \ \ \ \ \ }+\ \int_{E}c\left(  s,X_{-}\left(  s\right)
,\varphi\left(  s,X_{-}\left(  s\right)  \right)  ,e\right)  \tilde{N}\left(
ds,de\right)  \text{,\ }s\in\left[  t,T\right]  \text{,}\\
X\left(  t\right)  =y
\end{array}
\right.  \label{Eq4}%
\end{equation}
admits a unique strong solution $X\left(  \cdot\right)  =X^{t,y,\varphi\left(
\cdot,\cdot\right)  }\left(  \cdot\right)  \in\mathcal{S}_{\mathcal{F}}%
^{4}\left(  t,T;\mathbb{R}^{n}\right)  $ and the control process
$u^{t,y,\varphi}\left(  \cdot\right)  =\varphi\left(  \cdot,X_{-}%
^{t,y,\varphi}\left(  \cdot\right)  \right)  \in\mathcal{U}\left[  t,T\right]
$. The class of admissible strategies is denoted by $\mathcal{S}$.
\end{definition}

Subsequently, we will use the notation $\varphi^{s}\left(  z\right)  $ instead
of $\varphi\left(  s,z\right)  $, whenever no confusion arises. Our stochastic
optimal control problem can be stated as follows.\newline\textbf{Problem (N).
}\textit{For any given initial pair} $\left(  t,y\right)  \in\left[
0,T\right]  \times%
\mathbb{R}
^{n}$, \textit{find a} $\bar{u}^{t,y}\left(  \cdot\right)  \in\mathcal{U}%
\left[  t,T\right]  $\textit{ such that}%
\begin{equation}
\mathbf{J}\left(  t,y;\bar{u}^{t,y}\left(  \cdot\right)  \right)
=\inf_{u\left(  \cdot\right)  \in\mathcal{U}\left[  t,T\right]  }%
\mathbf{J}\left(  t,y;u\left(  \cdot\right)  \right)  \text{.} \label{Eq5}%
\end{equation}

For a given $\left(  t,y\right)  \in\left[  0,T\right]  \times%
\mathbb{R}
^{n}$, any $\bar{u}^{t,y}\left(  \cdot\right)  \in$ $\mathcal{U}\left[
t,T\right]  $ satisfying the above is called \textit{a pre-commitment optimal}
control for Problem (N) at $\left(  t,y\right)  $.

As well documented in the review paper \cite{Yan}, the appearance of the terms
$t$, $y$ and $\mathbb{E}_{t}\left[  \Psi\left(  X\left(  T\right)  \right)
\right]  $ in the performance functional (\ref{Eq3}) destroys the
time-consistency of pre-committed optimal controls of Problem \textbf{(N)}: If
we find some $\bar{u}^{t,y}\left(  \cdot\right)  $ satisfying (\ref{Eq5}) for
some initial pair $\left(  t,y\right)  $, we do not have%
\[
\mathbf{J}\left(  s,\bar{X}^{t,y}\left(  s\right)  ;\bar{u}|_{\left[
s,T\right]  }^{t,y}\left(  \cdot\right)  \right)  =\inf_{u\left(
\cdot\right)  \in\mathcal{U}\left[  s,T\right]  }\mathbf{J}\left(  s,\bar
{X}^{t,y}\left(  s\right)  ;u\left(  \cdot\right)  \right)  \text{,}%
\]
in general, for any $s\in\left(  t,T\right)  $, where $\bar{X}^{t,y}\left(
\cdot\right)  $ is the corresponding solution of (\ref{Eq2}) with $u\left(
\cdot\right)  =\bar{u}^{t,y}\left(  \cdot\right)  $ and $\bar{u}|_{\left[
s,T\right]  }^{t,y}\left(  \cdot\right)  $ denotes the restriction of $\bar
{u}^{t,y}\left(  \cdot\right)  $ on time interval $\left[  s,T\right]  $.

Since time-consistency is important for a rational controller, the concept
\textquotedblleft optimality\textquotedblright\ needs to be reconsidered in a
more sophisticated way. Here we adopt the concept of equilibrium strategy
which is, at any $t\in\left[  0,T\right]  $, optimal \textquotedblleft
infinitesimally\textquotedblright\ via spike variation.

Let $\hat{\varphi}\left(  \cdot,\cdot\right)  \in\mathcal{S}$ be a given
admissible strategy and $\hat{X}^{x_{0}}\left(  \cdot\right)  =X^{x_{0}%
,\hat{\varphi}\left(  \cdot,\cdot\right)  }\left(  \cdot\right)  $ be the
corresponding solution of (\ref{Eq4}) with $\varphi\left(  \cdot,\cdot\right)
=\hat{\varphi}\left(  \cdot,\cdot\right)  $ and $\left(  t,y\right)  =\left(
0,x_{0}\right)  $. For any $t\in\left[  0,T\right]  $, $v\in\mathbb{L}%
^{4}\left(  \Omega,\mathcal{F}_{t},\mathbb{P};U\right)  $ and for any
$\varepsilon\in\left[  0,T-t\right)  $ define%
\begin{equation}
u^{\varepsilon}\left(  s\right)  =\left\{
\begin{array}
[c]{l}%
v\text{, for }s\in\left[  t,t+\varepsilon\right]  \text{,}\\
\hat{\varphi}^{s}\left(  \hat{X}_{-}^{x_{0}}\left(  s\right)  \right)  \text{,
for }s\in\left(  t+\varepsilon,T\right)  .
\end{array}
\right.  \label{Eq6}%
\end{equation}

It is easy to verify that $u^{\varepsilon}\left(  \cdot\right)  \in
\mathcal{U}\left[  t,T\right]  $. The following definition is slightly
different from the original definition of open-loop equilibrium controls
defined by Hu et al. (\cite{Huetal2011}, \cite{Huetal2017}) in a
time-inconsistent LQ control problem.

\begin{definition}
[\textbf{Open-Loop Equilibrium Strategy}]\label{Def2}A strategy $\hat{\varphi
}\left(  \cdot,\cdot\right)  \in\mathcal{S}$ is called \textbf{an open-loop
equilibrium strategy}\textit{ for Problem (N) }if%
\begin{equation}
\underset{\varepsilon\downarrow0}{\lim\inf}\frac{1}{\varepsilon}\left\{
\mathbf{J}\left(  t,\hat{X}^{x_{0}}\left(  t\right)  ;u^{\varepsilon}\left(
\cdot\right)  \right)  -\mathbf{J}\left(  t,\hat{X}^{x_{0}}\left(  t\right)
;\hat{\varphi}\left(  \cdot,\cdot\right)  \right)  \right\}  \geq0\text{,}
\label{Eq7}%
\end{equation}
\textit{where }$u^{\varepsilon}\left(  \cdot\right)  $\textit{ is defined by
}(\ref{Eq6}), for any $t\in\left[  0,T\right]  $, $v\in\mathbb{L}^{4}\left(
\Omega,\mathcal{F}_{t},\mathbb{P};U\right)  $ and%
\begin{align*}
\mathbf{J}\left(  t,\hat{X}^{x_{0}}\left(  t\right)  ;\hat{\varphi}\left(
\cdot,\cdot\right)  \right)   &  :=\mathbb{E}_{t}\left[  \int_{t}^{T}f\left(
t,s,\hat{X}^{x_{0}}\left(  s\right)  ,\hat{\varphi}^{s}\left(  \hat{X}^{x_{0}%
}\left(  s\right)  \right)  \right)  ds+F\left(  t,\hat{X}^{x_{0}}\left(
T\right)  \right)  \right] \\
&  +G\left(  t,\hat{X}^{x_{0}}\left(  t\right)  ,\mathbb{E}_{t}\left[
\Psi\left(  \hat{X}^{x_{0}}\left(  T\right)  \right)  \right]  \right)
\text{.}%
\end{align*}

\end{definition}

\begin{remark}
The intuition behind this definition is similar to that of \cite{Bjork4}: For
any $t\in\left[  0,T\right]  $, on the premise that for every $s>t$, the
optimal decision for the controller at $s$ is $\hat{\varphi}^{s}\left(
\hat{X}_{-}^{x_{0}}\left(  s\right)  \right)  $, then the optimal choice for
the controller at $t$ is that he/she also uses the decision $\hat{\varphi}%
^{t}\left(  \hat{X}_{-}^{x_{0}}\left(  t\right)  \right)  $. However, even
thought the above-described equilibrium concept concerns closed-loop controls
only, there is a fundamental difference between the definition here and the
ones in \cite{Bjork4} and \cite{Yong2012}; indeed, in the above definition,
the perturbation of $\hat{\varphi}\left(  \cdot,\cdot\right)  $ in $\left[
t,t+\varepsilon\right]  $ will not affect the control $\hat{\varphi}%
^{s}\left(  \hat{X}_{-}^{x_{0}}\left(  s\right)  \right)  $ in $\left(
t+\varepsilon,T\right)  $, which is not the case with closed-loop equilibrium
concepts considered in \cite{Bjork4} and \cite{Yong2012}.
\end{remark}

\begin{remark}
In the work of Yan and Yong \cite{Yan} the definition of an\textbf{ open-loop
equilibrium control} $\hat{u}\left(  \cdot\right)  \in\mathcal{U}\left[
0,T\right]  $ is given by%
\begin{equation}
\underset{\varepsilon\downarrow0}{\lim\inf}\frac{1}{\varepsilon}\left\{
\mathbf{J}\left(  t,\hat{X}\left(  t\right)  ;u^{t,\varepsilon,u}\left(
\cdot\right)  \right)  -\mathbf{J}\left(  t,\hat{X}\left(  t\right)  ;\hat
{u}\left(  \cdot\right)  \right)  \right\}  \geq0\text{,} \label{Eq8**}%
\end{equation}
where $u^{t,\varepsilon,u}\left(  \cdot\right)  :=u\chi_{\left[
t,t+\varepsilon\right]  }\left(  \cdot\right)  +\hat{u}\left(  \cdot\right)
\chi_{\left(  t+\varepsilon,T\right)  }\left(  \cdot\right)  $ for any
$t\in\left[  0,T\right]  $, $u\in U$; $\hat{X}\left(  \cdot\right)  $ is the
corresponding solution of (\ref{Eq1}) with $u\left(  \cdot\right)  =\hat
{u}\left(  \cdot\right)  $. It is easy to see that if $\hat{\varphi}\left(
\cdot,\cdot\right)  \in\mathcal{S}$ is an open-loop equilibrium strategy then
the control process $\hat{u}^{x_{0}}\left(  \cdot\right)  $ defined by%
\[
\hat{u}^{x_{0}}\left(  s\right)  :=\hat{\varphi}^{s}\left(  \hat{X}_{-}%
^{x_{0}}\left(  s\right)  \right)  \text{, for }s\in\left[  0,T\right]
\text{,}%
\]
is an open-loop equilibrium control. Note that the open-loop equilibrium
strategy $\hat{\varphi}\left(  \cdot,\cdot\right)  $ does not depend on the
initial state $x_{0}$, while the open-loop equilibrium control $\hat{u}%
^{x_{0}}\left(  \cdot\right)  $ does. It is also worth mentioning that
$\hat{\varphi}\left(  \cdot,\cdot\right)  $ is a generalization of the notion
"closed-loop representation to an open-loop equilibrium control" introduced by
Wang [\cite{Wang}, \cite{Wang12}] in the framework of time-inconsistent LQ models.
\end{remark}

We close this section with the following definition of a function space that
will be used below.

\begin{definition}
\label{Def3}Let $\hat{\varphi}\left(  \cdot,\cdot\right)  \in\mathcal{S}$ be a
given strategy. A function $h\left(  \cdot,\cdot,\cdot\right)  \in
\mathcal{C}^{1,2,2}\left(  \left[  t,T\right]  \times%
\mathbb{R}
^{n}\times%
\mathbb{R}
^{n};%
\mathbb{R}
\right)  $ is said to belong to the space $\mathcal{C}_{2}^{1,2,2}\left(
\left[  t,T\right]  \times%
\mathbb{R}
^{n}\times%
\mathbb{R}
^{n};%
\mathbb{R}
\right)  $ if it satisfies the condition: There exists a constant $C>0$ such
that for all $\left(  s,u,x,z\right)  \in\left[  t,T\right]  \times
U\times\mathbb{R}^{n}\times\mathbb{R}^{n}$ we have%
\begin{align*}
&  C\left(  1+\left\vert x\right\vert ^{2}+\left\vert u\right\vert
^{2}+\left\vert \hat{\varphi}^{s}\left(  z\right)  \right\vert ^{2}+\left\vert
z\right\vert ^{2}\right) \\
&  \geq\left\vert h\left(  s,x,z\right)  \right\vert +\left\vert h_{x}\left(
s,x,z\right)  \right\vert \left\vert \mu\left(  s,x,u\right)  \right\vert
+\left\vert h_{x}\left(  s,x,z\right)  \right\vert \left\vert \sigma\left(
s,x,u\right)  \right\vert \\
&  +\left\vert h_{z}\left(  s,x,z\right)  \right\vert \left\vert \sigma\left(
s,z,\hat{\varphi}^{s}\left(  z\right)  \right)  \right\vert +\left\vert
h_{xx}\left(  s,x,z\right)  \right\vert \left\vert \sigma\left(  s,x,u\right)
\right\vert ^{2}\\
&  +\left\vert h_{xz}\left(  s,x,z\right)  \right\vert \left\vert
\sigma\left(  s,x,u\right)  \right\vert \left\vert \sigma\left(
s,z,\hat{\varphi}^{s}\left(  z\right)  \right)  \right\vert \\
&  +%
{\displaystyle\int_{E}}
\left\vert h_{x}\left(  s,x,z\right)  \right\vert \left\vert c\left(
\tau,x,u,e\right)  \right\vert \vartheta\left(  de\right)  \text{.}%
\end{align*}

\end{definition}

\section{PDEs and Verification Theorem\label{section4}}

A large number of existing studies in the literature have examined open-loop
equilibrium controls by adopting FBSDEs-approaches. For instance, in the
series of papers \cite{Huetal2011}, \cite{Djehiche}, \cite{Hamaguchi2},
\cite{Wang}, \cite{Wang11}, \cite{Wang12}, \cite{Huetal2019}, \cite{Sun} and
\cite{Zhang} the open-loop equilibrium control process $\hat{u}^{x_{0}}\left(
\cdot\right)  $ is constructed by solving a flow of FBSDEs. Recently, the work
\cite{Alia2} investigated open-loop equilibrium controls by using a
BSPDEs-approach. Therein, an open-loop equilibrium control can be constructed
by solving a system of forward-backward stochastic partial differential
equations (FBSPDEs, for short).

Unlike to the above mentioned papers, our aim in the present article is to
characterize open-loop equilibrium controls by using a PDEs-approach. The main
idea is to introduce two deterministic functions $\hat{\theta}^{t}\left(
s,x,z\right)  $ and $\hat{g}\left(  s,x,z\right)  $ as solutions to two
systems of parabolic integro-partial differential equations. A verification
theorem for open-loop equilibrium strategies can be derived in a natural way
via these functions.

As in \cite{Alia2}, before providing the precise statement of the main results
in this paper, let us begin by establishing some heuristic derivations. Let
$\hat{\varphi}\left(  \cdot,\cdot\right)  \in\mathcal{S}$ be a given
admissible strategy and $\hat{X}^{x_{0}}\left(  \cdot\right)  $ be the unique
strong solution of the SDE%
\begin{equation}
\left\{
\begin{array}
[c]{l}%
d\hat{X}^{x_{0}}\left(  s\right)  =\mu\left(  s,\hat{X}^{x_{0}}\left(
s\right)  ,\hat{\varphi}^{s}\left(  \hat{X}^{x_{0}}\left(  s\right)  \right)
\right)  ds+\sigma\left(  s,\hat{X}^{x_{0}}\left(  s\right)  ,\hat{\varphi
}^{s}\left(  \hat{X}^{x_{0}}\left(  s\right)  \right)  \right)  dW\left(
s\right) \\
\text{ \ \ \ \ \ \ \ \ \ \ \ \ \ \ }+%
{\displaystyle\int_{E}}
c\left(  s,\hat{X}_{-}^{x_{0}}\left(  s\right)  ,\hat{\varphi}^{s}\left(
\hat{X}_{-}^{x_{0}}\left(  s\right)  \right)  ,e\right)  \tilde{N}\left(
ds,de\right)  \text{, for }s\in\left[  0,T\right]  \text{,}\\
\hat{X}^{x_{0}}\left(  0\right)  =x_{0}\text{.}%
\end{array}
\right.  \label{Eq8*}%
\end{equation}
Consider the perturbed strategy $u^{\varepsilon}\left(  \cdot\right)  $
defined by the spike variation (\ref{Eq6}) for some fixed arbitrarily
$t\in\left[  0,T\right]  ,$ $v\in\mathbb{L}^{4}\left(  \Omega,\mathcal{F}%
_{t},\mathbb{P};U\right)  $ and $\varepsilon\in\left[  0,T-t\right)  $. We
would like to determine a suitable expression to the difference%
\[
\mathbf{\Delta\hat{J}}^{\varepsilon}\left(  t\right)  \mathbf{:=J}\left(
t,\hat{X}^{x_{0}}\left(  t\right)  ;u^{\varepsilon}\left(  \cdot\right)
\right)  -\mathbf{J}\left(  t,\hat{X}^{x_{0}}\left(  t\right)  ;\hat{\varphi
}\left(  \cdot,\cdot\right)  \right)  \text{,}%
\]
in order to be able to evaluate the limit in (\ref{Eq7}). To this end, let
$X^{\varepsilon}\left(  \cdot\right)  $ be the corresponding solution of
$u^{\varepsilon}\left(  \cdot\right)  $, i.e.%
\begin{equation}
\left\{
\begin{array}
[c]{l}%
dX^{\varepsilon}\left(  s\right)  =\mu\left(  s,X^{\varepsilon}\left(
s\right)  ,v\right)  ds+\sigma\left(  s,X^{\varepsilon}\left(  s\right)
,v\right)  dW\left(  s\right) \\
\text{ \ \ \ \ \ \ \ \ \ \ \ \ \ \ \ }+%
{\displaystyle\int_{E}}
c\left(  s,X_{-}^{\varepsilon}\left(  s\right)  ,v,e\right)  \tilde{N}\left(
ds,de\right)  ,\text{ for }s\in\left[  t,t+\varepsilon\right]  \text{,}\\
dX^{\varepsilon}\left(  s\right)  =\mu\left(  s,X^{\varepsilon}\left(
s\right)  ,\hat{\varphi}^{s}\left(  \hat{X}^{x_{0}}\left(  s\right)  \right)
\right)  ds+\sigma\left(  s,X^{\varepsilon}\left(  s\right)  ,\hat{\varphi
}^{s}\left(  \hat{X}^{x_{0}}\left(  s\right)  \right)  \right)  dW\left(
s\right) \\
\text{ \ \ \ \ \ \ \ \ \ \ \ \ \ \ \ }+%
{\displaystyle\int_{E}}
c\left(  s,X_{-}^{\varepsilon}\left(  s\right)  ,\hat{\varphi}^{s}\left(
\hat{X}_{-}^{x_{0}}\left(  s\right)  \right)  ,e\right)  \tilde{N}\left(
ds,de\right)  ,\text{ for }s\in\left(  t+\varepsilon,T\right]  \text{,}\\
X^{\varepsilon}\left(  t\right)  =\hat{X}^{x_{0}}\left(  t\right)
\end{array}
\right.  \label{Eq8}%
\end{equation}
and for each fixed $\left(  s,\mathbf{x},\mathbf{\mathbf{z}}\right)
\in\left[  t,T\right]  \times\mathbb{L}^{4}\left(  \Omega,\mathcal{F}%
_{s},\mathbb{P};%
\mathbb{R}
^{n}\right)  \times\mathbb{L}^{4}\left(  \Omega,\mathcal{F}_{s},\mathbb{P};%
\mathbb{R}
^{n}\right)  $, denote by $\hat{X}^{s,\mathbf{z}}\left(  \cdot\right)  $ and
$X^{s,\mathbf{x},\mathbf{z}}\left(  \cdot\right)  $\ the corresponding
solutions of the following coupled system of SDEs:%
\begin{equation}
\left\{
\begin{array}
[c]{l}%
d\hat{X}^{s,\mathbf{z}}\left(  \tau\right)  =b\left(  \tau,\hat{X}%
^{s,\mathbf{z}}\left(  \tau\right)  ,\hat{\varphi}^{\tau}\left(  \hat
{X}^{s,\mathbf{z}}\left(  \tau\right)  \right)  \right)  d\tau+\sigma\left(
\tau,\hat{X}^{s,\mathbf{z}}\left(  \tau\right)  ,\hat{\varphi}^{\tau}\left(
\hat{X}^{s,\mathbf{z}}\left(  \tau\right)  \right)  \right)  dW\left(
\tau\right) \\
\text{ \ \ \ \ \ \ \ \ \ \ \ \ \ \ \ }+%
{\displaystyle\int_{E}}
c\left(  \tau,\hat{X}_{-}^{s,\mathbf{z}}\left(  \tau\right)  ,\hat{\varphi
}^{\tau}\left(  \hat{X}_{-}^{s,\mathbf{z}}\left(  \tau\right)  \right)
,e\right)  \tilde{N}\left(  d\tau,de\right)  \text{,}\\
dX^{s,\mathbf{x},\mathbf{z}}\left(  \tau\right)  =b\left(  \tau
,X^{s,\mathbf{x},\mathbf{z}}\left(  \tau\right)  ,\hat{\varphi}^{\tau}\left(
\hat{X}^{s,\mathbf{z}}\left(  \tau\right)  \right)  \right)  d\tau
+\sigma\left(  \tau,X^{s,\mathbf{x},\mathbf{z}}\left(  \tau\right)
,\hat{\varphi}^{\tau}\left(  \hat{X}^{s,\mathbf{z}}\left(  \tau\right)
\right)  \right)  dW\left(  \tau\right) \\
\text{ \ \ \ \ \ \ \ \ \ \ \ \ \ \ \ }+%
{\displaystyle\int_{E}}
c\left(  \tau,X_{-}^{s,\mathbf{x},\mathbf{z}}\left(  \tau\right)
,\hat{\varphi}^{\tau}\left(  \hat{X}_{-}^{s,\mathbf{z}}\left(  \tau\right)
\right)  ,e\right)  \tilde{N}\left(  d\tau,de\right)  \text{, for }\tau
\in\left[  s,T\right]  \text{,}\\
\hat{X}^{s,\mathbf{z}}\left(  s\right)  =\mathbf{z}\text{, }X^{s,\mathbf{x}%
,\mathbf{z}}\left(  s\right)  =\mathbf{x}\text{.}%
\end{array}
\right.  \label{Eq9}%
\end{equation}
Observe that for each fixed $\left(  s,\mathbf{z}\right)  \in\left[
t,T\right]  \times\mathbb{L}^{4}\left(  \Omega,\mathcal{F}_{s},\mathbb{P};%
\mathbb{R}
^{n}\right)  $, we have%
\begin{equation}
X^{s,\mathbf{z},\mathbf{z}}\left(  \cdot\right)  =\hat{X}^{s,\mathbf{z}%
}\left(  \cdot\right)  \text{.} \label{Eq10}%
\end{equation}

By the definition of the cost functional, we have%
\begin{align}
\mathbf{J}\left(  t,\hat{X}^{x_{0}}\left(  t\right)  ;\hat{\varphi}\left(
\cdot,\cdot\right)  \right)   &  :=\mathbb{E}_{t}\left[  \int_{t}^{T}f\left(
t,\tau,\hat{X}^{x_{0}}\left(  \tau\right)  ,\hat{\varphi}^{s}\left(  \hat
{X}^{x_{0}}\left(  \tau\right)  \right)  \right)  ds+F\left(  t,\hat{X}%
^{x_{0}}\left(  T\right)  \right)  \right] \nonumber\\
&  +G\left(  t,\hat{X}^{x_{0}}\left(  t\right)  ,\mathbb{E}_{t}\left[
\Psi\left(  \hat{X}^{x_{0}}\left(  T\right)  \right)  \right]  \right)
\label{Eq11}%
\end{align}
and since $u^{\varepsilon}\left(  \tau\right)  \equiv v\chi_{\left[
t,t+\varepsilon\right]  }\left(  \tau\right)  +\hat{\varphi}^{\tau}\left(
\hat{X}_{-}^{x_{0}}\left(  \tau\right)  \right)  \chi_{\left(  t+\varepsilon
,T\right)  }\left(  \tau\right)  $\footnote{Observe that the control
$\hat{\varphi}^{\tau}\left(  \hat{X}_{-}^{x_{0}}\left(  \tau\right)  \right)
$ in $\left(  t+\varepsilon,T\right)  $ is not influenced by the initial pair
$\left(  t+\varepsilon,X^{\varepsilon}\left(  t+\varepsilon\right)  \right)
$; this is quite different from the closed-loop control concept adopted by
Bj\"{o}rk et al. \cite{Bjork4}.}, the objective functional of $u^{\varepsilon
}\left(  \cdot\right)  $ can be written as%
\begin{align*}
\mathbf{J}\left(  t,\hat{X}^{x_{0}}\left(  t\right)  ;u^{\varepsilon}\left(
\cdot\right)  \right)   &  =\mathbb{E}_{t}\left[  \int_{t}^{T}f\left(
t,\tau,X^{\varepsilon}\left(  \tau\right)  ,u^{\varepsilon}\left(
\tau\right)  \right)  d\tau+F\left(  t,X^{\varepsilon}\left(  T\right)
\right)  \right] \\
&  +G\left(  t,\hat{X}^{x_{0}}\left(  t\right)  ,\mathbb{E}_{t}\left[
\Psi\left(  X^{\varepsilon}\left(  T\right)  \right)  \right]  \right) \\
&  =\mathbb{E}_{t}\left[  \int_{t}^{t+\varepsilon}f\left(  t,\tau
,X^{\varepsilon}\left(  \tau\right)  ,v\right)  d\tau+\int_{t+\varepsilon}%
^{T}f\left(  t,\tau,X^{\varepsilon}\left(  \tau\right)  ,\hat{\varphi}^{\tau
}\left(  \hat{X}^{x_{0}}\left(  \tau\right)  \right)  \right)  d\tau\right. \\
&  \left.  +F\left(  t,X^{\varepsilon}\left(  T\right)  \right)  +G\left(
t,\hat{X}^{x_{0}}\left(  t\right)  ,\mathbb{E}_{t}\left[  \Psi\left(
X^{\varepsilon}\left(  T\right)  \right)  \right]  \right)  \right]  .
\end{align*}

On the other hand, it follows from the uniqueness of solutions to SDEs
(\ref{Eq8*})-(\ref{Eq8}) that for all $\tau\in\left[  t+\varepsilon,T\right]
,$%
\[
X^{\varepsilon}\left(  \tau\right)  =X^{t+\varepsilon,X^{\varepsilon}\left(
t+\varepsilon\right)  ,\hat{X}^{x_{0}}\left(  t+\varepsilon\right)  }\left(
\tau\right)  ,\text{ }a.s.
\]
and%
\[
\hat{X}^{x_{0}}\left(  \tau\right)  =\hat{X}^{t+\varepsilon,\hat{X}^{x_{0}%
}\left(  t+\varepsilon\right)  }\left(  \tau\right)  ,\text{ }a.s.,
\]
where $\hat{X}^{t+\varepsilon,\hat{X}^{x_{0}}\left(  t+\varepsilon\right)
}\left(  \cdot\right)  $ and $X^{t+\varepsilon,X^{\varepsilon}\left(
t+\varepsilon\right)  ,\hat{X}^{x_{0}}\left(  t+\varepsilon\right)  }\left(
\cdot\right)  $ are the corresponding solutions of (\ref{Eq9}) with
$s=t+\varepsilon$, $\mathbf{x}=X^{\varepsilon}\left(  t+\varepsilon\right)  $
and $\mathbf{z}=\hat{X}^{x_{0}}\left(  t+\varepsilon\right)  $.

Accordingly, the objective functional of $u^{\varepsilon}\left(  \cdot\right)
$ can be rewritten as%
\begin{align}
&  \mathbf{J}\left(  t,\hat{X}^{x_{0}}\left(  t\right)  ;u^{\varepsilon
}\left(  \cdot\right)  \right) \nonumber\\
&  =\mathbb{E}_{t}\left[  \int_{t}^{t+\varepsilon}f\left(  t,\tau
,X^{\varepsilon}\left(  \tau\right)  ,v\right)  d\tau\right. \nonumber\\
&  +\int_{t+\varepsilon}^{T}f\left(  t,\tau,X^{t+\varepsilon,X^{\varepsilon
}\left(  t+\varepsilon\right)  ,\hat{X}^{x_{0}}\left(  t+\varepsilon\right)
}\left(  \tau\right)  ,\hat{\varphi}^{\tau}\left(  \hat{X}^{t+\varepsilon
,\hat{X}^{x_{0}}\left(  t+\varepsilon\right)  }\left(  \tau\right)  \right)
\right)  d\tau\nonumber\\
&  \left.  +F\left(  t,X^{t+\varepsilon,X^{\varepsilon}\left(  t+\varepsilon
\right)  ,\hat{X}^{x_{0}}\left(  t+\varepsilon\right)  }\left(  T\right)
\right)  +G\left(  t,\hat{X}^{x_{0}}\left(  t\right)  ,\mathbb{E}_{t}\left[
\Psi\left(  X^{t+\varepsilon,X^{\varepsilon}\left(  t+\varepsilon\right)
,\hat{X}^{x_{0}}\left(  t+\varepsilon\right)  }\left(  T\right)  \right)
\right]  \right)  \right]  \text{.} \label{Eq12}%
\end{align}

Now, if we define for all $\left(  t,s,\mathbf{x},\mathbf{z}\right)  \in
D\left[  0,T\right]  \times\mathbb{L}^{4}\left(  \Omega,\mathcal{F}%
_{s},\mathbb{P};%
\mathbb{R}
^{n}\right)  \times\mathbb{L}^{4}\left(  \Omega,\mathcal{F}_{s},\mathbb{P};%
\mathbb{R}
^{n}\right)  $,%
\begin{equation}
h^{t}\left(  s,\mathbf{x},\mathbf{z}\right)  :=\mathbb{E}_{s}\left[  \int
_{s}^{T}f\left(  t,\tau,X^{s,\mathbf{x},\mathbf{z}}\left(  \tau\right)
,\hat{\varphi}^{\tau}\left(  \hat{X}^{s,\mathbf{z}}\left(  \tau\right)
\right)  \right)  d\tau+F\left(  t,X^{s,\mathbf{x},\mathbf{z}}\left(
T\right)  \right)  \right]  \label{Eq13}%
\end{equation}
and%
\begin{equation}
\eta\left(  s,\mathbf{x},\mathbf{z}\right)  :=\mathbb{E}_{s}\left[
\Psi\left(  X^{s,\mathbf{x},\mathbf{z}}\left(  T\right)  \right)  \right]
\text{,} \label{Eq14}%
\end{equation}
then the objective values (\ref{Eq11}) and (\ref{Eq12}) can be rewritten as%
\begin{align*}
\mathbf{J}\left(  t,\hat{X}^{x_{0}}\left(  t\right)  ;\hat{\varphi}\left(
\cdot,\cdot\right)  \right)   &  =h^{t}\left(  t,\hat{X}^{x_{0}}\left(
t\right)  ,\hat{X}^{x_{0}}\left(  t\right)  \right)  +G\left(  t,\hat
{X}^{x_{0}}\left(  t\right)  ,\eta\left(  t,\hat{X}^{x_{0}}\left(  t\right)
,\hat{X}^{x_{0}}\left(  t\right)  \right)  \right) \\
&  =h^{t}\left(  t,X^{\varepsilon}\left(  t\right)  ,\hat{X}^{x_{0}}\left(
t\right)  \right)  +G\left(  t,\hat{X}^{x_{0}}\left(  t\right)  ,\eta\left(
t,X^{\varepsilon}\left(  t\right)  ,\hat{X}^{x_{0}}\left(  t\right)  \right)
\right)
\end{align*}
and%
\begin{align*}
\mathbf{J}\left(  t,\hat{X}^{x_{0}}\left(  t\right)  ;u^{\varepsilon}\left(
\cdot\right)  \right)   &  =\mathbb{E}_{t}\left[  \int\limits_{t}%
^{t+\varepsilon}f\left(  t,\tau,X^{\varepsilon}\left(  \tau\right)  ,v\right)
d\tau+h^{t}\left(  t+\varepsilon,X^{\varepsilon}\left(  t+\varepsilon\right)
,\hat{X}^{x_{0}}\left(  t+\varepsilon\right)  \right)  \right] \\
&  +G\left(  t,\hat{X}^{x_{0}}\left(  t\right)  ,\mathbb{E}_{t}\left[
\eta\left(  t+\varepsilon,X^{\varepsilon}\left(  t+\varepsilon\right)
,\hat{X}^{x_{0}}\left(  t+\varepsilon\right)  \right)  \right]  \right)
\text{,}%
\end{align*}
respectively. Consequently, $\mathbf{\Delta\hat{J}}^{\varepsilon}\left(
t\right)  $ takes the form%
\begin{equation}
\mathbf{\Delta\hat{J}}^{\varepsilon}\left(  t\right)  =\mathbb{E}_{t}\left[
\int_{t}^{t+\varepsilon}f\left(  t,\tau,X^{\varepsilon}\left(  \tau\right)
,v\right)  d\tau+\hat{\Delta}_{1}^{\varepsilon}+\hat{\Delta}_{2}^{\varepsilon
}\right]  , \label{Eq15}%
\end{equation}
where%
\[
\hat{\Delta}_{1}^{\varepsilon}:=h^{t}\left(  t+\varepsilon,X^{\varepsilon
}\left(  t+\varepsilon\right)  ,\hat{X}^{x_{0}}\left(  t+\varepsilon\right)
\right)  -h^{t}\left(  t,X^{\varepsilon}\left(  t\right)  ,\hat{X}^{x_{0}%
}\left(  t\right)  \right)
\]
and%
\begin{align*}
\hat{\Delta}_{2}^{\varepsilon}  &  =:G\left(  t,\hat{X}^{x_{0}}\left(
t\right)  ,\mathbb{E}_{t}\left[  \eta\left(  t+\varepsilon,X^{\varepsilon
}\left(  t+\varepsilon\right)  ,\hat{X}^{x_{0}}\left(  t+\varepsilon\right)
\right)  \right]  \right) \\
&  -G\left(  t,\hat{X}^{x_{0}}\left(  t\right)  ,\eta\left(  t,X^{\varepsilon
}\left(  t\right)  ,\hat{X}^{x_{0}}\left(  t\right)  \right)  \right)
\text{.}%
\end{align*}

However, the above expression of $\mathbf{\Delta\hat{J}}^{\varepsilon}\left(
t\right)  $ is still very difficult to handle, since the terms $\hat{\Delta
}_{1}^{\varepsilon}$ and $\hat{\Delta}_{2}^{\varepsilon}$ involved on the
right-hand side of (\ref{Eq15}) are too complicated. Thus, we need to explore
(\ref{Eq15}) further, in order to be able to evaluate the limit in
(\ref{Eq7}). Inspired by the idea of Four Step Scheme introduced in \cite{MY}
for FBSDEs with jumps, we proceed as follows. For each fixed $\left(
t,s,\mathbf{x},\mathbf{z}\right)  \in D\left[  0,T\right]  \times
\mathbb{L}^{4}\left(  \Omega,\mathcal{F}_{s},\mathbb{P};%
\mathbb{R}
^{n}\right)  \times\mathbb{L}^{4}\left(  \Omega,\mathcal{F}_{s},\mathbb{P};%
\mathbb{R}
^{n}\right)  $, we introduce the following system of FBSDEs:%

\begin{equation}
\left\{
\begin{array}
[c]{l}%
d\hat{X}^{s,\mathbf{z}}\left(  \tau\right)  =b\left(  \tau,\hat{X}%
^{s,\mathbf{z}}\left(  \tau\right)  ,\hat{\varphi}^{\tau}\left(  \hat
{X}^{s,\mathbf{z}}\left(  \tau\right)  \right)  \right)  d\tau+\sigma\left(
\tau,\hat{X}^{s,\mathbf{z}}\left(  \tau\right)  ,\hat{\varphi}^{\tau}\left(
\hat{X}^{s,\mathbf{z}}\left(  \tau\right)  \right)  \right)  dW\left(
\tau\right) \\
\text{ \ \ \ \ \ \ \ \ \ \ \ \ \ \ \ }+%
{\displaystyle\int_{E}}
c\left(  \tau,\hat{X}_{-}^{s,\mathbf{z}}\left(  \tau\right)  ,\hat{\varphi
}^{\tau}\left(  \hat{X}_{-}^{s,\mathbf{z}}\left(  \tau\right)  \right)
,e\right)  \tilde{N}\left(  ds,de\right)  \text{,}\\
dX^{s,\mathbf{x},\mathbf{z}}\left(  \tau\right)  =b\left(  \tau
,X^{s,\mathbf{x},\mathbf{z}}\left(  \tau\right)  ,\hat{\varphi}^{\tau}\left(
\hat{X}^{s,\mathbf{z}}\left(  \tau\right)  \right)  \right)  d\tau
+\sigma\left(  \tau,X^{s,\mathbf{x},\mathbf{z}}\left(  \tau\right)
,\hat{\varphi}^{\tau}\left(  \hat{X}^{s,\mathbf{z}}\left(  \tau\right)
\right)  \right)  dW\left(  \tau\right) \\
\text{ \ \ \ \ \ \ \ \ \ \ \ \ \ \ \ }+%
{\displaystyle\int_{E}}
c\left(  \tau,X_{-}^{s,\mathbf{x},\mathbf{z}}\left(  \tau\right)
,\hat{\varphi}^{\tau}\left(  \hat{X}_{-}^{s,\mathbf{z}}\left(  \tau\right)
\right)  ,e\right)  \tilde{N}\left(  ds,de\right)  \text{,}\\
dY^{s,\mathbf{x},\mathbf{z}}\left(  \tau;t\right)  =-f\left(  t,\tau
,X^{s,\mathbf{x},\mathbf{z}}\left(  \tau\right)  ,\hat{\varphi}^{\tau}\left(
\hat{X}^{s,\mathbf{z}}\left(  \tau\right)  \right)  \right)  d\tau
+Q^{s,\mathbf{x},\mathbf{z}}\left(  \tau;t\right)  ^{\top}dW\left(
\tau\right) \\
\text{ \ \ \ \ \ \ \ \ \ \ \ \ \ \ \ }+%
{\displaystyle\int_{E}}
R^{s,\mathbf{x},\mathbf{z}}\left(  \tau,e;t\right)  \tilde{N}\left(
d\tau,de\right)  \text{,}\\
dK^{s,\mathbf{x},\mathbf{z}}\left(  \tau\right)  =H^{s,\mathbf{x},\mathbf{z}%
}\left(  \tau\right)  dW\left(  \tau\right)  +%
{\displaystyle\int_{E}}
\Gamma^{s,\mathbf{x},\mathbf{z}}\left(  \tau,e\right)  \tilde{N}\left(
d\tau,de\right)  \text{, for }\tau\in\left[  s,T\right]  \text{,}%
\end{array}
\right.  \label{Eq17}%
\end{equation}
with the boundary conditions%
\[
\hat{X}^{s,\mathbf{z}}\left(  s\right)  =\mathbf{z}\text{, }X^{s,\mathbf{x}%
,\mathbf{z}}\left(  s\right)  =\mathbf{x}\text{, }K^{s,\mathbf{x},\mathbf{z}%
}\left(  T\right)  =\Psi\left(  X^{s,\mathbf{x},\mathbf{z}}\left(  T\right)
\right)  \text{, }Y^{s,\mathbf{x},\mathbf{z}}\left(  T;t\right)  =F\left(
t,X^{s,\mathbf{x},\mathbf{z}}\left(  T\right)  \right)  \text{.}%
\]
Under Assumptions \textbf{(H1)-(H2) }and by virtue of Definition
\ref{Def1}\textbf{, }the above system of FBSDEs admits a unique adapted
solution $\left(  \hat{X}^{s,\mathbf{z}}\left(  \cdot\right)  ,X^{s,\mathbf{x}%
,\mathbf{z}}\left(  \cdot\right)  \right)  \in\mathcal{S}_{%
\mathcal{F}%
}^{4}\left(  s,T;%
\mathbb{R}
^{n}\right)  ^{2}$, $\left(  Y^{s,\mathbf{x},\mathbf{z}}\left(  \cdot
;t\right)  ,Q^{s,\mathbf{x},\mathbf{z}}\left(  \cdot;t\right)
,R^{s,\mathbf{x},\mathbf{z}}\left(  \cdot,\cdot;t\right)  \right)
\in\mathcal{S}_{%
\mathcal{F}%
}^{2}\left(  s,T;%
\mathbb{R}
\right)  \times\mathcal{L}_{%
\mathcal{F}%
}^{2}\left(  s,T;%
\mathbb{R}
^{d}\right)  \times\mathcal{L}_{%
\mathcal{F}%
,p}^{\vartheta,q}\left(  \left[  s,T\right]  \times E;%
\mathbb{R}
\right)  $, $\left(  K^{s,\mathbf{x},\mathbf{z}}\left(  \cdot\right)
,H^{s,\mathbf{x},\mathbf{z}}\left(  \cdot\right)  ,\Gamma^{s,\mathbf{x}%
,\mathbf{z}}\left(  \cdot,\cdot\right)  \right)  \in\mathcal{S}_{%
\mathcal{F}%
}^{2}\left(  s,T;%
\mathbb{R}
^{m}\right)  $\newline$\times\mathcal{L}_{%
\mathcal{F}%
}^{2}\left(  s,T;%
\mathbb{R}
^{m\times d}\right)  \times\mathcal{L}_{%
\mathcal{F}%
,p}^{\vartheta,q}\left(  \left[  s,T\right]  \times E;%
\mathbb{R}
^{m}\right)  $ uniquely depending on $\left(  s,\mathbf{x},\mathbf{z}\right)
$ (see e.g. \cite{Romuald}). Further, one has%
\begin{align*}
Y^{s,\mathbf{x},\mathbf{z}}\left(  s;t\right)   &  =\mathbb{E}_{s}\left[
\int_{s}^{T}f\left(  t,\tau,X^{s,\mathbf{x},\mathbf{z}}\left(  \tau\right)
,\hat{\varphi}^{\tau}\left(  \hat{X}^{s,\mathbf{z}}\left(  \tau\right)
\right)  \right)  d\tau+F\left(  t,X^{s,\mathbf{x},\mathbf{z}}\left(
T\right)  \right)  \right] \\
&  =h^{t}\left(  s,\mathbf{x},\mathbf{z}\right)
\end{align*}
and%
\[
K^{s,\mathbf{x},\mathbf{z}}\left(  s\right)  =\mathbb{E}_{s}\left[
\Psi\left(  \hat{X}^{s,\mathbf{x},\mathbf{z}}\left(  T\right)  \right)
\right]  =\eta\left(  s,\mathbf{x},\mathbf{z}\right)  \text{.}%
\]

It is seen that (\ref{Eq17}) is a system of FBSDEs with deterministic
coefficients. Inspired by Theorem \ref{result00}, we introduce the following
systems of IPDEs:%
\begin{equation}
\left\{
\begin{array}
[c]{l}%
0=\hat{\theta}_{s}^{t}\left(  s,x,z\right)  +\left\langle \hat{\theta}_{x}%
^{t}\left(  s,x,z\right)  ,\mu\left(  s,x,\hat{\varphi}^{s}\left(  z\right)
\right)  \right\rangle +\left\langle \hat{\theta}_{z}^{t}\left(  s,x,z\right)
,\mu\left(  s,z,\hat{\varphi}^{s}\left(  z\right)  \right)  \right\rangle \\
\text{ \ \ }+\frac{1}{2}\text{\textbf{tr}}\left[  \sigma\sigma^{\top}\left(
s,x,\hat{\varphi}^{s}\left(  z\right)  \right)  \hat{\theta}_{xx}^{t}\left(
s,x,z\right)  \right]  +\frac{1}{2}\text{\textbf{tr}}\left[  \sigma
\sigma^{\top}\left(  s,z,\hat{\varphi}^{s}\left(  z\right)  \right)
\hat{\theta}_{zz}^{t}\left(  s,x,z\right)  \right] \\
\text{ \ }+\text{\textbf{tr}}\left[  \sigma\left(  s,x,\hat{\varphi}%
^{s}\left(  z\right)  \right)  ^{\top}\hat{\theta}_{xz}^{t}\left(
s,x,z\right)  \sigma\left(  s,z,\hat{\varphi}^{s}\left(  z\right)  \right)
\right]  +\text{\ }f\left(  t,s,x,\hat{\varphi}^{s}\left(  z\right)  \right)
\\
\text{ \ \ }+%
{\displaystyle\int_{E}}
\left\{  \hat{\theta}^{t}\left(  s,x+c\left(  s,x,\hat{\varphi}^{s}\left(
z\right)  ,e\right)  ,z+c\left(  s,z,\hat{\varphi}^{s}\left(  z\right)
,e\right)  \right)  -\hat{\theta}^{t}\left(  s,x,z\right)  \right. \\
\text{\ }%
\begin{array}
[c]{r}%
\left.  -\left\langle \hat{\theta}_{x}^{t}\left(  s,x,z\right)  ,c\left(
s,x,\hat{\varphi}^{s}\left(  z\right)  ,e\right)  \right\rangle -\left\langle
\hat{\theta}_{z}^{t}\left(  s,x,z\right)  ,c\left(  s,z,\hat{\varphi}%
^{s}\left(  z\right)  ,e\right)  \right\rangle \right\}  \vartheta\left(
de\right)  \text{, }\\
\text{for }\left(  s,x,z\right)  \in\left[  t,T\right]  \times\mathbb{R}%
^{n}\times\mathbb{R}^{n}\text{,}%
\end{array}
\\
\hat{\theta}^{t}\left(  T,x,z\right)  =F\left(  t,x\right)  \text{, for
}\left(  x,z\right)  \in\mathbb{R}^{n}\times\mathbb{R}^{n}%
\end{array}
\right.  \label{Eq19}%
\end{equation}
and%
\begin{equation}
\left\{
\begin{array}
[c]{l}%
0=\hat{g}_{s}^{i}\left(  s,x,z\right)  +\left\langle \hat{g}_{x}^{i}\left(
s,x,z\right)  ,\mu\left(  s,x,\hat{\varphi}^{s}\left(  z\right)  \right)
\right\rangle +\left\langle \hat{g}_{z}^{i}\left(  s,x,z\right)  ,\mu\left(
s,z,\hat{\varphi}^{s}\left(  z\right)  \right)  \right\rangle \\
\text{ \ \ }+\frac{1}{2}\text{\textbf{tr}}\left[  \sigma\sigma^{\top}\left(
s,x,\hat{\varphi}^{s}\left(  z\right)  \right)  \hat{g}_{xx}^{i}\left(
s,x,z\right)  \right]  +\frac{1}{2}\text{\textbf{tr}}\left[  \sigma
\sigma^{\top}\left(  s,z,\hat{\varphi}^{s}\left(  z\right)  \right)  \hat
{g}_{zz}^{i}\left(  s,x,z\right)  \right] \\
\text{ \ }+\text{\textbf{tr}}\left[  \sigma\left(  s,x,\hat{\varphi}%
^{s}\left(  z\right)  \right)  ^{\top}\hat{g}_{xz}^{i}\left(  s,x,z\right)
\sigma\left(  s,z,\hat{\varphi}^{s}\left(  z\right)  \right)  \right] \\
\text{ \ }+%
{\displaystyle\int_{E}}
\left\{  \hat{g}^{i}\left(  s,x+c\left(  s,x,\hat{\varphi}^{s}\left(
z\right)  ,e\right)  ,z+c\left(  s,z,\hat{\varphi}^{s}\left(  z\right)
,e\right)  \right)  -\hat{g}^{i}\left(  s,x,z\right)  \right. \\%
\begin{array}
[c]{r}%
\left.  -\left\langle \hat{g}_{x}^{i}\left(  s,x,z\right)  ,c\left(
s,x,\hat{\varphi}^{s}\left(  z\right)  ,e\right)  \right\rangle -\left\langle
\hat{g}_{z}^{i}\left(  s,x,z\right)  ,c\left(  s,z,\hat{\varphi}^{s}\left(
z\right)  ,e\right)  \right\rangle \right\}  \vartheta\left(  de\right)
\text{, }\\
\text{for }\left(  s,x,z\right)  \in\left[  0,T\right]  \times\mathbb{R}%
^{n}\times\mathbb{R}^{n}\text{,}%
\end{array}
\\
\hat{g}^{i}\left(  T,x,z\right)  =\Psi_{i}\left(  x\right)  \text{, for
}\left(  x,z\right)  \in\mathbb{R}^{n}\times\mathbb{R}^{n}\text{; }1\leq i\leq
m\text{,}%
\end{array}
\right.  \label{Eq21}%
\end{equation}
where $\Psi_{i}\left(  \cdot\right)  $ is the i-th coordinate of $\Psi\left(
\cdot\right)  =\left(  \Psi_{1}\left(  \cdot\right)  ,...,\Psi_{m}\left(
\cdot\right)  \right)  ^{\top}$. If the above admit unique classical solutions
$\hat{\theta}^{t}\left(  \cdot,\cdot,\cdot\right)  $ and $\hat{g}\left(
\cdot,\cdot,\cdot\right)  =\left(  \hat{g}^{1}\left(  \cdot,\cdot
,\cdot\right)  ,...,\hat{g}^{m}\left(  \cdot,\cdot,\cdot\right)  \right)
^{\top}$, respectively, such that $\hat{\theta}^{t}\left(  \cdot,\cdot
,\cdot\right)  \in\mathcal{C}_{2}^{1,2,2}\left(  \left[  t,T\right]  \times%
\mathbb{R}
^{n}\times%
\mathbb{R}
^{n};%
\mathbb{R}
\right)  $ and $\hat{g}\left(  \cdot,\cdot,\cdot\right)  \in\mathcal{C}%
_{2}^{1,2,2}\left(  \left[  0,T\right]  \times%
\mathbb{R}
^{n}\times%
\mathbb{R}
^{n};%
\mathbb{R}
\right)  ^{m}$, then the following representations hold: For all $\tau
\in\left[  s,T\right]  $, a.s.,%
\begin{equation}
Y^{s,\mathbf{x},\mathbf{z}}\left(  \tau;t\right)  =\hat{\theta}^{t}\left(
\tau,X^{s,\mathbf{x},\mathbf{z}}\left(  \tau\right)  ,\hat{X}^{s,\mathbf{z}%
}\left(  \tau\right)  \right)  \label{Eq18}%
\end{equation}
and%
\begin{align}
K^{s,\mathbf{x},\mathbf{z}}\left(  \tau\right)   &  =\hat{g}\left(
\tau,X^{s,\mathbf{x},\mathbf{z}}\left(  \tau\right)  ,\hat{X}^{s,\mathbf{z}%
}\left(  \tau\right)  \right) \nonumber\\
&  =\left(  \hat{g}^{1}\left(  \tau,X^{s,\mathbf{x},\mathbf{z}}\left(
\tau\right)  ,\hat{X}^{s,\mathbf{z}}\left(  \tau\right)  \right)  ,...,\hat
{g}^{m}\left(  \tau,X^{s,\mathbf{x},\mathbf{z}}\left(  \tau\right)  ,\hat
{X}^{s,\mathbf{z}}\left(  \tau\right)  \right)  \right)  ^{\top}\text{.}
\label{Eq20}%
\end{align}
Note that by virtue of Definition \ref{Def3}, it is not difficult to verify
that all needed integrability conditions to the above representations are satisfied.

Thus, setting $\tau=s$ in (\ref{Eq18})-(\ref{Eq20}), we get%
\begin{equation}
\hat{\theta}^{t}\left(  s,\mathbf{x},\mathbf{z}\right)  =Y\left(  s;t\right)
=h^{t}\left(  s,\mathbf{x},\mathbf{z}\right)  \label{Eq22}%
\end{equation}
and%
\begin{equation}
\hat{g}\left(  s,\mathbf{x},\mathbf{z}\right)  =K\left(  s\right)
=\eta\left(  s,\mathbf{x},\mathbf{z}\right)  \text{.} \label{Eq23}%
\end{equation}
Hence we can rewrite $\hat{\Delta}_{1}^{\varepsilon}$ and $\hat{\Delta}%
_{2}^{\varepsilon}$ as follows:%
\[
\hat{\Delta}_{1}^{\varepsilon}:=\hat{\theta}^{t}\left(  t+\varepsilon
,X^{\varepsilon}\left(  t+\varepsilon\right)  ,\hat{X}^{x_{0}}\left(
t+\varepsilon\right)  \right)  -\hat{\theta}^{t}\left(  t,X^{\varepsilon
}\left(  t\right)  ,\hat{X}^{x_{0}}\left(  t\right)  \right)
\]
and%
\begin{align*}
\hat{\Delta}_{2}^{\varepsilon}  &  :=G\left(  t,\hat{X}^{x_{0}}\left(
t\right)  ,\mathbb{E}_{t}\left[  \hat{g}\left(  t+\varepsilon,X^{\varepsilon
}\left(  t+\varepsilon\right)  ,\hat{X}^{x_{0}}\left(  t+\varepsilon\right)
\right)  \right]  \right) \\
&  -G\left(  t,\hat{X}^{x_{0}}\left(  t\right)  ,\hat{g}\left(
t,X^{\varepsilon}\left(  t\right)  ,\hat{X}^{x_{0}}\left(  t\right)  \right)
\right)  \text{.}%
\end{align*}
Invoking this into (\ref{Eq15}), we obtain that%
\begin{align}
\mathbf{\Delta\hat{J}}^{\varepsilon}\left(  t\right)   &  =\mathbb{E}%
_{t}\left[  \int_{t}^{t+\varepsilon}f\left(  t,s,X^{\varepsilon}\left(
s\right)  ,v\right)  ds\right] \nonumber\\
&  +\mathbb{E}_{t}\left[  \hat{\theta}^{t}\left(  t+\varepsilon,X^{\varepsilon
}\left(  t+\varepsilon\right)  ,\hat{X}^{x_{0}}\left(  t+\varepsilon\right)
\right)  -\hat{\theta}^{t}\left(  t,X^{\varepsilon}\left(  t\right)  ,\hat
{X}^{x_{0}}\left(  t\right)  \right)  \right] \nonumber\\
&  +G\left(  t,\hat{X}^{x_{0}}\left(  t\right)  ,\mathbb{E}_{t}\left[  \hat
{g}\left(  t+\varepsilon,X^{\varepsilon}\left(  t+\varepsilon\right)  ,\hat
{X}^{x_{0}}\left(  t+\varepsilon\right)  \right)  \right]  \right)  -G\left(
t,\hat{X}^{x_{0}}\left(  t\right)  ,\hat{g}\left(  t,X^{\varepsilon}\left(
t\right)  ,\hat{X}^{x_{0}}\left(  t\right)  \right)  \right)  \text{.}%
\nonumber
\end{align}

Before going further, let us introduce an $\mathcal{H}$-function associated to
the 3-tuple $\left(  \hat{\varphi}\left(  \cdot,\cdot\right)  ,\hat{\theta
}^{t}\left(  \cdot,\cdot,\cdot\right)  ,\hat{g}\left(  \cdot,\cdot
,\cdot\right)  \right)  $ as follows: For each $\left(  t,s\right)  \in
D\left[  0,T\right]  $, $\left(  u,y\right)  \in U\times%
\mathbb{R}
^{n}$ and for each $X,Z\in\mathbb{L}^{4}\left(  \Omega,\mathcal{F}%
_{s},\mathbb{P};%
\mathbb{R}
^{n}\right)  $%
\begin{align}
&  \mathcal{H}\left(  t,y,s,X,Z,u\right) \nonumber\\
&  :=\left\langle \hat{\theta}_{x}^{t}\left(  s,X,Z\right)  +{%
{\textstyle\sum\limits_{i=1}^{m}}
}G_{\bar{x}_{i}}\left(  t,y,\mathbb{E}_{t}\left[  \hat{g}\left(  s,X,Z\right)
\right]  \right)  \hat{g}_{x}^{i}\left(  s,X,Z\right)  ,\mu\left(
s,X,u\right)  \right\rangle \nonumber\\
&  +\frac{1}{2}\text{\textbf{tr}}\left[  \left(  \hat{\theta}_{xx}^{t}\left(
s,X,Z\right)  +{%
{\textstyle\sum\limits_{i=1}^{m}}
}G_{\bar{x}_{i}}\left(  t,y,\mathbb{E}_{t}\left[  \hat{g}\left(  s,X,Z\right)
\right]  \right)  \hat{g}_{xx}^{i}\left(  s,X,Z\right)  \right)  \left(
\sigma\sigma^{\top}\right)  \left(  s,X,u\right)  \right] \nonumber\\
&  +\text{\textbf{tr}}\left[  \sigma\left(  s,X,u\right)  ^{\top}\left(
\hat{\theta}_{xz}^{t}\left(  s,X,Z\right)  +{%
{\textstyle\sum\limits_{i=1}^{m}}
}G_{\bar{x}_{i}}\left(  t,y,\mathbb{E}_{t}\left[  \hat{g}\left(  s,X,Z\right)
\right]  \right)  \hat{g}_{xz}^{i}\left(  s,X,Z\right)  \right)  \sigma\left(
s,Z,\hat{\varphi}^{s}\left(  Z\right)  \right)  \right] \nonumber\\
&  +%
{\displaystyle\int_{E}}
\left\{  \hat{\theta}^{t}\left(  s,X+c\left(  s,X,u,e\right)  ,Z+c\left(
s,Z,\hat{\varphi}^{s}\left(  Z\right)  ,e\right)  \right)  -\left\langle
\hat{\theta}_{x}^{t}\left(  s,X,Z\right)  ,c\left(  s,X,u,e\right)
\right\rangle \right\}  \vartheta\left(  de\right) \nonumber\\
&  +{%
{\textstyle\sum\limits_{i=1}^{m}}
}G_{\bar{x}_{i}}\left(  t,y,\mathbb{E}_{t}\left[  \hat{g}\left(  s,X,Z\right)
\right]  \right) \nonumber\\
&  \times%
{\displaystyle\int_{E}}
\left\{  \hat{g}^{i}\left(  s,X+c\left(  s,X,u,e\right)  ,Z+c\left(
s,Z,\hat{\varphi}^{s}\left(  Z\right)  ,e\right)  \right)  -\left\langle
\hat{g}_{x}^{i}\left(  s,X,Z\right)  ,c\left(  s,X,u,e\right)  \right\rangle
\right\}  \vartheta\left(  de\right) \nonumber\\
&  +f\left(  t,y,s,X,u\right)  \text{,} \label{Eq25}%
\end{align}
where for all $\left(  t,y,\bar{x}\right)  \in\left[  0,T\right]
\times\mathbb{R}^{n}\times\mathbb{R}^{m},$ $G_{\bar{x}_{i}}\left(  t,y,\bar
{x}\right)  =\frac{\partial}{\partial\bar{x}_{i}}G\left(  t,y,\bar{x}_{1}%
,\bar{x}_{2},...,\bar{x}_{m}\right)  $ and $\hat{g}^{i}\left(  \cdot
,\cdot,\cdot\right)  $ the i-th coordinate of $\hat{g}\left(  \cdot
,\cdot,\cdot\right)  =\left(  \hat{g}^{1}\left(  \cdot,\cdot,\cdot\right)
,...,\hat{g}^{m}\left(  \cdot,\cdot,\cdot\right)  \right)  ^{\top}$.

The above derivation can be summarized as follows.

\begin{proposition}
\label{result1}Let \textbf{(H1)-(H2)} hold. Given an admissible strategy
$\hat{\varphi}\left(  \cdot,\cdot\right)  \in\mathcal{S}$, suppose that for
each $t\in\left[  0,T\right]  $, the PDEs (\ref{Eq19}) and (\ref{Eq21}) admit
the classical solutions $\hat{\theta}^{t}\left(  \cdot,\cdot,\cdot\right)  $
and $\hat{g}\left(  \cdot,\cdot,\cdot\right)  $, respectively, such that%
\[
\hat{\theta}^{t}\left(  \cdot,\cdot,\cdot\right)  \in\mathcal{C}_{2}%
^{1,2,2}\left(  \left[  t,T\right]  \times%
\mathbb{R}
^{n}\times%
\mathbb{R}
^{n};%
\mathbb{R}
\right)
\]
and
\[
\hat{g}\left(  \cdot,\cdot,\cdot\right)  \in\mathcal{C}_{2}^{1,2,2}\left(
\left[  0,T\right]  \times%
\mathbb{R}
^{n}\times%
\mathbb{R}
^{n};%
\mathbb{R}
\right)  ^{m}.
\]
Then for each $\left(  t,s,\mathbf{x},\mathbf{z}\right)  \in D\left[
0,T\right]  \times\mathbb{L}^{4}\left(  \Omega,\mathcal{F}_{s},\mathbb{P};%
\mathbb{R}
^{n}\right)  \times\mathbb{L}^{4}\left(  \Omega,\mathcal{F}_{s},\mathbb{P};%
\mathbb{R}
^{n}\right)  $, $\hat{\theta}^{t}\left(  s,\mathbf{x},\mathbf{z}\right)  $ and
$\hat{g}\left(  s,\mathbf{x},\mathbf{z}\right)  $ have the following
probabilistic representations:%
\begin{equation}
\hat{\theta}^{t}\left(  s,\mathbf{x},\mathbf{z}\right)  =\mathbb{E}_{s}\left[
\int_{s}^{T}f\left(  t,\tau,X^{s,\mathbf{x},\mathbf{z}}\left(  \tau\right)
,\hat{\varphi}^{\tau}\left(  \hat{X}^{s,\mathbf{z}}\left(  \tau\right)
\right)  \right)  d\tau+F\left(  t,X^{s,\mathbf{x},\mathbf{z}}\left(
T\right)  \right)  \right]  \label{Eq25*}%
\end{equation}
and%
\begin{equation}
\hat{g}\left(  s,\mathbf{x},\mathbf{z}\right)  =\mathbb{E}_{s}\left[
\Psi\left(  X^{s,\mathbf{x},\mathbf{z}}\left(  T\right)  \right)  \right]
\text{,} \label{Eq26*}%
\end{equation}
respectively, where $\hat{X}^{s,\mathbf{z}}\left(  \cdot\right)  $ and
$X^{s,\mathbf{x},\mathbf{z}}\left(  \cdot\right)  $ are the corresponding
solutions of (\ref{Eq9}). Furthermore, for any $t\in\left[  0,T\right]  $,
$v\in\mathbb{L}^{4}\left(  \Omega,\mathcal{F}_{t},\mathbb{P};U\right)  $ and
for any $\varepsilon\in\left[  0,T-t\right)  $, the following equality holds%
\begin{align}
&  \mathbf{J}\left(  t,\hat{X}^{x_{0}}\left(  t\right)  ;u^{\varepsilon
}\left(  \cdot\right)  \right)  -\mathbf{J}\left(  t,\hat{X}^{x_{0}}\left(
t\right)  ;\hat{\varphi}\left(  \cdot,\cdot\right)  \right) \nonumber\\
&  =\int_{t}^{t+\varepsilon}\mathbb{E}_{t}\left[  \mathcal{H}\left(  t,\hat
{X}^{x_{0}}\left(  t\right)  ,s,X^{\varepsilon}\left(  s\right)  ,\hat
{X}^{x_{0}}\left(  s\right)  ,v\right)  -\mathcal{H}\left(  t,\hat{X}^{x_{0}%
}\left(  t\right)  ,s,X^{\varepsilon}\left(  s\right)  ,\hat{X}^{x_{0}}\left(
s\right)  ,\hat{\varphi}^{s}\left(  \hat{X}^{x_{0}}\left(  s\right)  \right)
\right)  \right]  ds\text{,} \label{Eq27}%
\end{align}
where $\mathcal{H}$ is as introduced in (\ref{Eq25}), $u^{\varepsilon}\left(
\cdot\right)  $\textit{ is defined by }(\ref{Eq6}) and $X^{\varepsilon}\left(
\cdot\right)  $ is unique strong solution of (\ref{Eq8}).
\end{proposition}

\bop The probabilistic representations in (\ref{Eq25*})-(\ref{Eq26*}) follow
from Equalities (\ref{Eq22})-(\ref{Eq23}). So we only need to show
(\ref{Eq27}). Consider the difference%
\begin{align}
&  \mathbf{J}\left(  t,\hat{X}^{x_{0}}\left(  t\right)  ;u^{\varepsilon
}\left(  \cdot\right)  \right)  -\mathbf{J}\left(  t,\hat{X}^{x_{0}}\left(
t\right)  ;\hat{\varphi}\left(  \cdot,\cdot\right)  \right) \nonumber\\
&  =\mathbb{E}_{t}\left[  \int_{t}^{t+\varepsilon}f\left(  t,s,X^{\varepsilon
}\left(  s\right)  ,v\right)  ds\right] \nonumber\\
&  +\mathbb{E}_{t}\left[  \hat{\theta}^{t}\left(  t+\varepsilon,X^{\varepsilon
}\left(  t+\varepsilon\right)  ,\hat{X}^{x_{0}}\left(  t+\varepsilon\right)
\right)  -\hat{\theta}^{t}\left(  t,X^{\varepsilon}\left(  t\right)  ,\hat
{X}^{x_{0}}\left(  t\right)  \right)  \right] \nonumber\\
&  +G\left(  t,\hat{X}^{x_{0}}\left(  t\right)  ,\mathbb{E}_{t}\left[  \hat
{g}\left(  t+\varepsilon,X^{\varepsilon}\left(  t+\varepsilon\right)  ,\hat
{X}^{x_{0}}\left(  t+\varepsilon\right)  \right)  \right]  \right) \nonumber\\
&  -G\left(  t,\hat{X}^{x_{0}}\left(  t\right)  ,\hat{g}\left(
t,X^{\varepsilon}\left(  t\right)  ,\hat{X}^{x_{0}}\left(  t\right)  \right)
\right)  \text{.} \label{Eq28}%
\end{align}

Applying It\^{o}'s formula to $\hat{\theta}^{t}\left(  \cdot,X^{\varepsilon
}\left(  \cdot\right)  ,\hat{X}^{x_{0}}\left(  \cdot\right)  \right)  $ on
time interval $\left[  t,t+\varepsilon\right]  $ and taking conditional
expectations, we obtain that%
\begin{align}
&  \mathbb{E}_{t}\left[  \hat{\theta}^{t}\left(  t+\varepsilon,X^{\varepsilon
}\left(  t+\varepsilon\right)  ,\hat{X}^{x_{0}}\left(  t+\varepsilon\right)
\right)  -\hat{\theta}^{t}\left(  t,X^{\varepsilon}\left(  t\right)  ,\hat
{X}^{x_{0}}\left(  t\right)  \right)  \right] \nonumber\\
&  =\int\limits_{t}^{t+\varepsilon}\mathbb{E}_{t}\left[  \hat{\theta}_{s}%
^{t}\left(  s,X^{\varepsilon}\left(  s\right)  ,\hat{X}^{x_{0}}\left(
s\right)  \right)  +\left\langle \hat{\theta}_{x}^{t}\left(  s,X^{\varepsilon
}\left(  s\right)  ,\hat{X}^{x_{0}}\left(  s\right)  \right)  ,\mu\left(
s,X^{\varepsilon}\left(  s\right)  ,v\right)  \right\rangle \right.
\nonumber\\
&  +\left\langle \hat{\theta}_{z}^{t}\left(  s,X^{\varepsilon}\left(
s\right)  ,\hat{X}^{x_{0}}\left(  s\right)  \right)  ,\mu\left(  s,\hat
{X}^{x_{0}}\left(  s\right)  ,\hat{\varphi}^{s}\left(  \hat{X}^{x_{0}}\left(
s\right)  \right)  \right)  \right\rangle \nonumber\\
&  +\frac{1}{2}\text{\textbf{tr}}\left[  \sigma\sigma^{\top}\left(
s,X^{\varepsilon}\left(  s\right)  ,v\right)  \hat{\theta}_{xx}^{t}\left(
s,X^{\varepsilon}\left(  s\right)  ,\hat{X}^{x_{0}}\left(  s\right)  \right)
\right] \nonumber\\
&  +\frac{1}{2}\text{\textbf{tr}}\left[  \sigma\sigma^{\top}\left(  s,\hat
{X}^{x_{0}}\left(  s\right)  ,\hat{\varphi}^{s}\left(  \hat{X}^{x_{0}}\left(
s\right)  \right)  \right)  \hat{\theta}_{zz}^{t}\left(  s,X^{\varepsilon
}\left(  s\right)  ,\hat{X}^{x_{0}}\left(  s\right)  \right)  \right]
\nonumber\\
&  +\text{\textbf{tr}}\left[  \sigma\left(  s,X^{\varepsilon}\left(  s\right)
,v\right)  ^{\top}\hat{\theta}_{xz}^{t}\left(  s,X^{\varepsilon}\left(
s\right)  ,\hat{X}^{x_{0}}\left(  s\right)  \right)  \sigma\left(  s,\hat
{X}^{x_{0}}\left(  s\right)  ,\hat{\varphi}^{s}\left(  \hat{X}^{x_{0}}\left(
s\right)  \right)  \right)  \right] \nonumber\\
&  +%
{\displaystyle\int\limits_{E}}
\left\{  \hat{\theta}^{t}\left(  s,X^{\varepsilon}\left(  s\right)  +c\left(
s,X^{\varepsilon}\left(  s\right)  ,v,e\right)  ,\hat{X}^{x_{0}}\left(
s\right)  +c\left(  s,\hat{X}^{x_{0}}\left(  s\right)  ,\hat{\varphi}%
^{s}\left(  \hat{X}^{x_{0}}\left(  s\right)  \right)  ,e\right)  \right)
\right. \nonumber\\
&  -\hat{\theta}^{t}\left(  s,X^{\varepsilon}\left(  s\right)  ,\hat{X}%
^{x_{0}}\left(  s\right)  \right)  -\left\langle \hat{\theta}_{x}^{t}\left(
s,X^{\varepsilon}\left(  s\right)  ,\hat{X}^{x_{0}}\left(  s\right)  \right)
,c\left(  s,X^{\varepsilon}\left(  s\right)  ,v,e\right)  \right\rangle
\nonumber\\
&  \left.  \left.  -\left\langle \hat{\theta}_{z}^{t}\left(  s,X^{\varepsilon
}\left(  s\right)  ,\hat{X}^{x_{0}}\left(  s\right)  \right)  ,c\left(
s,\hat{X}^{x_{0}}\left(  s\right)  ,\hat{\varphi}^{s}\left(  \hat{X}^{x_{0}%
}\left(  s\right)  \right)  \right)  \right\rangle \right\}  \vartheta\left(
de\right)  \right]  ds. \label{Eq29}%
\end{align}

On the other hand, it follows from the PDE (\ref{Eq19}) that%
\begin{align*}
&  \hat{\theta}_{s}^{t}\left(  s,X^{\varepsilon}\left(  s\right)  ,\hat
{X}^{x_{0}}\left(  s\right)  \right) \\
&  =-f\left(  t,s,X^{\varepsilon}\left(  s\right)  ,\hat{\varphi}^{s}\left(
\hat{X}^{x_{0}}\left(  s\right)  \right)  \right)  -\left\langle \hat{\theta
}_{x}^{t}\left(  s,X^{\varepsilon}\left(  s\right)  ,\hat{X}^{x_{0}}\left(
s\right)  \right)  ,\mu\left(  s,X^{\varepsilon}\left(  s\right)
,\hat{\varphi}^{s}\left(  \hat{X}^{x_{0}}\left(  s\right)  \right)  \right)
\right\rangle \ \\
&  -\left\langle \hat{\theta}_{z}^{t}\left(  s,X^{\varepsilon}\left(
s\right)  ,\hat{X}^{x_{0}}\left(  s\right)  \right)  ,\mu\left(  s,\hat
{X}^{x_{0}}\left(  s\right)  ,\hat{\varphi}^{s}\left(  \hat{X}^{x_{0}}\left(
s\right)  \right)  \right)  \right\rangle \\
&  -\frac{1}{2}\text{\textbf{tr}}\left[  \sigma\sigma^{\top}\left(
s,X^{\varepsilon}\left(  s\right)  ,\hat{\varphi}^{s}\left(  \hat{X}^{x_{0}%
}\left(  s\right)  \right)  \right)  \hat{\theta}_{xx}^{t}\left(
s,X^{\varepsilon}\left(  s\right)  ,\hat{X}^{x_{0}}\left(  s\right)  \right)
\right] \\
&  -\frac{1}{2}\text{\textbf{tr}}\left[  \sigma\sigma^{\top}\left(  s,\hat
{X}^{x_{0}}\left(  s\right)  ,\hat{\varphi}^{s}\left(  \hat{X}^{x_{0}}\left(
s\right)  \right)  \right)  ^{\top}\hat{\theta}_{zz}^{t}\left(
s,X^{\varepsilon}\left(  s\right)  ,\hat{X}^{x_{0}}\left(  s\right)  \right)
\right] \\
&  -\text{\textbf{tr}}\left[  \sigma\left(  s,X^{\varepsilon}\left(  s\right)
,\hat{\varphi}^{s}\left(  \hat{X}^{x_{0}}\left(  s\right)  \right)  \right)
^{\top}\hat{\theta}_{xz}^{t}\left(  s,X^{\varepsilon}\left(  s\right)
,\hat{X}^{x_{0}}\left(  s\right)  \right)  \sigma\left(  s,\hat{X}^{x_{0}%
}\left(  s\right)  ,\hat{\varphi}^{s}\left(  \hat{X}^{x_{0}}\left(  s\right)
\right)  \right)  \right]
\end{align*}%
\begin{align*}
&  -%
{\displaystyle\int\limits_{E}}
\left\{  \hat{\theta}^{t}\left(  s,X^{\varepsilon}\left(  s\right)  +c\left(
s,X^{\varepsilon}\left(  s\right)  ,\hat{\varphi}^{s}\left(  \hat{X}^{x_{0}%
}\left(  s\right)  \right)  ,e\right)  ,\hat{X}^{x_{0}}\left(  s\right)
+c\left(  s,\hat{X}^{x_{0}}\left(  s\right)  ,\hat{\varphi}^{s}\left(  \hat
{X}^{x_{0}}\left(  s\right)  \right)  ,e\right)  \right)  \right. \\
&  -\hat{\theta}^{t}\left(  s,X^{\varepsilon}\left(  s\right)  ,\hat{X}%
^{x_{0}}\left(  s\right)  \right)  -\left\langle \theta_{x}^{t}\left(
s,X^{\varepsilon}\left(  s\right)  ,\hat{X}^{x_{0}}\left(  s\right)  \right)
,c\left(  s,X^{\varepsilon}\left(  s\right)  ,\hat{\varphi}^{s}\left(  \hat
{X}^{x_{0}}\left(  s\right)  \right)  ,e\right)  \right\rangle \\
&  \left.  -\left\langle \theta_{z}^{t}\left(  s,X^{\varepsilon}\left(
s\right)  ,\hat{X}^{x_{0}}\left(  s\right)  \right)  ,c\left(  s,\hat
{X}^{x_{0}}\left(  s\right)  ,\hat{\varphi}^{s}\left(  \hat{X}^{x_{0}}\left(
s\right)  \right)  \right)  \right\rangle \right\}  \vartheta\left(
de\right)  .
\end{align*}
Invoking this into (\ref{Eq29}), we get%
\begin{align}
&  \mathbb{E}_{t}\left[  \hat{\theta}^{t}\left(  t+\varepsilon,X^{\varepsilon
}\left(  t+\varepsilon\right)  ,\hat{X}^{x_{0}}\left(  t+\varepsilon\right)
\right)  -\hat{\theta}^{t}\left(  t,\hat{X}^{x_{0}}\left(  t\right)  ,\hat
{X}^{x_{0}}\left(  t\right)  \right)  \right] \nonumber\\
&  =\int\limits_{t}^{t+\varepsilon}\mathbb{E}_{t}\left[  \left\{  \left\langle
\hat{\theta}_{x}^{t}\left(  s,X^{\varepsilon}\left(  s\right)  ,\hat{X}%
^{x_{0}}\left(  s\right)  \right)  ,\mu\left(  s,X^{\varepsilon}\left(
s\right)  ,v\right)  -\mu\left(  s,X^{\varepsilon}\left(  s\right)
,\hat{\varphi}^{s}\left(  \hat{X}^{x_{0}}\left(  s\right)  \right)  \right)
\right\rangle \right.  \right. \nonumber\\
&  +\frac{1}{2}\text{\textbf{tr}}\left[  \left(  \sigma\sigma^{\top}\left(
s,X^{\varepsilon}\left(  s\right)  ,v\right)  -\sigma\sigma^{\top}\left(
s,X^{\varepsilon}\left(  s\right)  ,\hat{\varphi}^{s}\left(  \hat{X}^{x_{0}%
}\left(  s\right)  \right)  \right)  \right)  \hat{\theta}_{xx}^{t}\left(
s,X^{\varepsilon}\left(  s\right)  ,\hat{X}^{x_{0}}\left(  s\right)  \right)
\right] \nonumber\\
&  +\text{\textbf{tr}}\left[  \left(  \sigma\left(  s,X^{\varepsilon}\left(
s\right)  ,v\right)  -\sigma\left(  s,X^{\varepsilon}\left(  s\right)
,\hat{\varphi}^{s}\left(  \hat{X}^{x_{0}}\left(  s\right)  \right)  \right)
\right)  ^{\top}\hat{\theta}_{xz}^{t}\left(  s,X^{\varepsilon}\left(
s\right)  ,\hat{X}^{x_{0}}\left(  s\right)  \right)  \right. \nonumber\\
&  \left.  \times\sigma\left(  s,\hat{X}^{x_{0}}\left(  s\right)
,\hat{\varphi}^{s}\left(  \hat{X}^{x_{0}}\left(  s\right)  \right)  \right)
\right]  -f\left(  t,s,X^{\varepsilon}\left(  s\right)  ,\hat{\varphi}%
^{s}\left(  \hat{X}^{x_{0}}\left(  s\right)  \right)  \right) \nonumber\\
&  +%
{\displaystyle\int\limits_{E}}
\left\{  \hat{\theta}^{t}\left(  s,X^{\varepsilon}\left(  s\right)  +c\left(
s,X^{\varepsilon}\left(  s\right)  ,v,e\right)  ,\hat{X}^{x_{0}}\left(
s\right)  +c\left(  s,\hat{X}^{x_{0}}\left(  s\right)  ,\hat{\varphi}%
^{s}\left(  \hat{X}^{x_{0}}\left(  s\right)  \right)  ,e\right)  \right)
\right. \nonumber\\
&  -\hat{\theta}^{t}\left(  s,X^{\varepsilon}\left(  s\right)  +c\left(
s,X^{\varepsilon}\left(  s\right)  ,\hat{\varphi}^{s}\left(  \hat{X}^{x_{0}%
}\left(  s\right)  \right)  ,e\right)  ,\hat{X}^{x_{0}}\left(  s\right)
+c\left(  s,\hat{X}^{x_{0}}\left(  s\right)  ,\hat{\varphi}^{s}\left(  \hat
{X}^{x_{0}}\left(  s\right)  \right)  ,e\right)  \right) \nonumber\\
&  \left.  \left.  -\left\langle \theta_{x}^{t}\left(  s,X^{\varepsilon
}\left(  s\right)  ,\hat{X}^{x_{0}}\left(  s\right)  \right)  ,c\left(
s,X^{\varepsilon}\left(  s\right)  ,v,e\right)  -c\left(  s,X^{\varepsilon
}\left(  s\right)  ,\hat{\varphi}^{s}\left(  \hat{X}^{x_{0}}\left(  s\right)
\right)  ,e\right)  \right\rangle \right\}  \vartheta\left(  de\right)
\right]  ds\text{.} \label{Eq30}%
\end{align}

Next, for each $1\leq i\leq m$, applying It\^{o}'s formula to $\hat{g}%
^{i}\left(  \cdot,X^{\varepsilon}\left(  \cdot\right)  ,\hat{X}^{x_{0}}\left(
\cdot\right)  \right)  $ on time interval $\left[  t,t+\varepsilon\right]  $,
and taking the conditional expectations, we get%
\begin{align*}
&  d\mathbb{E}_{t}\left[  \hat{g}^{i}\left(  s,X^{\varepsilon}\left(
s\right)  ,\hat{X}^{x_{0}}\left(  s\right)  \right)  \right] \\
&  =\mathbb{E}_{t}\left[  \hat{g}_{s}^{i}\left(  s,X^{\varepsilon}\left(
s\right)  ,\hat{X}^{x_{0}}\left(  s\right)  \right)  +\left\langle \hat{g}%
_{x}^{i}\left(  s,X^{\varepsilon}\left(  s\right)  ,\hat{X}^{x_{0}}\left(
s\right)  \right)  ,\mu\left(  s,X^{\varepsilon}\left(  s\right)  ,v\right)
\right\rangle \right.  \ \\
&  +\left\langle \hat{g}_{z}^{i}\left(  s,X^{\varepsilon}\left(  s\right)
,\hat{X}^{x_{0}}\left(  s\right)  \right)  ,\mu\left(  s,\hat{X}^{x_{0}%
}\left(  s\right)  ,\hat{\varphi}^{s}\left(  \hat{X}^{x_{0}}\left(  s\right)
\right)  \right)  \right\rangle \\
&  +\frac{1}{2}\text{\textbf{tr}}\left[  \sigma\sigma^{\top}\left(
s,X^{\varepsilon}\left(  s\right)  ,v\right)  \hat{g}_{xx}^{i}\left(
s,X^{\varepsilon}\left(  s\right)  ,\hat{X}^{x_{0}}\left(  s\right)  \right)
\right] \\
&  +\frac{1}{2}\text{\textbf{tr}}\left[  \sigma\sigma^{\top}\left(  s,\hat
{X}^{x_{0}}\left(  s\right)  ,\hat{\varphi}^{s}\left(  \hat{X}^{x_{0}}\left(
s\right)  \right)  \right)  \hat{g}_{zz}^{i}\left(  s,X^{\varepsilon}\left(
s\right)  ,\hat{X}^{x_{0}}\left(  s\right)  \right)  \right] \\
&  +\text{\textbf{tr}}\left[  \sigma\left(  s,X^{\varepsilon}\left(  s\right)
,v\right)  ^{\top}\hat{g}_{xz}^{i}\left(  s,X^{\varepsilon}\left(  s\right)
,\hat{X}^{x_{0}}\left(  s\right)  \right)  \sigma\left(  s,\hat{X}^{x_{0}%
}\left(  s\right)  ,\hat{\varphi}^{s}\left(  \hat{X}^{x_{0}}\left(  s\right)
\right)  \right)  \right] \\
&  +%
{\displaystyle\int\limits_{E}}
\left\{  \hat{g}^{i}\left(  s,X^{\varepsilon}\left(  s\right)  +c\left(
s,X^{\varepsilon}\left(  s\right)  ,v,e\right)  ,\hat{X}^{x_{0}}\left(
s\right)  +c\left(  s,\hat{X}^{x_{0}}\left(  s\right)  ,\hat{\varphi}%
^{s}\left(  \hat{X}^{x_{0}}\left(  s\right)  \right)  ,e\right)  \right)
\right. \\
&  -\text{ }\hat{g}^{i}\left(  s,X^{\varepsilon}\left(  s\right)  ,\hat
{X}^{x_{0}}\left(  s\right)  \right)  -\left\langle \hat{g}_{x}^{i}\left(
s,X^{\varepsilon}\left(  s\right)  ,\hat{X}^{x_{0}}\left(  s\right)  \right)
,c\left(  s,X^{\varepsilon}\left(  s\right)  ,v,e\right)  \right\rangle \\
&  \left.  \left.  -\left\langle \hat{g}_{z}^{i}\left(  s,X^{\varepsilon
}\left(  s\right)  ,\hat{X}^{x_{0}}\left(  s\right)  \right)  ,c\left(
s,\hat{X}^{x_{0}}\left(  s\right)  ,\hat{\varphi}^{s}\left(  \hat{X}^{x_{0}%
}\left(  s\right)  \right)  \right)  \right\rangle \right\}  \vartheta\left(
de\right)  \right]  ds.
\end{align*}
Thus, using the system of PDEs (\ref{Eq21}) yields%
\begin{align*}
&  d\mathbb{E}_{t}\left[  \hat{g}^{i}\left(  s,X^{\varepsilon}\left(
s\right)  ,\hat{X}^{x_{0}}\left(  s\right)  \right)  \right] \\
&  =\mathbb{E}_{t}\left[  \left\langle \hat{g}_{x}^{i}\left(  s,X^{\varepsilon
}\left(  s\right)  ,\hat{X}^{x_{0}}\left(  s\right)  \right)  ,\mu\left(
s,X^{\varepsilon}\left(  s\right)  ,v\right)  -\mu\left(  s,X^{\varepsilon
}\left(  s\right)  ,\hat{\varphi}^{s}\left(  \hat{X}^{x_{0}}\left(  s\right)
\right)  \right)  \right\rangle \right. \\
&  +\frac{1}{2}\text{\textbf{tr}}\left[  \left(  \sigma\sigma^{\top}\left(
s,X^{\varepsilon}\left(  s\right)  ,v\right)  -\sigma\sigma^{\top}\left(
s,X^{\varepsilon}\left(  s\right)  ,\hat{\varphi}^{s}\left(  \hat{X}^{x_{0}%
}\left(  s\right)  \right)  \right)  \right)  \hat{g}_{xx}^{i}\left(
s,X^{\varepsilon}\left(  s\right)  ,\hat{X}^{x_{0}}\left(  s\right)  \right)
\right]
\end{align*}%
\begin{align*}
&  +\text{\textbf{tr}}\left[  \left(  \sigma\left(  s,X^{\varepsilon}\left(
s\right)  ,v\right)  -\sigma\left(  s,X^{\varepsilon}\left(  s\right)
,\hat{\varphi}^{s}\left(  \hat{X}^{x_{0}}\left(  s\right)  \right)  \right)
\right)  ^{\top}\hat{g}_{xz}^{i}\left(  s,X^{\varepsilon}\left(  s\right)
,\hat{X}^{x_{0}}\left(  s\right)  \right)  \right. \\
&  \left.  \times\sigma\left(  s,\hat{X}^{x_{0}}\left(  s\right)
,\hat{\varphi}^{s}\left(  \hat{X}^{x_{0}}\left(  s\right)  \right)  \right)
\right] \\
&  +%
{\displaystyle\int\limits_{E}}
\left\{  \hat{g}^{i}\left(  s,X^{\varepsilon}\left(  s\right)  +c\left(
s,X^{\varepsilon}\left(  s\right)  ,v,e\right)  ,\hat{X}^{x_{0}}\left(
s\right)  +c\left(  s,\hat{X}^{x_{0}}\left(  s\right)  ,\hat{\varphi}%
^{s}\left(  \hat{X}^{x_{0}}\left(  s\right)  \right)  ,e\right)  \right)
\right. \\
&  -\hat{g}^{i}\left(  s,X^{\varepsilon}\left(  s\right)  +c\left(
s,X^{\varepsilon}\left(  s\right)  ,\hat{\varphi}^{s}\left(  \hat{X}^{x_{0}%
}\left(  s\right)  \right)  ,e\right)  ,\hat{X}^{x_{0}}\left(  s\right)
+c\left(  s,\hat{X}^{x_{0}}\left(  s\right)  ,\hat{\varphi}^{s}\left(  \hat
{X}^{x_{0}}\left(  s\right)  \right)  ,e\right)  \right) \\
&  \left.  \left.  -\left\langle \hat{g}_{x}^{i}\left(  s,X^{\varepsilon
}\left(  s\right)  ,\hat{X}^{x_{0}}\left(  s\right)  \right)  ,c\left(
s,X^{\varepsilon}\left(  s\right)  ,v,e\right)  -c\left(  s,X^{\varepsilon
}\left(  s\right)  ,\hat{\varphi}^{s}\left(  \hat{X}^{x_{0}}\left(  s\right)
\right)  ,e\right)  \right\rangle \right\}  \vartheta\left(  de\right)
\right]  ds.
\end{align*}

By the chain rule we obtain that%
\begin{align}
&  G\left(  t,\hat{X}^{x_{0}}\left(  t\right)  ,\mathbb{E}_{t}\left[  \hat
{g}\left(  t+\varepsilon,X^{\varepsilon}\left(  t+\varepsilon\right)  ,\hat
{X}^{x_{0}}\left(  t+\varepsilon\right)  \right)  \right]  \right)  -G\left(
t,\hat{X}^{x_{0}}\left(  t\right)  ,\hat{g}\left(  t,X^{\varepsilon}\left(
t\right)  ,\hat{X}^{x_{0}}\left(  t\right)  \right)  \right) \nonumber\\
&  =\int_{t}^{t+\varepsilon}\sum_{i=1}^{m}G_{\bar{x}_{i}}\left(  t,\hat
{X}^{x_{0}}\left(  t\right)  ,\mathbb{E}_{t}\left[  \hat{g}\left(
s,X^{\varepsilon}\left(  s\right)  ,\hat{X}^{x_{0}}\left(  s\right)  \right)
\right]  \right)  d\mathbb{E}_{t}\left[  \hat{g}^{i}\left(  s,X^{\varepsilon
}\left(  s\right)  ,\hat{X}^{x_{0}}\left(  s\right)  \right)  \right]
\text{.} \label{Eq31}%
\end{align}

Combining (\ref{Eq28}) together with (\ref{Eq30}) and (\ref{Eq31}), it follows
that%
\begin{align*}
&  \mathbf{J}\left(  t,\hat{X}^{x_{0}}\left(  t\right)  ;u^{\varepsilon
}\left(  \cdot\right)  \right)  -\mathbf{J}\left(  t,\hat{X}^{x_{0}}\left(
t\right)  ;\hat{\varphi}\left(  \cdot,\cdot\right)  \right) \\
&  =\left[  \int\limits_{t}^{t+\varepsilon}\mathbb{E}_{t}\left[  \left\langle
\hat{\theta}_{x}^{t}\left(  s,X^{\varepsilon}\left(  s\right)  ,\hat{X}%
^{x_{0}}\left(  s\right)  \right)  ,\mu\left(  s,X^{\varepsilon}\left(
s\right)  ,v\right)  -\mu\left(  s,X^{\varepsilon}\left(  s\right)
,\hat{\varphi}^{s}\left(  \hat{X}^{x_{0}}\left(  s\right)  \right)  \right)
\right\rangle \right.  \right. \\
&  +\frac{1}{2}\text{\textbf{tr}}\left[  \left(  \sigma\sigma^{\top}\left(
s,X^{\varepsilon}\left(  s\right)  ,v\right)  -\sigma\sigma^{\top}\left(
s,X^{\varepsilon}\left(  s\right)  ,\hat{\varphi}^{s}\left(  \hat{X}^{x_{0}%
}\left(  s\right)  \right)  \right)  \right)  \hat{\theta}_{xx}^{t}\left(
s,X^{\varepsilon}\left(  s\right)  ,\hat{X}^{x_{0}}\left(  s\right)  \right)
\right] \\
&  +\text{\textbf{tr}}\left[  \left(  \sigma\left(  s,X^{\varepsilon}\left(
s\right)  ,v\right)  -\sigma\left(  s,X^{\varepsilon}\left(  s\right)
,\hat{\varphi}^{s}\left(  \hat{X}^{x_{0}}\left(  s\right)  \right)  \right)
\right)  \hat{\theta}_{xz}^{t}\left(  s,X^{\varepsilon}\left(  s\right)
,\hat{X}^{x_{0}}\left(  s\right)  \right)  \right. \\
&  \left.  \times\sigma\left(  s,\hat{X}^{x_{0}}\left(  s\right)
,\hat{\varphi}^{s}\left(  \hat{X}^{x_{0}}\left(  s\right)  \right)  \right)
\right] \\
&  +%
{\displaystyle\int\limits_{E}}
\left\{  \hat{\theta}^{t}\left(  s,X^{\varepsilon}\left(  s\right)  +c\left(
s,X^{\varepsilon}\left(  s\right)  ,v,e\right)  ,\hat{X}^{x_{0}}\left(
s\right)  +c\left(  s,\hat{X}^{x_{0}}\left(  s\right)  ,\hat{\varphi}%
^{s}\left(  \hat{X}^{x_{0}}\left(  s\right)  \right)  ,e\right)  \right)
\right. \\
&  -\hat{\theta}^{t}\left(  s,X^{\varepsilon}\left(  s\right)  +c\left(
s,X^{\varepsilon}\left(  s\right)  ,\hat{\varphi}^{s}\left(  \hat{X}^{x_{0}%
}\left(  s\right)  \right)  ,e\right)  ,\hat{X}^{x_{0}}\left(  s\right)
+c\left(  s,\hat{X}^{x_{0}}\left(  s\right)  ,\hat{\varphi}^{s}\left(  \hat
{X}^{x_{0}}\left(  s\right)  \right)  ,e\right)  \right) \\
&  \left.  -\left\langle \hat{\theta}_{x}^{t}\left(  s,X^{\varepsilon}\left(
s\right)  ,\hat{X}^{x_{0}}\left(  s\right)  \right)  ,c\left(
s,X^{\varepsilon}\left(  s\right)  ,v,e\right)  -c\left(  s,X^{\varepsilon
}\left(  s\right)  ,\hat{\varphi}^{s}\left(  \hat{X}^{x_{0}}\left(  s\right)
\right)  ,e\right)  \right\rangle \right\}  \vartheta\left(  de\right) \\
&  +f\left(  t,s,X^{\varepsilon}\left(  s\right)  ,v\right)  -f\left(
t,s,X^{\varepsilon}\left(  s\right)  ,\hat{\varphi}^{s}\left(  \hat{X}^{x_{0}%
}\left(  s\right)  \right)  \right) \\
&  +\sum_{i=1}^{m}G_{\bar{x}_{i}}\left(  t,\hat{X}^{x_{0}}\left(  t\right)
,\mathbb{E}_{t}\left[  \hat{g}\left(  s,X^{\varepsilon}\left(  s\right)
,\hat{X}^{x_{0}}\left(  s\right)  \right)  \right]  \right) \\
&  \times\left(  \left\langle \hat{g}_{x}^{i}\left(  s,X^{\varepsilon}\left(
s\right)  ,\hat{X}^{x_{0}}\left(  s\right)  \right)  ,\mu\left(
s,X^{\varepsilon}\left(  s\right)  ,v\right)  -\mu\left(  s,X^{\varepsilon
}\left(  s\right)  ,\hat{\varphi}^{s}\left(  \hat{X}^{x_{0}}\left(  s\right)
\right)  \right)  \right\rangle \right. \\
&  +\frac{1}{2}\text{\textbf{tr}}\left[  \left(  \sigma\sigma^{\top}\left(
s,X^{\varepsilon}\left(  s\right)  ,v\right)  -\sigma\sigma^{\top}\left(
s,X^{\varepsilon},\hat{\varphi}^{s}\left(  \hat{X}^{x_{0}}\left(  s\right)
\right)  \right)  \right)  \hat{g}_{xx}^{i}\left(  s,X^{\varepsilon}\left(
s\right)  ,\hat{X}^{x_{0}}\left(  s\right)  \right)  \right] \\
&  +\text{\textbf{tr}}\left[  \left(  \sigma\left(  s,X^{\varepsilon}\left(
s\right)  ,v\right)  -\sigma\left(  s,X^{\varepsilon}\left(  s\right)
,\hat{\varphi}^{s}\left(  \hat{X}^{x_{0}}\left(  s\right)  \right)  \right)
\right)  ^{\top}\hat{g}_{xz}^{i}\left(  s,X^{\varepsilon}\left(  s\right)
,\hat{X}^{x_{0}}\left(  s\right)  \right)  \right. \\
&  \left.  \times\sigma\left(  s,\hat{X}^{x_{0}}\left(  s\right)
,\hat{\varphi}^{s}\left(  \hat{X}^{x_{0}}\left(  s\right)  \right)  \right)
\right] \\
&  +%
{\displaystyle\int\limits_{E}}
\left\{  \hat{g}^{i}\left(  s,X^{\varepsilon}\left(  s\right)  +c\left(
s,X^{\varepsilon}\left(  s\right)  ,v,e\right)  ,\hat{X}^{x_{0}}\left(
s\right)  +c\left(  s,\hat{X}^{x_{0}}\left(  s\right)  ,\hat{\varphi}%
^{s}\left(  \hat{X}^{x_{0}}\left(  s\right)  \right)  ,e\right)  \right)
\right. \\
&  -\hat{g}^{i}\left(  s,X^{\varepsilon}\left(  s\right)  +c\left(
s,X^{\varepsilon}\left(  s\right)  ,\hat{\varphi}^{s}\left(  \hat{X}^{x_{0}%
}\left(  s\right)  \right)  ,e\right)  ,\hat{X}^{x_{0}}\left(  s\right)
+c\left(  s,\hat{X}^{x_{0}}\left(  s\right)  ,\hat{\varphi}^{s}\left(  \hat
{X}^{x_{0}}\left(  s\right)  \right)  ,e\right)  \right) \\
&  \left.  \left.  \left.  -\left\langle \hat{g}_{x}^{i}\left(
s,X^{\varepsilon}\left(  s\right)  ,\hat{X}^{x_{0}}\left(  s\right)  \right)
,c\left(  s,X^{\varepsilon}\left(  s\right)  ,v,e\right)  -c\left(
s,X^{\varepsilon}\left(  s\right)  ,\hat{\varphi}^{s}\left(  \hat{X}^{x_{0}%
}\left(  s\right)  \right)  ,e\right)  \right\rangle \right\}  \vartheta
\left(  de\right)  \right)  \right]  ds.
\end{align*}
This completes the proof.\eop

\begin{remark}
For each $t\in\left[  0,T\right]  $ fixed, setting $\left(  s,\mathbf{x}%
,\mathbf{z}\right)  =\left(  t,\hat{X}^{x_{0}}\left(  t\right)  ,\hat
{X}^{x_{0}}\left(  t\right)  \right)  $ in (\ref{Eq25*})-(\ref{Eq26*}), we
obtain that%
\begin{align}
&  \hat{\theta}^{t}\left(  t,\hat{X}^{x_{0}}\left(  t\right)  ,\hat{X}^{x_{0}%
}\left(  t\right)  \right) \nonumber\\
&  =\mathbb{E}_{t}\left[  \int_{t}^{T}f\left(  t,s,X^{t,\hat{X}^{x_{0}}\left(
t\right)  ,\hat{X}^{x_{0}}\left(  t\right)  }\left(  s\right)  ,\hat{\varphi
}^{s}\left(  \hat{X}^{t,\hat{X}^{x_{0}}\left(  t\right)  }\left(  s\right)
\right)  \right)  ds+F\left(  t,X^{t,\hat{X}^{x_{0}}\left(  t\right)  ,\hat
{X}^{x_{0}}\left(  t\right)  }\left(  T\right)  \right)  \right] \nonumber\\
&  =\mathbb{E}_{t}\left[  \int_{t}^{T}f\left(  t,s,\hat{X}^{x_{0}}\left(
s\right)  ,\hat{\varphi}^{s}\left(  \hat{X}^{x_{0}}\left(  s\right)  \right)
\right)  ds+F\left(  t,\hat{X}^{x_{0}}\left(  T\right)  \right)  \right]
\label{Eq32}%
\end{align}
and%
\begin{equation}
\hat{g}\left(  t,\hat{X}^{x_{0}}\left(  t\right)  ,\hat{X}^{x_{0}}\left(
t\right)  \right)  =\mathbb{E}_{t}\left[  \Psi\left(  X^{t,\hat{X}^{x_{0}%
}\left(  t\right)  ,\hat{X}^{x_{0}}\left(  t\right)  }\left(  T\right)
\right)  \right]  =\mathbb{E}_{t}\left[  \Psi\left(  \hat{X}^{x_{0}}\left(
T\right)  \right)  \right]  \text{,} \label{Eq33}%
\end{equation}
where we have used Equality (\ref{Eq10}).
\end{remark}

\subsection{Verification Theorem}

Now, we are ready to state the main result of this work.

\begin{theorem}
[Verification Theorem]\label{result2}Let \textbf{(H1)-(H2)} hold. Given an
admissible strategy $\hat{\varphi}\left(  \cdot,\cdot\right)  \in\mathcal{S}$,
suppose that for each $t\in\left[  0,T\right]  $, the IPDEs (\ref{Eq19}) and
(\ref{Eq21}) admit the classical solutions $\hat{\theta}^{t}\left(
\cdot,\cdot,\cdot\right)  \in\mathcal{C}^{1,2,2}\left(  \left[  t,T\right]
\times%
\mathbb{R}
^{n}\times%
\mathbb{R}
^{n};%
\mathbb{R}
\right)  $ and $\hat{g}\left(  \cdot,\cdot,\cdot\right)  \in\mathcal{C}%
^{1,2,2}\left(  \left[  0,T\right]  \times%
\mathbb{R}
^{n}\times%
\mathbb{R}
^{n};%
\mathbb{R}
\right)  ^{m}$, respectively, such that the following hold:\textbf{\newline(i)
}$\hat{\varphi}\left(  \cdot,\cdot\right)  $ is a continuous function and
there exists a constant $C>0$ such that for all $\left(  s,z\right)
\in\left[  0,T\right]  \times\mathbb{R}^{n}$ we have%
\begin{equation}
\left\vert \hat{\varphi}^{s}\left(  z\right)  \right\vert \leq C\left(
1+\left\vert z\right\vert \right)  . \label{Eq34}%
\end{equation}
\newline\textbf{(ii) }$\hat{\theta}^{t}\left(  \cdot,\cdot,\cdot\right)
\in\mathcal{C}_{2}^{1,2,2}\left(  \left[  t,T\right]  \times%
\mathbb{R}
^{n}\times%
\mathbb{R}
^{n};%
\mathbb{R}
\right)  $ and $\hat{g}\left(  \cdot,\cdot,\cdot\right)  \in\mathcal{C}%
_{2}^{1,2,2}\left(  \left[  0,T\right]  \times%
\mathbb{R}
^{n}\times%
\mathbb{R}
^{n};%
\mathbb{R}
\right)  ^{m}$.\newline\textbf{(iii) }For all $\left(  t,z,u\right)
\in\left[  0,T\right]  \times%
\mathbb{R}
^{n}\times U$,%
\begin{equation}
0\leq\mathcal{H}\left(  t,z,t,z,z,u\right)  -\mathcal{H}\left(  t,z,t,z,z,\hat
{\varphi}^{t}\left(  z\right)  \right)  \text{.} \label{Eq35}%
\end{equation}
Then $\hat{\varphi}\left(  \cdot,\cdot\right)  $ is an open-loop equilibrium
strategy. Furthermore, the objective value of $\hat{\varphi}\left(
\cdot,\cdot\right)  $ at time $t\in\left[  0,T\right]  $ is given by%
\begin{equation}
\mathbf{J}\left(  t,\hat{X}^{x_{0}}\left(  t\right)  ;\hat{\varphi}\left(
\cdot,\cdot\right)  \right)  =\theta^{t}\left(  t,\hat{X}^{x_{0}}\left(
t\right)  ,\hat{X}^{x_{0}}\left(  t\right)  \right)  +G\left(  t,\hat
{X}^{x_{0}}\left(  t\right)  ,g\left(  t,\hat{X}^{x_{0}}\left(  t\right)
,\hat{X}^{x_{0}}\left(  t\right)  \right)  \right)  \text{.} \label{Eq36}%
\end{equation}

\end{theorem}

\bop In this proof, $C>0$ denotes a universal constant which may vary from
line to line. Let $\hat{\varphi}\left(  \cdot,\cdot\right)  \in\mathcal{S}$ be
an admissible strategy for which Conditions (i)-(iii) in Theorem \ref{result2}
hold. We are going to show that $\hat{\varphi}\left(  \cdot,\cdot\right)  $ is
an open-loop equilibrium strategy. Consider the perturbed strategy
$u^{\varepsilon}\left(  \cdot\right)  $ defined by the spike variation
(\ref{Eq6}) for some fixed arbitrarily $t\in\left[  0,T\right]  $,
$v\in\mathbb{L}^{4}\left(  \Omega,\mathcal{F}_{t},\mathbb{P};U\right)  $ and
$\varepsilon\in\left[  0,T-t\right)  $. Let $X^{\varepsilon}\left(
\cdot\right)  $ be the unique strong solution of (\ref{Eq8}) and
$X^{t,v}\left(  \cdot\right)  $ be the solution of the following SDE%
\[
\left\{
\begin{array}
[c]{l}%
dX^{t,v}\left(  s\right)  =b\left(  s,X^{t,v}\left(  s\right)  ,v\right)
ds+\sigma\left(  s,X^{t,v}\left(  s\right)  ,v\right)  dW\left(  s\right) \\
\text{ \ \ \ \ \ \ \ \ \ \ \ \ \ \ }+%
{\displaystyle\int_{E}}
c\left(  s,X_{-}^{t,v}\left(  s\right)  ,v,e\right)  \tilde{N}\left(
ds,de\right)  ,\text{ for }s\in\left[  t,T\right]  ,\\
X^{t,v}\left(  t\right)  =\hat{X}^{x_{0}}\left(  t\right)  .
\end{array}
\right.
\]
Under Assumption \textbf{(H1) }the above SDE admits a unique solution
$X^{t,v}\left(  \cdot\right)  $ and there exists a constant $C>0$ such that%
\[
\mathbb{E}\left[  \sup_{s\in\left[  t,T\right]  }\left\vert X^{t,v}\left(
s\right)  \right\vert ^{4}\right]  \leq C\left(  1+\mathbb{E}\left[
\left\vert \hat{X}^{x_{0}}\left(  t\right)  \right\vert ^{4}\right]  \right)
\text{.}%
\]
Moreover, since $u^{\varepsilon}\left(  s\right)  =v$ on $\left[
t,t+\varepsilon\right]  $, we have%
\begin{equation}
X^{\varepsilon}\left(  s\right)  =X^{t,v}\left(  s\right)  ,\text{ a.s., for
}s\in\left[  t,t+\varepsilon\right]  \text{.} \label{Eq37}%
\end{equation}

First, we begin by showing that%
\begin{align*}
&  \underset{\varepsilon\downarrow0}{\lim}\frac{1}{\varepsilon}\left\{
\mathbf{J}\left(  t,\hat{X}^{x_{0}}\left(  t\right)  ;u^{\varepsilon}\left(
\cdot\right)  \right)  -\mathbf{J}\left(  t,\hat{X}^{x_{0}}\left(  t\right)
;\hat{\varphi}\left(  \cdot,\cdot\right)  \right)  \right\} \\
&  =\mathcal{H}\left(  t,\hat{X}^{x_{0}}\left(  t\right)  ,t,\hat{X}^{x_{0}%
}\left(  t\right)  ,\hat{X}^{x_{0}}\left(  t\right)  ,v\right) \\
&  -\mathcal{H}\left(  t,\hat{X}^{x_{0}}\left(  t\right)  ,t,\hat{X}^{x_{0}%
}\left(  t\right)  ,\hat{X}^{x_{0}}\left(  t\right)  ,\hat{\varphi}^{t}\left(
\hat{X}^{x_{0}}\left(  t\right)  \right)  \right)  \text{.}%
\end{align*}
To this end, let us introduce in the time interval $\left[  t,T\right]  $ the
process $\Lambda\left(  \cdot;t\right)  $ defined by%
\[
\Lambda\left(  s;t\right)  :=\mathcal{H}\left(  t,\hat{X}^{x_{0}}\left(
t\right)  ,s,X^{t,v}\left(  s\right)  ,\hat{X}^{x_{0}}\left(  s\right)
,v\right)  -\mathcal{H}\left(  t,\hat{X}^{x_{0}}\left(  t\right)
,s,X^{t,v}\left(  s\right)  ,\hat{X}^{x_{0}}\left(  s\right)  ,\hat{\varphi
}^{s}\left(  \hat{X}^{x_{0}}\left(  s\right)  \right)  \right)  .
\]
Under Assumptions \textbf{(H1)-(H2)} together with Condition (i), it is easy
to see that $\Lambda\left(  \cdot;t\right)  $ is a right-continuous $\left(
\mathcal{F}_{s}\right)  _{s\in\left[  t,T\right]  }$-progressively measurable
process. Moreover, it follows from Assumptions \textbf{(H1)-(H2) }together
with\textbf{ }Condition (ii) and Definition \ref{Def3} that there exists a
constant $C>0$, such that%
\begin{align*}
&  \left\vert \mathcal{H}\left(  t,y,s,X,Z,u\right)  \right\vert \\
&  \leq\left\vert \hat{\theta}_{x}^{t}\left(  s,X,Z\right)  \right\vert
\left\vert \mu\left(  s,X,u\right)  \right\vert +\frac{1}{2}\left\vert
\hat{\theta}_{xx}^{t}\left(  s,X,Z\right)  \right\vert \left\vert
\sigma\left(  s,X,u\right)  \right\vert ^{2}\\
&  +\left\vert \hat{\theta}_{xz}^{t}\left(  s,X,Z\right)  \right\vert
\left\vert \sigma\left(  s,X,u\right)  \right\vert \left\vert \sigma\left(
s,Z,\hat{\varphi}^{s}\left(  Z\right)  \right)  \right\vert +\left\vert
f\left(  t,s,X,u\right)  \right\vert \\
&  +%
{\displaystyle\int_{E}}
\left\vert \hat{\theta}^{t}\left(  s,X+c\left(  s,X,u,e\right)  ,Z+c\left(
s,Z,\hat{\varphi}^{s}\left(  Z\right)  ,e\right)  \right)  \right\vert
\vartheta\left(  de\right) \\
&  +%
{\displaystyle\int_{E}}
\left\vert \hat{\theta}_{x}^{t}\left(  s,X,Z\right)  \right\vert \left\vert
c\left(  s,X,u,e\right)  \right\vert \vartheta\left(  de\right) \\
&  +\sum_{i=1}^{m}\left\vert G_{\bar{x}_{i}}\left(  t,y,\mathbb{E}_{t}\left[
\hat{g}\left(  s,X,Z\right)  \right]  \right)  \right\vert \times\left(
\left\vert \hat{g}_{x}^{i}\left(  s,X,Z\right)  \right\vert \left\vert
\mu\left(  s,X,u\right)  \right\vert \right. \\
&  +\frac{1}{2}\left\vert \hat{g}_{xx}^{i}\left(  s,X,Z\right)  \right\vert
\left\vert \sigma\left(  s,X,u\right)  \right\vert ^{2}+\left\vert g_{xz}%
^{i}\left(  s,X,Z\right)  \right\vert \left\vert \sigma\left(  s,X,u\right)
\right\vert \left\vert \sigma\left(  s,Z,\hat{\varphi}^{s}\left(  Z\right)
\right)  \right\vert \\
&  +%
{\displaystyle\int_{E}}
\left\vert \hat{g}^{i}\left(  s,X+c\left(  s,X,u,e\right)  ,Z+c\left(
s,Z,\hat{\varphi}^{s}\left(  Z\right)  ,e\right)  \right)  \right\vert
\vartheta\left(  de\right) \\
&  \left.  +%
{\displaystyle\int_{E}}
\left\vert \hat{g}_{x}^{i}\left(  s,X,Z\right)  \right\vert \left\vert
c\left(  s,X,u,e\right)  \right\vert \vartheta\left(  de\right)  \right) \\
&  \leq C\left(  1+\left\vert X\right\vert ^{2}+\left\vert u\right\vert
^{2}+\left\vert \varphi^{s}\left(  Z\right)  \right\vert ^{2}+\left\vert
Z\right\vert ^{2}\right) \\
&  +\sum_{i=1}^{m}C\left(  1+\left\vert y\right\vert +\mathbb{E}_{t}\left[
\left\vert X\right\vert ^{2}\right]  +\left\vert u\right\vert ^{2}%
+\mathbb{E}_{t}\left[  \left\vert \varphi^{s}\left(  Z\right)  \right\vert
^{2}\right]  +\mathbb{E}_{t}\left[  \left\vert Z\right\vert ^{2}\right]
\right)  \left(  1+\left\vert X\right\vert ^{2}+\left\vert u\right\vert
^{2}+\left\vert \varphi^{s}\left(  Z\right)  \right\vert ^{2}+\left\vert
Z\right\vert ^{2}\right) \\
&  \leq C\left(  1+\left\vert y\right\vert ^{4}+\left\vert Z\right\vert
^{4}+\mathbb{E}_{t}\left[  \left\vert Z\right\vert ^{4}\right]  +\left\vert
X\right\vert ^{4}+\mathbb{E}_{t}\left[  \left\vert X\right\vert ^{4}\right]
+\mathbb{E}_{t}\left[  \left\vert \hat{\varphi}^{s}\left(  Z\right)
\right\vert ^{4}\right]  +\left\vert \hat{\varphi}^{s}\left(  Z\right)
\right\vert ^{4}+\left\vert u\right\vert ^{4}\right)  \text{,}%
\end{align*}
for all $\left(  t,s\right)  \in D\left[  0,T\right]  ,$ $\left(  u,y\right)
\in U\times%
\mathbb{R}
^{n}$ and for each $X,Z\in\mathbb{L}^{4}\left(  \Omega,\mathcal{F}%
_{s},\mathbb{P};%
\mathbb{R}
^{n}\right)  $, where we have used Jensen inequality together with the
inequalities $a_{1}^{2}\leq\frac{1}{2}\left(  a_{1}^{4}+1\right)  ,$
$a_{1}a_{2}\leq\frac{1}{2}\left(  a_{1}^{2}+b_{2}^{2}\right)  $ and $\left(
{\displaystyle\sum_{i=1}^{n}}
a_{i}^{2}\right)  ^{2}\leq n%
{\displaystyle\sum_{i=1}^{n}}
a_{i}^{4}$ for any $n\geq2,$ $a_{1},...,a_{n}\in\mathbb{R}.$

Accordingly, we obtain that%
\begin{align*}
\left\vert \Lambda\left(  s;t\right)  \right\vert  &  \leq C\left(
1+\left\vert \hat{X}^{x_{0}}\left(  t\right)  \right\vert ^{4}+\left\vert
v\right\vert ^{4}+\left\vert X^{t,v}\left(  s\right)  \right\vert
^{4}+\left\vert \hat{X}^{x_{0}}\left(  s\right)  \right\vert ^{4}%
+\mathbb{E}_{t}\left[  \left\vert X^{t,v}\left(  s\right)  \right\vert
^{4}\right]  \right. \\
&  \left.  +\mathbb{E}_{t}\left[  \left\vert \hat{X}^{x_{0}}\left(  s\right)
\right\vert ^{4}\right]  +\left\vert \hat{\varphi}^{s}\left(  \hat{X}^{x_{0}%
}\left(  s\right)  \right)  \right\vert ^{4}+\mathbb{E}_{t}\left[  \left\vert
\hat{\varphi}^{s}\left(  \hat{X}^{x_{0}}\left(  s\right)  \right)  \right\vert
^{4}\right]  \right) \\
&  \leq C\left(  1+\left\vert \hat{X}^{x_{0}}\left(  t\right)  \right\vert
^{4}+\left\vert v\right\vert ^{4}+\left\vert X^{t,v}\left(  s\right)
\right\vert ^{4}+\left\vert \hat{X}^{x_{0}}\left(  s\right)  \right\vert
^{4}+\mathbb{E}_{t}\left[  \left\vert X^{t,v}\left(  s\right)  \right\vert
^{4}\right]  \right. \\
&  \left.  +\mathbb{E}_{t}\left[  \left\vert \hat{X}^{x_{0}}\left(  s\right)
\right\vert ^{4}\right]  \right)  \text{,}%
\end{align*}
where we have used (\ref{Eq34}). Thus noting that $\sup\limits_{s\in\left[
t,T\right]  }\mathbb{E}_{t}\left[  \left\vert X\left(  s\right)  \right\vert
^{4}\right]  \leq\mathbb{E}_{t}\left[  \sup\limits_{s\in\left[  t,T\right]
}\left\vert X\left(  s\right)  \right\vert ^{4}\right]  $, for any $X\left(
\cdot\right)  \in\mathcal{S}_{%
\mathcal{F}%
}^{4}\left(  t,T;%
\mathbb{R}
^{n}\right)  $, we get
\begin{align*}
&  \mathbb{E}\left[  \sup_{s\in\left[  t,T\right]  }\left\vert \Lambda\left(
s;t\right)  \right\vert \right] \\
&  \leq C\left(  1+\mathbb{E}\left[  \left\vert \hat{X}^{x_{0}}\left(
t\right)  \right\vert ^{4}\right]  +\mathbb{E}\left[  \left\vert v\right\vert
^{4}\right]  +\mathbb{E}\left[  \sup_{s\in\left[  t,T\right]  }\left\vert
X^{t,v}\left(  s\right)  \right\vert ^{4}\right]  \right. \\
&  \left.  +\mathbb{E}\left[  \sup_{s\in\left[  t,T\right]  }\left\vert
\hat{X}^{x_{0}}\left(  s\right)  \right\vert ^{4}\right]  \right) \\
&  <\infty\text{.}%
\end{align*}

Applying Dominated Convergence Theorem together with (\ref{Eq27}) and
(\ref{Eq37}), we get%
\begin{align*}
&  \lim_{\varepsilon\downarrow0}\frac{1}{\varepsilon}\left\{  \mathbf{J}%
\left(  t,\hat{X}^{x_{0}}\left(  t\right)  ;u^{\varepsilon}\left(
\cdot\right)  \right)  -\mathbf{J}\left(  t,\hat{X}^{x_{0}}\left(  t\right)
;\hat{\varphi}\left(  \cdot,\cdot\right)  \right)  \right\} \\
&  =\lim_{\varepsilon\downarrow0}\frac{1}{\varepsilon}\int_{t}^{t+\varepsilon
}\mathbb{E}_{t}\left[  \Lambda\left(  s;t\right)  \right]  ds\\
&  =\lim_{s\downarrow t}\mathbb{E}_{t}\left[  \Lambda\left(  s;t\right)
\right] \\
&  =\Lambda\left(  t;t\right) \\
&  =\mathcal{H}\left(  t,\hat{X}^{x_{0}}\left(  t\right)  ,t,\hat{X}^{x_{0}%
}\left(  t\right)  ,\hat{X}^{x_{0}}\left(  t\right)  ,v\right) \\
&  -\mathcal{H}\left(  t,\hat{X}^{x_{0}}\left(  t\right)  ,t,\hat{X}^{x_{0}%
}\left(  t\right)  ,\hat{X}^{x_{0}}\left(  t\right)  ,\hat{\varphi}^{t}\left(
\hat{X}^{x_{0}}\left(  t\right)  \right)  \right)  \text{.}%
\end{align*}
Hence, in view of (\ref{Eq36}), we obtain that%
\begin{align*}
&  \underset{\varepsilon\downarrow0}{\lim\inf}\frac{1}{\varepsilon}\left\{
\mathbf{J}\left(  t,\hat{X}^{x_{0}}\left(  t\right)  ;u^{\varepsilon}\left(
\cdot\right)  \right)  -\mathbf{J}\left(  t,\hat{X}^{x_{0}}\left(  t\right)
;\hat{\varphi}\left(  \cdot,\cdot\right)  \right)  \right\} \\
&  =\mathcal{H}\left(  t,\hat{X}^{x_{0}}\left(  t\right)  ,t,\hat{X}^{x_{0}%
}\left(  t\right)  ,\hat{X}^{x_{0}}\left(  t\right)  ,v\right)  -\mathcal{H}%
\left(  t,\hat{X}^{x_{0}}\left(  t\right)  ,t,\hat{X}^{x_{0}}\left(  t\right)
,\hat{X}^{x_{0}}\left(  t\right)  ,\hat{\varphi}^{t}\left(  \hat{X}^{x_{0}%
}\left(  t\right)  \right)  \right) \\
&  \geq0\text{, a.s..}%
\end{align*}
Since $t\in\left[  0,T\right]  $ and $v\in\mathbb{L}^{4}\left(  \Omega
,\mathcal{F}_{t},\mathbb{P};U\right)  $ are arbitrary, we deduce that
$\hat{\varphi}\left(  \cdot,\cdot\right)  $ is an equilibrium strategy. The
equality in (\ref{Eq36}) is an immediate consequence of (\ref{Eq32}%
)-(\ref{Eq33}). This completes the proof.\eop

Let us make some remarks on Theorem \ref{result2}.

\begin{remark}
Note that, the verification argument in Theorem \ref{result2} permits us to
construct an open-loop equilibriums strategy by solving the following coupled
system of IPDEs%
\begin{equation}
\left\{
\begin{array}
[c]{l}%
0=\hat{\theta}_{s}^{t}\left(  s,x,z\right)  +\left\langle \hat{\theta}_{x}%
^{t}\left(  s,x,z\right)  ,\mu\left(  s,x,\hat{\varphi}^{s}\left(  z\right)
\right)  \right\rangle +\left\langle \hat{\theta}_{z}^{t}\left(  s,x,z\right)
,\mu\left(  s,z,\hat{\varphi}^{s}\left(  z\right)  \right)  \right\rangle \\
\text{ \ \ }+\frac{1}{2}\text{\textbf{tr}}\left[  \sigma\sigma^{\top}\left(
s,x,\hat{\varphi}^{s}\left(  z\right)  \right)  \hat{\theta}_{xx}^{t}\left(
s,x,z\right)  \right]  +\frac{1}{2}\text{\textbf{tr}}\left[  \sigma
\sigma^{\top}\left(  s,z,\hat{\varphi}^{s}\left(  z\right)  \right)
\hat{\theta}_{zz}^{t}\left(  s,x,z\right)  \right] \\
\text{ \ }+\text{\textbf{tr}}\left[  \sigma\left(  s,x,\hat{\varphi}%
^{s}\left(  z\right)  \right)  ^{\top}\hat{\theta}_{xz}^{t}\left(
s,x,z\right)  \sigma\left(  s,z,\hat{\varphi}^{s}\left(  z\right)  \right)
\right]  +\text{\ }f\left(  t,s,x,\hat{\varphi}^{s}\left(  z\right)  \right)
\\
\text{ \ \ }+%
{\displaystyle\int\limits_{E}}
\left\{  \hat{\theta}^{t}\left(  s,x+c\left(  s,x,\hat{\varphi}^{s}\left(
z\right)  ,e\right)  ,z+c\left(  s,z,\hat{\varphi}^{s}\left(  z\right)
,e\right)  \right)  -\hat{\theta}^{t}\left(  s,x,z\right)  \right. \\
\text{\ }\left.  -\left\langle \hat{\theta}_{x}^{t}\left(  s,x,z\right)
,c\left(  s,x,\hat{\varphi}^{s}\left(  z\right)  ,e\right)  \right\rangle
-\left\langle \hat{\theta}_{z}^{t}\left(  s,x,z\right)  ,c\left(
s,z,\hat{\varphi}^{s}\left(  z\right)  ,e\right)  \right\rangle \right\}
\vartheta\left(  de\right)  \text{, }\\
0=\hat{g}_{s}^{i}\left(  s,x,z\right)  +\left\langle \hat{g}_{x}^{i}\left(
s,x,z\right)  ,\mu\left(  s,x,\hat{\varphi}^{s}\left(  z\right)  \right)
\right\rangle +\left\langle \hat{g}_{z}^{i}\left(  s,x,z\right)  ,\mu\left(
s,z,\hat{\varphi}^{s}\left(  z\right)  \right)  \right\rangle \\
\text{ \ \ \ }+\frac{1}{2}\text{\textbf{tr}}\left[  \sigma\sigma^{\top}\left(
s,x,\hat{\varphi}^{s}\left(  z\right)  \right)  \hat{g}_{xx}^{i}\left(
s,x,z\right)  \right]  \text{\ }+\frac{1}{2}\text{\textbf{tr}}\left[
\sigma\sigma^{\top}\left(  s,z,\hat{\varphi}^{s}\left(  z\right)  \right)
\hat{g}_{zz}^{i}\left(  s,x,z\right)  \right] \\
\text{ \ \ \ }+\text{\textbf{tr}}\left[  \sigma\left(  s,x,\hat{\varphi}%
^{s}\left(  z\right)  \right)  ^{\top}\hat{g}_{xz}^{i}\left(  s,x,z\right)
\sigma\left(  s,z,\hat{\varphi}^{s}\left(  z\right)  \right)  \right] \\
\text{ \ \ \ }+%
{\displaystyle\int\limits_{E}}
\left\{  \hat{g}^{i}\left(  s,x+c\left(  s,x,\hat{\varphi}^{s}\left(
z\right)  ,e\right)  ,z+c\left(  s,z,\hat{\varphi}^{s}\left(  z\right)
,e\right)  \right)  -\hat{g}^{i}\left(  s,x,z\right)  \right. \\
\text{ \ \ }%
\begin{array}
[c]{r}%
\left.  -\left\langle \hat{g}_{x}^{i}\left(  s,x,z\right)  ,c\left(
s,x,\hat{\varphi}^{s}\left(  z\right)  ,e\right)  \right\rangle -\left\langle
\hat{g}_{z}^{i}\left(  s,x,z\right)  ,c\left(  s,z,\hat{\varphi}^{s}\left(
z\right)  ,e\right)  \right\rangle \right\}  \vartheta\left(  de\right)
\text{, }\\
\text{for }\left(  t,s,x,z\right)  \in D\left[  t,T\right]  \times\left(
\mathbb{R}^{n}\right)  ^{2}\text{,}%
\end{array}
\\
\hat{\varphi}^{s}\left(  z\right)  \in\arg\min\mathcal{H}\left(
s,z,s,z,z,\cdot\right)  ,\text{ for }\left(  s,z\right)  \in\left[
0,T\right]  \times%
\mathbb{R}
^{n}.
\end{array}
\text{ }\right.  \label{Eq38}%
\end{equation}
with the boundary conditions%
\[%
\begin{array}
[c]{r}%
\hat{\theta}^{t}\left(  T,x,z\right)  =F\left(  t,x\right)  \text{, }\hat
{g}^{i}\left(  T,x,z\right)  =\Psi_{i}\left(  x\right)  \text{, for }\left(
t,x,z\right)  \in\left[  0,T\right]  \times\left(  \mathbb{R}^{n}\right)
^{2}\text{,}\\
\text{for }1\leq i\leq m.
\end{array}
\]
To the best of our knowledge, coupled IPDEs of the above form appear for the
first time in the literature. However, proving the general existence for these
IPDEs remains an outstanding open problem even for the simplest case when
$n=m=d=1$.
\end{remark}

\begin{remark}
\label{Rem1}Compared with some existing studies in the literature, our
PDEs-method demonstrates several new advantages on the treatment of open-loop
equilibrium controls. Some of them are listed as follows:\newline\textbf{(i)}
Unlike \cite{Huetal2011}, \cite{Huetal2017}, \cite{Djehiche}, \cite{Sun}, our
approach avoids the complex variation arguments of the second order expansion
in the spike variation. Moreover, our results do not require any
differentiability assumptions on the involved coefficients, except the
function $G\left(  \cdot,\cdot,\cdot\right)  $.\newline\textbf{(ii)} By
solving the system of IPDEs (\ref{Eq19})-(\ref{Eq21}), one can derive
simultaneously the equilibrium strategy $\hat{\varphi}\left(  \cdot
,\cdot\right)  $, as well as, its corresponding equilibrium objective value
$\mathbf{J}\left(  t,\hat{X}^{x_{0}}\left(  t\right)  ;\hat{\varphi}\left(
\cdot,\cdot\right)  \right)  $, at each $t\in\left[  0,T\right]  $; this is
different from most of the existing methodologies in the literature, since
they does not provide the objective value $\mathbf{J}\left(  t,\hat{X}^{x_{0}%
}\left(  t\right)  ;\hat{u}\left(  \cdot\right)  \right)  $.\newline%
\textbf{(iii) }Our methodology in constructing an open-loop equilibrium
strategy is purely deterministic; indeed an equilibrium strategy $\hat
{\varphi}\left(  \cdot,\cdot\right)  $ can be constructed by solving
(\ref{Eq38}), which is actually a deterministic coupled system of
IPDEs.\newline\textbf{(iv)} Even thought we are concerned with the open-loop
equilibrium framework, our approach looks more like a dynamical programming
than a maximum principle.\newline
\end{remark}

\begin{remark}
Note that some of the above-listed "advantages" are also satisfied by the
BSPDEs-approach introduced by Alia~\cite{Alia2}. However, since Problem (N) is
formulated in the Markovian framework, it seems to us that it is more natural
to use PDEs, instead of BSPDEs, to characterize open-loop equilibriums.
\end{remark}

\begin{remark}
The PDEs-approach is closely linked to the strong Markovian property of the
controlled process. Thus when the problem is not Markovian, the PDEs-approach
does not apply while the variational method of Hu et al. (\cite{Huetal2011}%
,\cite{Huetal2017}) does.
\end{remark}

\section{Connection between the PDEs-approach and the variational
approach\label{section5}}

As previously mentioned in Introduction section, the idea of defining the
equilibrium control within the whole class of open-loop controls goes back to
Hu et al. (\cite{Huetal2011}, \cite{Huetal2017}), where the authors undertook
a deep study of a time-inconsistent stochastic LQ model. They performed a
variational method in the spirit of Peng's stochastic maximum principle and
constructed an equilibrium solution by solving a \textit{"flow"} of coupled
FBSDEs. In this section, we briefly discuss the connection between the PDEs
approach of the present paper and the variational approach of Hu et al.
(\cite{Huetal2011}, \cite{Huetal2017}). For sake of simplicity and clarity of
the presentation, we suppose that $\Psi\left(  x\right)  \equiv x$ and all the
coefficients are assumed to be one dimensional (i.e. $n=d=l=1$). Specially,
the PDEs (\ref{Eq19})-(\ref{Eq21})\ reduce, respectively, to%
\begin{equation}
\left\{
\begin{array}
[c]{l}%
0=\hat{\theta}_{s}^{t}\left(  s,x,z\right)  +\hat{\theta}_{x}^{t}\left(
s,x,z\right)  \mu\left(  s,x,\hat{\varphi}^{s}\left(  z\right)  \right)
+\hat{\theta}_{z}^{t}\left(  s,x,z\right)  \mu\left(  s,z,\hat{\varphi}%
^{s}\left(  z\right)  \right) \\
\text{ \ \ }+\frac{1}{2}\sigma\left(  s,x,\hat{\varphi}^{s}\left(  z\right)
\right)  ^{2}\hat{\theta}_{xx}^{t}\left(  s,x,z\right)  +\frac{1}{2}%
\sigma\left(  s,z,\hat{\varphi}^{s}\left(  z\right)  \right)  ^{2}\hat{\theta
}_{zz}^{t}\left(  s,x,z\right) \\
\text{ \ }+\sigma\left(  s,x,\hat{\varphi}^{s}\left(  z\right)  \right)
\sigma\left(  s,z,\hat{\varphi}^{s}\left(  z\right)  \right)  \hat{\theta
}_{xz}^{t}\left(  s,x,z\right)  +f\left(  t,s,x,\hat{\varphi}^{s}\left(
z\right)  \right) \\
\text{ \ \ }+%
{\displaystyle\int_{E}}
\left\{  \hat{\theta}^{t}\left(  s,x+c\left(  s,x,\hat{\varphi}^{s}\left(
z\right)  ,e\right)  ,z+c\left(  s,z,\hat{\varphi}^{s}\left(  z\right)
,e\right)  \right)  -\hat{\theta}^{t}\left(  s,x,z\right)  \right. \\
\text{\ }%
\begin{array}
[c]{r}%
\left.  -\hat{\theta}_{x}^{t}\left(  s,x,z\right)  c\left(  s,x,\hat{\varphi
}^{s}\left(  z\right)  ,e\right)  -\hat{\theta}_{z}^{t}\left(  s,x,z\right)
c\left(  s,z,\hat{\varphi}^{s}\left(  z\right)  ,e\right)  \right\}
\vartheta\left(  de\right)  \text{, }\\
\text{for }\left(  s,x,z\right)  \in\left[  t,T\right]  \times\mathbb{R}%
\times\mathbb{R}\text{,}%
\end{array}
\\
\hat{\theta}^{t}\left(  T,x,z\right)  =F\left(  t,x\right)  \text{, for
}\left(  x,z\right)  \in\mathbb{R}\times\mathbb{R}\text{,}%
\end{array}
\right.  \label{Eq70}%
\end{equation}
and%
\begin{equation}
\left\{
\begin{array}
[c]{l}%
0=\hat{g}_{s}\left(  s,x,z\right)  +\hat{g}_{x}\left(  s,x,z\right)
\mu\left(  s,x,\hat{\varphi}^{s}\left(  z\right)  \right)  +\hat{g}_{z}%
^{i}\left(  s,x,z\right)  \mu\left(  s,z,\hat{\varphi}^{s}\left(  z\right)
\right) \\
\text{ \ \ }+\frac{1}{2}\sigma\left(  s,x,\hat{\varphi}^{s}\left(  z\right)
\right)  ^{2}\hat{g}_{xx}^{i}\left(  s,x,z\right)  +\frac{1}{2}\sigma\left(
s,z,\hat{\varphi}^{s}\left(  z\right)  \right)  ^{2}\hat{g}_{zz}^{i}\left(
s,x,z\right) \\
\text{ \ }+\sigma\left(  s,x,\hat{\varphi}^{s}\left(  z\right)  \right)
\hat{g}_{xz}^{i}\left(  s,x,z\right)  \sigma\left(  s,z,\hat{\varphi}%
^{s}\left(  z\right)  \right) \\
\text{ \ }+%
{\displaystyle\int_{E}}
\left\{  \hat{g}\left(  s,x+c\left(  s,x,\hat{\varphi}^{s}\left(  z\right)
,e\right)  ,z+c\left(  s,z,\hat{\varphi}^{s}\left(  z\right)  ,e\right)
\right)  -\hat{g}\left(  s,x,z\right)  \right. \\%
\begin{array}
[c]{r}%
\left.  -\hat{g}_{x}\left(  s,x,z\right)  c\left(  s,x,\hat{\varphi}%
^{s}\left(  z\right)  ,e\right)  -\hat{g}_{z}\left(  s,x,z\right)  c\left(
s,z,\hat{\varphi}^{s}\left(  z\right)  ,e\right)  \right\}  \vartheta\left(
de\right)  \text{, }\\
\text{for }\left(  s,x,z\right)  \in\left[  0,T\right]  \times\mathbb{R}%
^{n}\times\mathbb{R}^{n}\text{,}%
\end{array}
\\
\hat{g}\left(  T,x,z\right)  =x\text{, for }\left(  x,z\right)  \in
\mathbb{R}\times\mathbb{R}\text{.}%
\end{array}
\right.  \label{Eq71}%
\end{equation}

In order to be able to establish the link between our results and the
variational approach of Hu et al. (\cite{Huetal2011}, \cite{Huetal2017}), we
first need to derive a version of Peng's type stochastic maximum principle
that characterizes open-loop equilibrium controls. To this end, let
$\hat{\varphi}\left(  \cdot,\cdot\right)  \in\mathcal{S}$ be a fixed
admissible strategy and $\hat{X}^{x_{0}}\left(  \cdot\right)  $ be the unique
strong solution of the SDE (\ref{Eq8*}). For some fixed arbitrary $u\in U$, we
put for $\rho=\mu,\sigma$:%
\[
\left\{
\begin{array}
[c]{l}%
\rho\left(  t\right)  =\rho\left(  t,\hat{X}^{x_{0}}\left(  t\right)
,\hat{\varphi}^{t}\left(  \hat{X}^{x_{0}}\left(  t\right)  \right)  \right)
\text{, }c\left(  t,e\right)  =c\left(  t,\hat{X}^{x_{0}}\left(  t\right)
,\hat{\varphi}^{t}\left(  \hat{X}^{x_{0}}\left(  t\right)  \right)  ,e\right)
\text{,}\\
\rho_{x}\left(  t\right)  =\rho_{x}\left(  t,\hat{X}^{x_{0}}\left(  t\right)
,\hat{\varphi}^{t}\left(  \hat{X}^{x_{0}}\left(  t\right)  \right)  \right)
\text{, }\rho_{xx}\left(  t\right)  =\rho_{xx}\left(  t,\hat{X}^{x_{0}}\left(
t\right)  ,\hat{\varphi}^{t}\left(  \hat{X}^{x_{0}}\left(  t\right)  \right)
\right)  \text{,}\\
\delta\rho\left(  t;u\right)  =\rho\left(  t,\hat{X}^{x_{0}}\left(  t\right)
,u\right)  -\rho\left(  t,\hat{X}^{x_{0}}\left(  t\right)  ,\hat{\varphi}%
^{t}\left(  \hat{X}^{x_{0}}\left(  t\right)  \right)  \right)  \text{,}\\
\delta f\left(  t,s;u\right)  =f\left(  t,s,\hat{X}^{x_{0}}\left(  t\right)
,u\right)  -f\left(  t,s,\hat{X}^{x_{0}}\left(  t\right)  ,\hat{\varphi}%
^{t}\left(  \hat{X}^{x_{0}}\left(  t\right)  \right)  \right)  \text{,}\\
c_{x}\left(  t,e\right)  =c_{x}\left(  t,\hat{X}^{x_{0}}\left(  t\right)
,\hat{\varphi}^{t}\left(  \hat{X}^{x_{0}}\left(  t\right)  \right)  ,e\right)
\text{, }c_{xx}\left(  t,e\right)  =c_{xx}\left(  t,\hat{X}^{x_{0}}\left(
t\right)  ,\hat{\varphi}^{t}\left(  \hat{X}^{x_{0}}\left(  t\right)  \right)
,e\right)  \text{,}\\
\delta c\left(  t,e;u\right)  =c\left(  t,\hat{X}^{x_{0}}\left(  t\right)
,u,e\right)  -c\left(  t,\hat{X}^{x_{0}}\left(  t\right)  ,\hat{\varphi}%
^{t}\left(  \hat{X}^{x_{0}}\left(  t\right)  \right)  ,e\right)  \text{.}%
\end{array}
\right.
\]

We impose the following assumptions.

\begin{enumerate}
\item[(\textbf{H1*)}] The maps $\mu\left(  s,x,u\right)  $, $\sigma\left(
s,x,u\right)  $ and $c\left(  s,x,u,e\right)  $ are twice continuously
differentiable with respect to $x.$ They and their derivatives in $x$ are
continuous in $\left(  s,x,u\right)  $, and bounded.

\item[\textbf{(H2*)}] \textbf{(i)} The functions $F\left(  t,x\right)  $ and
$f\left(  t,s,x,u\right)  $ are twice continuously differentiable with respect
to $x.$ They and their derivatives in $x$ are continuous in $\left(
s,x,u\right)  $, and bounded.\newline\textbf{(ii)} The function $G\left(
t,y,\bar{x}\right)  $ is twice continuously differentiable with respect to
$\bar{x}$. $G$ and its derivatives in $\bar{x}$ are continuous in $\bar{x}$,
and bounded.
\end{enumerate}

We point out that Assumptions \textbf{(H1*)-(H2*)} can be substantially
relaxed, but we do not focus on this here. For any $t\in\left[  0,T\right]  $,
define in the time interval $\left[  t,T\right]  $ the processes $\left(
p^{t}\left(  \cdot\right)  ,q^{t}\left(  \cdot\right)  ,r^{t}\left(
\cdot,\cdot\right)  \right)  \in\mathcal{S}_{%
\mathcal{F}%
}^{2}\left(  t,T;\mathbb{%
\mathbb{R}
}\right)  \times\mathcal{L}_{%
\mathcal{F}%
}^{2}\left(  t,T;\mathbb{%
\mathbb{R}
}\right)  \times\mathcal{L}_{%
\mathcal{F}%
,p}^{\vartheta,q}\left(  \left[  t,T\right]  \times E;\mathbb{%
\mathbb{R}
}\right)  $ and $\left(  P^{t}\left(  \cdot\right)  ,\Phi^{t}\left(
\cdot\right)  ,\Upsilon^{t}\left(  \cdot,\cdot\right)  \right)  \in
\mathcal{S}_{%
\mathcal{F}%
}^{2}\left(  t,T;\mathbb{%
\mathbb{R}
}\right)  \times\mathcal{L}_{%
\mathcal{F}%
}^{2}\left(  t,T;\mathbb{%
\mathbb{R}
}\right)  \times\mathcal{L}_{%
\mathcal{F}%
,p}^{\vartheta,2}\left(  \left[  t,T\right]  \times E;\mathbb{%
\mathbb{R}
}\right)  $ as the solution of the following equations:%
\begin{equation}
\left\{
\begin{array}
[c]{l}%
dp^{t}\left(  s\right)  =-\left\{  \mu_{x}\left(  s\right)  p^{t}\left(
s\right)  +\sigma_{x}\left(  s\right)  q^{t}\left(  s\right)  +\int_{E}%
c_{x}\left(  s,e\right)  r^{t}\left(  s,e\right)  \vartheta\left(  de\right)
\right. \\
\text{ \ \ \ \ \ \ \ \ \ \ \ \ }\left.  +f_{x}\left(  t,s\right)  \right\}
ds+q^{t}\left(  s\right)  dW\left(  s\right)  +\int_{E}r\left(  s,e\right)
\tilde{N}\left(  ds,de\right)  \text{, }s\in\left[  t,T\right]  \text{,}\\
p^{t}\left(  T\right)  =F_{x}\left(  t,\hat{X}^{x_{0}}\left(  T\right)
\right)  +G_{\bar{x}}\left(  t,\hat{X}^{x_{0}}\left(  t\right)  ,\mathbb{E}%
_{t}\left[  \hat{X}^{x_{0}}\left(  T\right)  \right]  \right)  \text{;}%
\end{array}
\right.  \label{Eq72}%
\end{equation}%
\begin{equation}
\left\{
\begin{array}
[c]{l}%
dP^{t}\left(  s\right)  =-\left\{  \left(  2\mu_{x}\left(  s\right)
+\sigma_{x}\left(  s\right)  ^{2}\right)  P^{t}\left(  s\right)  +2\sigma
_{x}\left(  s\right)  \Phi^{t}\left(  s\right)  \right. \\
\text{ \ \ \ \ \ \ \ \ \ \ \ \ \ \ \ \ \ }\left.  +%
{\displaystyle\int_{E}}
\left\{  \left(  \Upsilon^{t}\left(  s,e\right)  +P^{t}\left(  s\right)
\right)  c_{x}\left(  s,e\right)  ^{2}+2c_{x}\left(  s,e\right)  \Upsilon
^{t}\left(  s,e\right)  \right\}  \vartheta\left(  de\right)  \right. \\
\text{ \ \ \ \ \ \ \ \ \ \ \ \ \ \ \ \ \ }+\left.  \mathbb{H}_{xx}\left(
s,\hat{X}^{x_{0}}\left(  s\right)  ,\hat{\varphi}^{s}\left(  \hat{X}^{x_{0}%
}\left(  s\right)  \right)  ,p^{t}\left(  s\right)  ,q^{t}\left(  s\right)
,r^{t}\left(  s,\cdot\right)  \right)  \right\}  ds\\
\text{ \ \ \ \ \ \ \ \ \ \ \ \ \ \ \ \ \ }+\Phi^{t}\left(  s\right)  dW\left(
s\right)  +%
{\displaystyle\int_{E}}
\Upsilon^{t}\left(  s,e\right)  \tilde{N}\left(  ds,de\right)  \text{, }%
s\in\left[  t,T\right]  \text{,}\\
P^{t}\left(  T\right)  =F_{xx}\left(  t,\hat{X}^{x_{0}}\left(  T\right)
\right)  \text{,}%
\end{array}
\right.  \label{Eq73}%
\end{equation}
where the Hamiltonian $\mathbb{H}$ is defined by%
\begin{align*}
\mathbb{H}\left(  t,s,x,u,p,q,r\left(  \cdot\right)  \right)   &  :=\mu\left(
s,x,u\right)  p+\sigma\left(  s,x,u\right)  q\\
&  +\int_{E}c\left(  s,x,u,e\right)  r\left(  e\right)  \vartheta\left(
de\right)  +f\left(  t,s,x,u\right)  \text{,}%
\end{align*}
for any $\left(  t,s,x,u,p,q,r\left(  \cdot\right)  \right)  \in\left[
0,T\right]  \times\left[  0,T\right]  \times%
\mathbb{R}
\times U\times%
\mathbb{R}
\times%
\mathbb{R}
\times\mathbb{L}^{2}\left(  E,\mathcal{B}\left(  E\right)  ,\vartheta;%
\mathbb{R}
\right)  $.

We also need to introduce an $\mathcal{\bar{H}}$-function associated with the
family of processes $\left(  \hat{\varphi}\left(  \cdot,\hat{X}^{x_{0}}\left(
\cdot\right)  \right)  ,\hat{X}^{x_{0}}\left(  \cdot\right)  \right.
$\newline$\left.  ,\left\{  \left(  p^{t}\left(  \cdot\right)  ,q^{t}\left(
\cdot\right)  ,r^{t}\left(  \cdot,\cdot\right)  ,P^{t}\left(  \cdot\right)
,\Gamma^{t}\left(  \cdot,\cdot\right)  \right)  \right\}  _{t\in\left[
0,T\right]  }\right)  $ as follows:%
\begin{align}
\mathcal{\bar{H}}\left(  t,s,x,u\right)   &  :=\mathbb{H}\left(
t,s,x,u,p^{t}\left(  s\right)  ,q^{t}\left(  s\right)  ,r^{t}\left(
s,\cdot\right)  \right) \nonumber\\
&  +\frac{1}{2}P^{t}\left(  s\right)  \left(  \sigma\left(  s,x,u\right)
-\sigma\left(  s,\hat{X}^{x_{0}}\left(  s\right)  ,\hat{\varphi}^{s}\left(
\hat{X}^{x_{0}}\left(  s\right)  \right)  \right)  \right)  ^{2}\label{Eq74}\\
&  +\frac{1}{2}%
{\displaystyle\int_{E}}
\left(  \Upsilon^{t}\left(  s,e\right)  +P^{t}\left(  s\right)  \right)
\left(  c\left(  s,x,u,z\right)  -c\left(  s,\hat{X}^{x_{0}}\left(  s\right)
,\hat{\varphi}^{s}\left(  \hat{X}^{x_{0}}\left(  s\right)  \right)  ,e\right)
\right)  ^{2}\vartheta\left(  de\right)  \text{,}\nonumber
\end{align}
for any $\left(  t,s,u,x\right)  \in D\left[  0,T\right]  \times U\times%
\mathbb{R}
$.

The following theorem is comparable with Theorem \ref{result2}.

\begin{theorem}
[Stochastic Maximum Principle]\label{result3}Let \textbf{(H1*)-(H2*)
}hold\textbf{. }Given an admissible strategy $\hat{\varphi}\left(  \cdot
,\cdot\right)  \in\mathcal{S}$, let for each $t\in\left[  0,T\right]  $,
$\left(  p^{t}\left(  \cdot\right)  ,q^{t}\left(  \cdot\right)  ,r^{t}\left(
\cdot,\cdot\right)  \right)  \in\mathcal{S}_{%
\mathcal{F}%
}^{2}\left(  t,T;\mathbb{%
\mathbb{R}
}\right)  \times\mathcal{L}_{%
\mathcal{F}%
}^{2}\left(  t,T;\mathbb{%
\mathbb{R}
}\right)  \times\mathcal{L}_{%
\mathcal{F}%
,p}^{\vartheta,q}\left(  \left[  t,T\right]  \times E;\mathbb{%
\mathbb{R}
}\right)  $ and $\left(  P^{t}\left(  \cdot\right)  ,\Phi^{t}\left(
\cdot\right)  ,\Upsilon^{t}\left(  \cdot,\cdot\right)  \right)  \in
\mathcal{S}_{%
\mathcal{F}%
}^{2}\left(  t,T;\mathbb{%
\mathbb{R}
}\right)  \times\mathcal{L}_{%
\mathcal{F}%
}^{2}\left(  t,T;\mathbb{%
\mathbb{R}
}\right)  \times\mathcal{L}_{%
\mathcal{F}%
,p}^{\vartheta,q}\left(  \left[  t,T\right]  \times E;\mathbb{%
\mathbb{R}
}\right)  $ be the unique solutions to the BSDEs (\ref{Eq72}) and
(\ref{Eq73}), respectively. Suppose that the following hold:\newline%
\textbf{(i) }$\hat{\varphi}\left(  \cdot,\cdot\right)  $ is a continuous
function.\newline\textbf{(ii) }For each $\left(  t,e\right)  \in\left[
0,T\right]  \times E,$ $s\rightarrow$\textbf{ }$\left(  q^{t}\left(  s\right)
,r^{t}\left(  s,e\right)  ,\Upsilon^{t}\left(  s,e\right)  \right)  $ is
right-continuous at $s=t$ a.s..\newline\textbf{(iii)} For each $t\in\left[
0,T\right]  $, there exists a constant $C>0$ such that%
\begin{align*}
C  &  \geq\mathbb{E}\left[  \sup_{s\in\left[  t,T\right]  }\left\vert
q^{t}\left(  s\right)  \right\vert ^{2}\right]  +\mathbb{E}\left[  \int
_{E}\sup_{s\in\left[  t,T\right]  }\left\vert r^{t}\left(  s,e\right)
\right\vert ^{2}\vartheta\left(  de\right)  \right] \\
&  +\mathbb{E}\left[  \int_{E}\sup_{s\in\left[  t,T\right]  }\left\vert
\Upsilon^{t}\left(  s,e\right)  \right\vert ^{2}\vartheta\left(  de\right)
\right]  \text{, a.s..}%
\end{align*}
\newline\textbf{(iv)} For any $t\in\left[  0,T\right]  $,%
\[
\mathcal{\bar{H}}\left(  t,t,\hat{X}_{-}^{x_{0}}\left(  t\right)
,\hat{\varphi}^{t}\left(  \hat{X}_{-}^{x_{0}}\left(  t\right)  \right)
\right)  =\min_{u\in U}\mathcal{\bar{H}}\left(  t,t,\hat{X}_{-}^{x_{0}}\left(
t\right)  ,u\right)  \text{.}%
\]
\newline Then $\hat{\varphi}\left(  \cdot,\cdot\right)  $ is an open-loop
equilibrium strategy.
\end{theorem}

\bop The proof follows in a similar way as in [\cite{Alia2}, Theorem 4.1]. We
omit the details here.\eop

We are now ready to state the main result of this section. For brevity, we put
for $\varrho=\theta^{t},g$%
\[
\left\{
\begin{array}
[c]{l}%
\varrho\left(  s\right)  =\varrho\left(  s,\hat{X}^{x_{0}}\left(  s\right)
,\hat{X}^{x_{0}}\left(  s\right)  \right)  \text{, \ }\varrho_{x}\left(
s\right)  =\varrho_{x}\left(  s,\hat{X}^{x_{0}}\left(  s\right)  ,\hat
{X}^{x_{0}}\left(  s\right)  \right)  \text{,}\\
\varrho_{xx}\left(  s\right)  =\varrho_{xx}\left(  s,\hat{X}^{x_{0}}\left(
s\right)  ,\hat{X}^{x_{0}}\left(  s\right)  \right)  \text{, }\varrho
_{xz}\left(  s\right)  =\varrho_{xz}\left(  s,\hat{X}^{x_{0}}\left(  s\right)
,\hat{X}^{x_{0}}\left(  s\right)  \right)  \text{,}\\
\varrho_{xxx}\left(  s\right)  =\varrho_{xxx}\left(  s,\hat{X}^{x_{0}}\left(
s\right)  ,\hat{X}^{x_{0}}\left(  s\right)  \right)  \text{, }\varrho
_{xxz}\left(  s\right)  =\varrho_{xxz}\left(  s,\hat{X}^{x_{0}}\left(
s\right)  ,\hat{X}^{x_{0}}\left(  s\right)  \right)  \text{,}\\
\Delta\varrho\left(  s,e\right)  =\varrho\left(  s,\hat{X}_{-}^{x_{0}}\left(
s\right)  +\hat{c}\left(  s,e\right)  ,\hat{X}_{-}^{x_{0}}\left(  s\right)
+\hat{c}\left(  s,e\right)  \right)  -\varrho\left(  s,\hat{X}_{-}^{x_{0}%
}\left(  s\right)  ,\hat{X}_{-}^{x_{0}}\left(  s\right)  \right)  \text{,}\\
\Delta\varrho_{x}\left(  s,e\right)  =\varrho_{x}\left(  s,\hat{X}_{-}^{x_{0}%
}\left(  s\right)  +\hat{c}\left(  s,e\right)  ,\hat{X}_{-}^{x_{0}}\left(
s\right)  +\hat{c}\left(  s,e\right)  \right)  -\varrho_{x}\left(  s,\hat
{X}_{-}^{x_{0}}\left(  s\right)  ,\hat{X}_{-}^{x_{0}}\left(  s\right)
\right)  \text{,}\\
\Delta\varrho_{xx}\left(  s,e\right)  =\varrho_{xx}\left(  s,\hat{X}%
_{-}^{x_{0}}\left(  s\right)  +\hat{c}\left(  s,e\right)  ,\hat{X}_{-}^{x_{0}%
}\left(  s\right)  +\hat{c}\left(  s,e\right)  \right)  -\varrho_{xx}\left(
s,\hat{X}_{-}^{x_{0}}\left(  s\right)  ,\hat{X}_{-}^{x_{0}}\left(  s\right)
\right)  \text{,}%
\end{array}
\right.
\]
where%
\[
\hat{c}\left(  s,e\right)  :=c\left(  s,\hat{X}_{-}^{x_{0}}\left(  s\right)
,\hat{\varphi}^{s}\left(  \hat{X}_{-}^{x_{0}}\left(  s\right)  \right)
,e\right)  \text{.}%
\]

\begin{theorem}
\label{result4}Suppose that \textbf{(H1*)-(H2*)} hold. Let $\hat{\varphi
}\left(  \cdot,\cdot\right)  \in\mathcal{S}$ be a fixed admissible control and
$\hat{X}^{x_{0}}\left(  \cdot\right)  $ be the unique strong solution of the
SDE (\ref{Eq8*}). Suppose that, for each $t\in\left[  0,T\right]  ,$ the
BSPDEs (\ref{Eq70})-(\ref{Eq71}) admit two unique classical solutions
$\hat{\theta}^{t}\left(  \cdot,\cdot,\cdot\right)  \in\mathcal{C}%
^{1,4,4}\left(  \left[  t,T\right]  \times%
\mathbb{R}
\times%
\mathbb{R}
;%
\mathbb{R}
\right)  $ and $\hat{g}\left(  \cdot,\cdot,\cdot\right)  \in\mathcal{C}%
^{1,4,4}\left(  \left[  t,T\right]  \times%
\mathbb{R}
\times%
\mathbb{R}
;%
\mathbb{R}
\right)  $, respectively, such that%
\[
\left(  \hat{\theta}_{x}^{t}\left(  \cdot,\cdot,\cdot\right)  ,\hat{\theta
}_{xx}^{t}\left(  \cdot,\cdot,\cdot\right)  \right)  \in\mathcal{C}%
_{2}^{1,2,2}\left(  \left[  t,T\right]  \times%
\mathbb{R}
\times%
\mathbb{R}
;%
\mathbb{R}
\right)  ^{2}%
\]
and%
\[
\left(  \hat{g}_{x}\left(  \cdot,\cdot,\cdot\right)  ,\hat{g}_{xx}\left(
\cdot,\cdot,\cdot\right)  \right)  \in\mathcal{C}_{2}^{1,2,2}\left(  \left[
0,T\right]  \times%
\mathbb{R}
\times%
\mathbb{R}
;%
\mathbb{R}
\right)  ^{2}\text{.}%
\]
For each $0\leq t\leq s\leq T$, $e\in E$, define%
\begin{equation}
\left\{
\begin{array}
[c]{l}%
p^{t}\left(  s\right)  :=\hat{\theta}_{x}^{t}\left(  s\right)  +G_{\bar{x}%
}\left(  t,\hat{X}^{x_{0}}\left(  t\right)  ,\mathbb{E}_{t}\left[  \hat
{X}^{x_{0}}\left(  T\right)  \right]  \right)  \hat{g}_{x}\left(  s\right)
,\\
q^{t}\left(  s\right)  :=\sigma\left(  s\right)  \left(  \hat{\theta}_{xx}%
^{t}\left(  s\right)  +G_{\bar{x}}\left(  t,\hat{X}^{x_{0}}\left(  t\right)
,\mathbb{E}_{t}\left[  \hat{X}^{x_{0}}\left(  T\right)  \right]  \right)
\hat{g}_{xx}\left(  s\right)  \right) \\
\text{ \ \ \ \ \ \ \ \ \ \ \ }+\sigma\left(  s\right)  \left(  \hat{\theta
}_{xz}^{t}\left(  s\right)  +G_{\bar{x}}\left(  t,\hat{X}^{x_{0}}\left(
t\right)  ,\mathbb{E}_{t}\left[  \hat{X}^{x_{0}}\left(  T\right)  \right]
\right)  \hat{g}_{xz}\left(  s\right)  \right)  \text{,}\\
r^{t}\left(  s,e\right)  :=\Delta\hat{\theta}_{x}\left(  s,e\right)
+G_{\bar{x}}\left(  t,\hat{X}^{x_{0}}\left(  t\right)  ,\mathbb{E}_{t}\left[
\hat{X}^{x_{0}}\left(  T\right)  \right]  \right)  \Delta\hat{g}_{x}\left(
s,e\right)  \text{,}%
\end{array}
\right.  \label{Eq76}%
\end{equation}%
\begin{equation}
\left\{
\begin{array}
[c]{l}%
P^{t}\left(  s\right)  :=\hat{\theta}_{xx}^{t}\left(  s\right)  +G_{\bar{x}%
}\left(  t,\hat{X}^{x_{0}}\left(  t\right)  ,\mathbb{E}_{t}\left[  \hat
{X}^{x_{0}}\left(  T\right)  \right]  \right)  \hat{g}_{xx}\left(  s\right)
\text{,}\\
\Phi^{t}\left(  s\right)  :=\sigma\left(  s\right)  \left(  \hat{\theta}%
_{xxx}^{t}\left(  s\right)  +G_{\bar{x}}\left(  t,\hat{X}^{x_{0}}\left(
t\right)  ,\mathbb{E}_{t}\left[  \hat{X}^{x_{0}}\left(  T\right)  \right]
\right)  \hat{g}_{xxx}\left(  s\right)  \right) \\
\text{ \ \ \ \ \ \ \ }+\sigma\left(  s\right)  \left(  \hat{\theta}_{xxz}%
^{t}\left(  s\right)  +G_{\bar{x}}\left(  t,\hat{X}^{x_{0}}\left(  t\right)
,\mathbb{E}_{t}\left[  \hat{X}^{x_{0}}\left(  T\right)  \right]  \right)
\hat{g}_{xxz}\left(  s\right)  \right)  \text{,}\\
\Upsilon^{t}\left(  s,e\right)  :=\Delta\hat{\theta}_{xx}\left(  s,e\right)
+G_{\bar{x}}\left(  t,\hat{X}^{x_{0}}\left(  t\right)  ,\mathbb{E}_{t}\left[
\hat{X}^{x_{0}}\left(  T\right)  \right]  \right)  \Delta\hat{g}_{xx}\left(
s,e\right)  \text{.}%
\end{array}
\right.  \label{Eq77}%
\end{equation}
Then $\left(  p^{t}\left(  \cdot\right)  ,q^{t}\left(  \cdot\right)
,r^{t}\left(  \cdot,\cdot\right)  \right)  $ satisfies the BSDE (\ref{Eq72})
and $\left(  P^{t}\left(  \cdot\right)  ,\Phi^{t}\left(  \cdot\right)
,\Upsilon^{t}\left(  \cdot,\cdot\right)  \right)  $ satisfies the BSDE
(\ref{Eq73}).
\end{theorem}

\bop Let $t\in\left[  0,T\right]  $ be fixed. Differentiate equations
(\ref{Eq70})-(\ref{Eq71}) in $x$, we obtain that $\hat{\theta}_{x}^{t}\left(
\cdot,\cdot,\cdot\right)  $ and $\hat{g}_{x}\left(  \cdot,\cdot,\cdot\right)
$ satisfy the following IPDEs:%

\[
\left\{
\begin{array}
[c]{l}%
0=\hat{\theta}_{sx}^{t}\left(  s,x,z\right)  +\hat{\theta}_{xx}^{t}\left(
s,x,z\right)  \mu\left(  s,x,\hat{\varphi}^{s}\left(  z\right)  \right)
+\hat{\theta}_{x}^{t}\left(  s,x,z\right)  \mu_{x}\left(  s,x,\hat{\varphi
}^{s}\left(  z\right)  \right) \\
\text{ \ }+\hat{\theta}_{zx}^{t}\left(  s,x,z\right)  \mu\left(
s,z,\hat{\varphi}^{s}\left(  z\right)  \right)  +\frac{1}{2}\sigma\left(
s,z,\hat{\varphi}^{s}\left(  z\right)  \right)  ^{2}\hat{\theta}_{zzx}%
^{t}\left(  s,x,z\right) \\
\text{ \ }+\sigma\left(  s,x,\hat{\varphi}^{s}\left(  z\right)  \right)
\sigma_{x}\left(  s,x,\hat{\varphi}^{s}\left(  z\right)  \right)  \hat{\theta
}_{xx}^{t}\left(  s,x,z\right)  +\frac{1}{2}\sigma\left(  s,x,\hat{\varphi
}^{s}\left(  z\right)  \right)  ^{2}\hat{\theta}_{xxx}^{t}\left(  s,x,z\right)
\\
\text{ \ }+\sigma_{x}\left(  s,x,\hat{\varphi}^{s}\left(  z\right)  \right)
\sigma\left(  s,z,\hat{\varphi}^{s}\left(  z\right)  \right)  \hat{\theta
}_{xz}^{t}\left(  s,x,z\right)  +\sigma\left(  s,x,\hat{\varphi}^{s}\left(
z\right)  \right)  \sigma\left(  s,z,\hat{\varphi}^{s}\left(  z\right)
\right)  \hat{\theta}_{xxz}^{t}\left(  s,x,z\right) \\
\text{ \ }+%
{\displaystyle\int_{E}}
\left\{  \hat{\theta}_{x}^{t}\left(  s,x+c\left(  s,x,\hat{\varphi}^{s}\left(
z\right)  ,e\right)  ,z+c\left(  s,z,\hat{\varphi}^{s}\left(  z\right)
,e\right)  \right)  \left(  1+c_{x}\left(  s,x,\hat{\varphi}^{s}\left(
z\right)  ,e\right)  \right)  \right. \\
\text{ \ \ }-\hat{\theta}_{x}^{t}\left(  s,x,z\right)  -\hat{\theta}_{xx}%
^{t}\left(  s,x,z\right)  c\left(  s,x,\hat{\varphi}^{s}\left(  z\right)
,e\right)  -\hat{\theta}_{x}^{t}\left(  s,x,z\right)  c_{x}\left(
s,x,\hat{\varphi}^{s}\left(  z\right)  ,e\right) \\
\text{ \ \ \ }%
\begin{array}
[c]{r}%
\left.  -\hat{\theta}_{zx}^{t}\left(  s,x,z\right)  c\left(  s,z,\hat{\varphi
}^{s}\left(  z\right)  ,e\right)  \right\}  \vartheta\left(  de\right)
+f_{x}\left(  t,s,x,\hat{\varphi}^{s}\left(  z\right)  \right)  \text{, }\\
\text{for }\left(  s,x,z\right)  \in\left[  t,T\right]  \times\mathbb{R}%
\times\mathbb{R}\text{,}%
\end{array}
\\
\hat{\theta}_{x}^{t}\left(  T,x,z\right)  =F_{x}\left(  t,x\right)  \text{,
for }\left(  x,z\right)  \in\mathbb{R}\times\mathbb{R}%
\end{array}
\right.
\]
and%
\[
\left\{
\begin{array}
[c]{l}%
0=\hat{g}_{sx}\left(  s,x,z\right)  +\hat{g}_{xx}\left(  s,x,z\right)
\mu\left(  s,x,\hat{\varphi}^{s}\left(  z\right)  \right)  +\hat{g}_{x}\left(
s,x,z\right)  \mu_{x}\left(  s,x,\hat{\varphi}^{s}\left(  z\right)  \right) \\
\text{ \ }+\hat{g}_{zx}\left(  s,x,z\right)  \mu\left(  s,z,\hat{\varphi}%
^{s}\left(  z\right)  \right)  +\frac{1}{2}\sigma\left(  s,z,\hat{\varphi}%
^{s}\left(  z\right)  \right)  ^{2}\hat{g}_{zzx}\left(  s,x,z\right) \\
\text{ \ }+\sigma\left(  s,x,\hat{\varphi}^{s}\left(  z\right)  \right)
\sigma_{x}\left(  s,x,\hat{\varphi}^{s}\left(  z\right)  \right)  \hat{g}%
_{xx}\left(  s,x,z\right)  +\frac{1}{2}\sigma\left(  s,x,\hat{\varphi}%
^{s}\left(  z\right)  \right)  ^{2}\hat{g}_{xxx}\left(  s,x,z\right) \\
\text{ \ }+\sigma_{x}\left(  s,x,\hat{\varphi}^{s}\left(  z\right)  \right)
\sigma\left(  s,z,\hat{\varphi}^{s}\left(  z\right)  \right)  \hat{g}%
_{xz}\left(  s,x,z\right)  +\sigma\left(  s,x,\hat{\varphi}^{s}\left(
z\right)  \right)  \sigma\left(  s,z,\hat{\varphi}^{s}\left(  z\right)
\right)  \hat{g}_{xxz}\left(  s,x,z\right) \\
\text{ \ }+%
{\displaystyle\int_{E}}
\left\{  \hat{g}_{x}\left(  s,x+c\left(  s,x,\hat{\varphi}^{s}\left(
z\right)  ,e\right)  ,z+c\left(  s,z,\hat{\varphi}^{s}\left(  z\right)
,e\right)  \right)  \left(  1+c_{x}\left(  s,x,\hat{\varphi}^{s}\left(
z\right)  ,e\right)  \right)  \right. \\
\text{ \ \ }-\hat{g}_{x}\left(  s,x,z\right)  -\hat{g}_{xx}\left(
s,x,z\right)  c\left(  s,x,\hat{\varphi}^{s}\left(  z\right)  ,e\right)
-\hat{g}_{x}\left(  s,x,z\right)  c_{x}\left(  s,x,\hat{\varphi}^{s}\left(
z\right)  ,e\right) \\
\text{ \ \ \ }%
\begin{array}
[c]{r}%
\left.  -\hat{g}_{zx}\left(  s,x,z\right)  c\left(  s,z,\hat{\varphi}%
^{s}\left(  z\right)  ,e\right)  \right\}  \vartheta\left(  de\right)  \text{,
}\\
\text{for }\left(  s,x,z\right)  \in\left[  0,T\right]  \times\mathbb{R}%
\times\mathbb{R}\text{,}%
\end{array}
\\
\hat{g}_{x}\left(  T,x,z\right)  =1\text{, for }\left(  x,z\right)
\in\mathbb{R}\times\mathbb{R}\text{.}%
\end{array}
\right.
\]

Then applying It\^{o} formula to $\theta_{x}^{t}\left(  \cdot,\hat{X}^{x_{0}%
}\left(  \cdot\right)  ,\hat{X}^{x_{0}}\left(  \cdot\right)  \right)  $,
$g_{x}\left(  \cdot,\hat{X}^{x_{0}}\left(  \cdot\right)  ,\hat{X}^{x_{0}%
}\left(  \cdot\right)  \right)  $ and using the above IPDEs, we get%
\begin{equation}
\left\{
\begin{array}
[c]{l}%
d\hat{\theta}_{x}^{t}\left(  s\right)  =-\left\{  \mu_{x}\left(  s\right)
\hat{\theta}^{t}\left(  s\right)  +\sigma_{x}\left(  s\right)  \sigma\left(
s\right)  \left(  \hat{\theta}_{xx}^{t}\left(  s\right)  +\hat{\theta}%
_{xz}\left(  s\right)  \right)  \right. \\
\text{ \ \ \ \ \ \ \ \ \ \ \ \ \ }\left.  +\int_{E}c_{x}\left(  s,e\right)
\Delta\hat{\theta}_{x}\left(  s,e\right)  \vartheta\left(  de\right)
+f_{x}\left(  t,s\right)  \right\}  ds\\
\text{ \ \ \ \ \ \ \ \ \ \ \ \ }+\sigma\left(  s\right)  \left(  \hat{\theta
}_{xx}^{t}\left(  s\right)  +\hat{\theta}_{xz}\left(  s\right)  \right)
dW\left(  s\right)  +\int_{E}\Delta\hat{\theta}_{x}\left(  s,e\right)
\tilde{N}\left(  ds,de\right)  \text{, }s\in\left[  t,T\right]  \text{,}\\
\hat{\theta}_{x}^{t}\left(  T\right)  =F_{x}\left(  t,\hat{X}^{x_{0}}\left(
T\right)  \right)
\end{array}
\right.  \label{Eq78}%
\end{equation}
and%
\begin{equation}
\left\{
\begin{array}
[c]{l}%
d\hat{g}_{x}\left(  s\right)  =-\left\{  \mu_{x}\left(  s\right)  \hat
{g}\left(  s\right)  +\sigma_{x}\left(  s\right)  \sigma\left(  s\right)
\left(  \hat{g}_{xx}\left(  s\right)  +\hat{g}_{xz}\left(  s\right)  \right)
\right. \\
\text{ \ \ \ \ \ \ \ \ \ \ \ \ \ }\left.  +\int_{E}c_{x}\left(  s,e\right)
\Delta\hat{g}_{x}\left(  s,e\right)  \vartheta\left(  de\right)  \right\}
ds\\
\text{ \ \ \ \ \ \ \ \ \ \ \ \ }+\sigma\left(  s\right)  \left(  \hat{g}%
_{xx}\left(  s\right)  +\hat{g}_{xz}\left(  s\right)  \right)  dW\left(
s\right)  +\int_{E}\Delta\hat{g}_{x}\left(  s,e\right)  \tilde{N}\left(
ds,de\right)  \text{, }s\in\left[  t,T\right]  \text{,}\\
\hat{g}_{x}\left(  T\right)  =1\text{.}%
\end{array}
\right.  \label{Eq79}%
\end{equation}

Now define $\left(  \bar{p}^{t}\left(  \cdot\right)  ,\bar{q}^{t}\left(
\cdot\right)  ,\bar{r}^{t}\left(  \cdot,\cdot\right)  \right)  $ as in
(\ref{Eq76}). Using (\ref{Eq78})-(\ref{Eq79}), we can easily verify that the
3-tuple of processes $\left(  \bar{p}^{t}\left(  \cdot\right)  ,\bar{q}%
^{t}\left(  \cdot\right)  ,\bar{r}^{t}\left(  \cdot,\cdot\right)  \right)  $
satisfies the following BSDE%
\[
\left\{
\begin{array}
[c]{l}%
d\bar{p}^{t}\left(  s\right)  =-\left\{  \mu_{x}\left(  s\right)  \bar{p}%
^{t}\left(  s\right)  +\sigma_{x}\left(  s\right)  \bar{q}^{t}\left(
s\right)  +\int_{E}c_{x}\left(  s,e\right)  \bar{r}^{t}\left(  s,e\right)
\vartheta\left(  de\right)  +f_{x}\left(  t,s\right)  \right\}  ds\\
\text{ \ \ \ \ \ \ \ \ \ \ \ \ }+\bar{q}^{t}\left(  s\right)  dW\left(
s\right)  +\int_{E}\bar{r}\left(  s,e\right)  \tilde{N}\left(  ds,de\right)
\text{, }s\in\left[  t,T\right]  \text{,}\\
\bar{p}^{t}\left(  T\right)  =F_{x}\left(  t,\hat{X}^{x_{0}}\left(  T\right)
\right)  +G_{\bar{x}}\left(  t,\hat{X}^{x_{0}}\left(  t\right)  ,\mathbb{E}%
_{t}\left[  \hat{X}^{x_{0}}\left(  T\right)  \right]  \right)  \text{.}%
\end{array}
\right.
\]
Hence, by the uniqueness of the solution to (\ref{Eq72}), we obtain that%
\[
\left(  \bar{p}^{t}\left(  \cdot\right)  ,\bar{q}^{t}\left(  \cdot\right)
,\bar{r}^{t}\left(  \cdot,\cdot\right)  \right)  =\left(  p^{t}\left(
\cdot\right)  ,q^{t}\left(  \cdot\right)  ,r^{t}\left(  \cdot,\cdot\right)
\right)  ,\text{ a.s., a.e. }\left[  s,e\right]  \in\left[  t,T\right]  \times
E\text{.}%
\]

Similarly, we can easily verify that the 3-tuple of processes defined by
(\ref{Eq77}) coincides with the solution of BSDE (\ref{Eq73}). This completes
the proof.\eop

\section{A Mean-variance portfolio problem \label{section6}}

As an application of the general theory, we consider a Markowitz mean-variance
portfolio selection problem\ associated to a jump-diffusion model with
deterministic coefficients. We apply the verification argument in Theorem
\ref{result2} to derive the equilibrium investment strategy in an explicit
form. Note that, this problem was discussed in \cite{Sun}, where the authors
derived the open-loop equilibrium solution by solving a\ flow of FBSDEs with
jumps. We emphasize that, we solve this problem for illustrative purposes
only. Our goal is to show that the solution we obtain via our verification
theorem coincides with the one derived in \cite{Huetal2011}; which indirectly
indicates that our PDEs-approach is correct.

In this section, we assume that $d=1$ (i.e., the Brownian motion is
one-dimensional). Suppose that there is a financial market, in which two
securities are traded continuously. One of them is a bond, with price
$S_{0}\left(  s\right)  $ at time $s\in\left[  0,T\right]  $ governed by%
\[
\frac{dS_{0}\left(  s\right)  }{S_{0}\left(  s\right)  }=r_{0}\left(
s\right)  ds\text{, }S_{0}\left(  0\right)  =s_{0}>0\text{,}%
\]
where $r_{0}:\left[  0,T\right]  \rightarrow\left(  0,+\infty\right)  $ is a
deterministic function which represents the risk-free rate. The other asset is
called the risky stock, whose price process $S_{1}\left(  \cdot\right)  $
satisfies the following stochastic differential equation%
\[
\dfrac{dS_{1}\left(  s\right)  }{S_{1}\left(  s\right)  }=r\left(  s\right)
ds+\sigma\left(  s\right)  dW\left(  s\right)  +%
{\displaystyle\int_{E}}
\phi\left(  s,e\right)  \tilde{N}\left(  ds,de\right)  \text{, }S_{1}\left(
0\right)  =s_{1}>0\text{,}%
\]
where $r:\left[  0,T\right]  \rightarrow\left(  0,+\infty\right)  $,
$\sigma:\left[  0,T\right]  \rightarrow\left(  0,+\infty\right)  $ and
$\phi:\left[  0,T\right]  \times E\rightarrow%
\mathbb{R}
$ are deterministic measurable functions; $r\left(  \cdot\right)  $,
$\sigma\left(  \cdot\right)  $ and $\phi\left(  \cdot,\cdot\right)  $
represent the appreciation rate, the volatility and the jump coefficient of
the risky stock, respectively.

We assume that $r_{0}\left(  \cdot\right)  $, $r\left(  \cdot\right)  $,
$\sigma\left(  \cdot\right)  $ and $\phi\left(  \cdot,\cdot\right)  $ are
continuous and uniformly bounded functions, such that $r\left(  s\right)
>r_{0}\left(  s\right)  $ and $\phi\left(  s,e\right)  \geq-1$, for all
$s\in\left[  0,T\right]  $ and $\vartheta-$a.e. $e\in E$. We also require a
uniform ellipticity condition as follows:%
\[
\sigma\left(  s\right)  ^{2}+%
{\displaystyle\int_{E}}
\phi\left(  s,e\right)  ^{2}\theta\left(  de\right)  \geq\epsilon\text{, for
all }s\in\left[  0,T\right]  \text{,}%
\]
for some $\epsilon>0$\textit{.}

Starting from an initial capital $x_{0}>0$, the investor is allowed to invest
in the financial market. A trading strategy is a one-dimensional stochastic
process $u\left(  \cdot\right)  $, where $u\left(  s\right)  $ represents the
amount invested in the risky stock at time $s\in\left[  0,T\right]  $. The
dollar amount invested in the bond at time $s$ is $X^{x_{0},u\left(  .\right)
}\left(  s\right)  -u\left(  s\right)  $, where $X^{x_{0},u\left(  .\right)
}\left(  \cdot\right)  $ is the wealth process associated with the strategy
$u\left(  \cdot\right)  $ and the initial capital $x_{0}$. Then the evolution
of $\ X^{x_{0},u\left(  .\right)  }\left(  \cdot\right)  $ can be described as%
\[
\left\{
\begin{array}
[c]{l}%
dX^{x_{0},u\left(  .\right)  }\left(  s\right)  =\left(  r_{0}\left(
s\right)  X^{x_{0},u\left(  .\right)  }\left(  s\right)  +u\left(  s\right)
\rho\left(  s\right)  \right)  ds+u\left(  s\right)  \sigma\left(  s\right)
dW\left(  s\right) \\
\text{ \ \ \ \ \ \ \ \ \ \ \ \ \ \ \ \ \ \ \ \ \ }+u\left(  s\right)
{\displaystyle\int_{E}}
\phi\left(  s,e\right)  \tilde{N}\left(  ds,de\right)  \text{, for }%
s\in\left[  0,T\right]  \text{,}\\
X^{x_{0},u\left(  .\right)  }\left(  0\right)  =x_{0}\text{,}%
\end{array}
\right.
\]
where $\rho\left(  \cdot\right)  :=r\left(  \cdot\right)  -r_{0}\left(
\cdot\right)  $.

when time evolves, we consider the controlled stochastic differential equation
parameterized by $\left(  t,y\right)  \in\left[  0,T\right]  \times%
\mathbb{R}
$,%
\begin{equation}
\left\{
\begin{array}
[c]{l}%
dX\left(  s\right)  =\left\{  r_{0}\left(  s\right)  X\left(  s\right)
+u\left(  s\right)  \rho\left(  s\right)  \right\}  ds+u\left(  s\right)
\sigma\left(  s\right)  dW\left(  s\right) \\
\text{ \ \ \ \ \ \ \ \ \ \ \ \ \ }+%
{\displaystyle\int_{E}}
u\left(  s\right)  \phi\left(  s,z\right)  \tilde{N}\left(  ds,dz\right)
\text{, for }s\in\left[  t,T\right]  \text{,}\\
X\left(  t\right)  =y\text{.}%
\end{array}
\right.  \label{Eq39}%
\end{equation}

A trading strategy $u\left(  \cdot\right)  $ is said to be admissible over
$\left[  t,T\right]  $, if it is a $%
\mathbb{R}
$-valued $\left(  \mathcal{F}_{s}\right)  _{s\in\left[  t,T\right]  }%
$-predictable process such that:%
\[
\mathbb{E}\left[  \int_{t}^{T}\left\vert u\left(  s\right)  \right\vert
^{4}ds\right]  <\infty\text{.}%
\]

For any fixed initial capital $\left(  t,y\right)  \in\left[  0,T\right]
\times%
\mathbb{R}
$, the investor's aim is to choose an investment strategy $u\left(
\cdot\right)  $ in order to maximize the conditional expectation of the
terminal wealth over the period $\left[  t,T\right]  $, while trying at the
same time to minimize the conditional variance of the terminal wealth (i.e.
$\mathbb{V}\mathrm{ar}_{t}\left[  X\left(  T\right)  \right]  =\mathbb{E}%
_{t}\left[  X\left(  T\right)  ^{2}\right]  -\mathbb{E}_{t}\left[  X\left(
T\right)  \right]  ^{2}$). The optimization problem is therefore to minimize%
\begin{align}
\mathbf{J}\left(  t,y;u\left(  \cdot\right)  \right)   &  :=\dfrac{1}%
{2}\mathbb{V}\mathrm{ar}_{t}\left[  X\left(  T\right)  \right]  -\mu
y\mathbb{E}_{t}\left[  X\left(  T\right)  \right] \nonumber\\
&  :=\dfrac{1}{2}\mathbb{E}_{t}\left[  X\left(  T\right)  ^{2}\right]
-\left(  \dfrac{1}{2}\mathbb{E}_{t}\left[  X\left(  T\right)  \right]
^{2}+\mu y\mathbb{E}_{t}\left[  X\left(  T\right)  \right]  \right)
\label{Eq40}%
\end{align}
over $\mathcal{L}_{%
\mathcal{F}%
,p}^{4}\left(  t,T;%
\mathbb{R}
\right)  $, where $\mu y$, with $\mu\geq0$, denotes the weight between the
conditional variance and the conditional expectation.

It is easy to see that, the above model is mathematically a special case of
Problem (N) formulated earlier in this paper, with%
\begin{align*}
n  &  =d=m=1,\text{ }U=\mathbb{R}\text{,}\\
\sigma\left(  s,x,u\right)   &  \equiv u\sigma\left(  s\right)  ,\text{
}c\left(  s,x,u,e\right)  \equiv u\phi\left(  s,e\right)  ,\\
\mu\left(  s,x,u\right)   &  \equiv r_{0}\left(  s\right)  x+u\rho\left(
s\right)  ,\\
\text{ }f\left(  t,y,s,x,u\right)   &  \equiv0,\text{ }F\left(  t,y,x\right)
\equiv\dfrac{1}{2}x^{2},\\
\text{ }\Psi\left(  x\right)   &  =x\text{ and }G\left(  t,y,\bar{x}\right)
\equiv-\dfrac{1}{2}\bar{x}^{2}-\mu y\bar{x}.
\end{align*}
Accordingly, the IPDEs associated to an admissible strategy $\hat{\varphi
}\left(  \cdot,\cdot\right)  $ are defined as follows:%
\[
\left\{
\begin{array}
[c]{l}%
0=\hat{\theta}_{s}\left(  s,x,z\right)  +\hat{\theta}_{x}\left(  s,x,z\right)
\left(  r_{0}\left(  s\right)  x+\hat{\varphi}^{s}\left(  z\right)
\rho\left(  s\right)  \right) \\
\text{ \ \ }+\hat{\theta}_{z}\left(  s,x,z\right)  \left(  r_{0}\left(
s\right)  z+\hat{\varphi}^{s}\left(  z\right)  \rho\left(  s\right)  \right)
\\
\text{ \ \ }+\frac{1}{2}\left\vert \hat{\varphi}^{s}\left(  z\right)
\sigma\left(  s\right)  \right\vert ^{2}\left(  \hat{\theta}_{xx}\left(
s,x,z\right)  +\hat{\theta}_{zz}\left(  s,x,z\right)  +2\hat{\theta}%
_{xz}\left(  s,x,z\right)  \right) \\
\text{ \ \ }+%
{\displaystyle\int\limits_{E}}
\left\{  \hat{\theta}\left(  s,x+\hat{\varphi}^{s}\left(  z\right)
\phi\left(  s,e\right)  ,z+\hat{\varphi}^{s}\left(  z\right)  \phi\left(
s,e\right)  \right)  -\hat{\theta}\left(  s,x,z\right)  \right. \\
\text{ \ \ }%
\begin{array}
[c]{r}%
\text{\ }\left.  -\left(  \hat{\theta}_{x}\left(  s,x,z\right)  +\hat{\theta
}_{z}\left(  s,x,z\right)  \right)  \hat{\varphi}^{s}\left(  z\right)
\phi\left(  s,e\right)  \right\}  \vartheta\left(  de\right)  \text{, }\\
\text{for }\left(  s,x,z\right)  \in\left[  0,T\right]  \times%
\mathbb{R}
\times%
\mathbb{R}
\text{,}%
\end{array}
\\
\hat{\theta}\left(  T,x,z\right)  =\frac{1}{2}x^{2}\text{, for }\left(
x,z\right)  \in%
\mathbb{R}
\times%
\mathbb{R}
,
\end{array}
\right.
\]%
\[
\left\{
\begin{array}
[c]{l}%
0=\hat{g}_{s}\left(  s,x,z\right)  +\hat{g}_{x}\left(  s,x,z\right)  \left(
r_{0}\left(  s\right)  x+\hat{\varphi}^{s}\left(  z\right)  \rho\left(
s\right)  \right) \\
\text{ \ \ }+\hat{g}_{z}\left(  s,x,z\right)  \left(  r_{0}\left(  s\right)
z+\hat{\varphi}^{s}\left(  z\right)  \rho\left(  s\right)  \right) \\
\text{ \ \ }+\frac{1}{2}\left\vert \hat{\varphi}^{s}\left(  z\right)
\sigma\left(  s\right)  \right\vert ^{2}\left(  \hat{g}_{xx}\left(
s,x,z\right)  +\hat{g}_{zz}\left(  s,x,z\right)  +2\hat{g}_{xz}\left(
s,x,z\right)  \right) \\
\text{ \ \ }+%
{\displaystyle\int\limits_{E}}
\left\{  \hat{g}\left(  s,x+\hat{\varphi}^{s}\left(  z\right)  \phi\left(
s,e\right)  ,z+\hat{\varphi}^{s}\left(  z\right)  \phi\left(  s,e\right)
\right)  -\hat{g}\left(  s,x,z\right)  \right. \\
\text{ \ \ }%
\begin{array}
[c]{r}%
\text{\ }\left.  -\left(  \hat{g}_{x}\left(  s,x,z\right)  +\hat{g}_{z}\left(
s,x,z\right)  \right)  \hat{\varphi}^{s}\left(  z\right)  \phi\left(
s,e\right)  \right\}  \vartheta\left(  de\right)  \text{, }\\
\text{for }\left(  s,x,z\right)  \in\left[  0,T\right]  \times%
\mathbb{R}
\times%
\mathbb{R}
\text{,}%
\end{array}
\\
\hat{g}\left(  T,x,z\right)  =x\text{, for }\left(  x,z\right)  \in%
\mathbb{R}
\times%
\mathbb{R}
,
\end{array}
\right.
\]
and the $\mathcal{H}$-function associated to $\left(  \hat{\varphi}\left(
\cdot,\cdot\right)  ,\hat{\theta}\left(  \cdot,\cdot,\cdot\right)  ,\hat
{g}\left(  \cdot,\cdot,\cdot\right)  \right)  $ defined in (\ref{Eq25}) takes
the form,%
\begin{align*}
&  \mathcal{H}\left(  t,y,s,X,Z,u\right) \\
&  :=\frac{1}{2}\left(  \hat{\theta}_{xx}\left(  s,X,Z\right)  -\left(  \mu
y+\mathbb{E}_{t}\left[  \hat{g}\left(  s,X,Z\right)  \right]  \right)  \hat
{g}_{xx}\left(  s,X,Z\right)  \right)  \left(  \sigma\left(  s\right)
u\right)  ^{2}\\
&  +\left(  \hat{\theta}_{x}\left(  s,X,Z\right)  -\left(  \mu y+\mathbb{E}%
_{t}\left[  \hat{g}\left(  s,X,Z\right)  \right]  \right)  \hat{g}_{x}\left(
s,X,Z\right)  \right)  \left(  r_{0}\left(  s\right)  x+u\rho\left(  s\right)
\right) \\
&  +\left(  \hat{\theta}_{xz}\left(  s,X,Z\right)  -\left(  \mu y+\mathbb{E}%
_{t}\left[  \hat{g}\left(  s,X,Z\right)  \right]  \right)  \hat{g}_{xz}\left(
s,X,Z\right)  \right)  u\hat{\varphi}^{s}\left(  z\right)  \sigma\left(
s\right)  ^{2}\\
&  +%
{\displaystyle\int\limits_{E}}
\left\{  \hat{\theta}\left(  s,X+u\phi\left(  s,e\right)  ,z+\hat{\varphi}%
^{s}\left(  Z\right)  \phi\left(  s,e\right)  \right)  -\hat{\theta}%
_{x}\left(  s,X,Z\right)  u\phi\left(  s,e\right)  \right\}  \vartheta\left(
de\right) \\
&  -\left(  \mu y+\mathbb{E}_{t}\left[  \hat{g}\left(  s,X,Z\right)  \right]
\right)
{\displaystyle\int\limits_{E}}
\left\{  \hat{g}\left(  s,X+u\phi\left(  s,e\right)  ,z+\hat{\varphi}%
^{s}\left(  Z\right)  \phi\left(  s,e\right)  \right)  \right. \\
&  \left.  -\hat{g}_{x}\left(  s,X,Z\right)  u\phi\left(  s,e\right)
\right\}  \vartheta\left(  de\right)  \text{.}%
\end{align*}

\subsection{Equilibrium solution}

In the next, we derive the equilibrium investment strategy in an explicit
form. First, letting for $s\in\left[  0,T\right]  $,%
\begin{equation}
\kappa\left(  s\right)  :=\frac{\rho\left(  s\right)  }{\left(  \sigma\left(
s\right)  ^{2}+%
{\displaystyle\int\limits_{E}}
\phi\left(  s,e\right)  ^{2}\vartheta\left(  de\right)  \right)  }.
\label{Eq40*}%
\end{equation}

\begin{theorem}
The mean-variance portfolio problem in (\ref{Eq39})-(\ref{Eq40}) has an
open-loop equilibrium strategy that can be represented by%
\begin{equation}%
\begin{array}
[c]{r}%
\hat{\varphi}^{s}\left(  z\right)  =\mu e^{-\int_{s}^{T}r_{0}\left(
\tau\right)  d\tau}\left(  1+\mu\int_{s}^{T}e^{-\int_{\tau}^{T}r_{0}\left(
v\right)  dv}\kappa\left(  \tau\right)  \rho\left(  \tau\right)  d\tau\right)
^{-1}\kappa\left(  s\right)  z\text{,}\\
\text{for all }\left(  s,z\right)  \in\left[  0,T\right]  \times%
\mathbb{R}
\text{.}%
\end{array}
\label{Eq41}%
\end{equation}

\end{theorem}

\bop Assume for the time being that the conditions of Theorem \ref{result2}
hold. The verification argument in Theorem \ref{result2} leads to the
following system of IPDEs,%
\begin{equation}
\left\{
\begin{array}
[c]{l}%
0=\hat{\theta}_{s}\left(  s,x,z\right)  +\hat{\theta}_{x}\left(  s,x,z\right)
\left(  r_{0}\left(  s\right)  x+\hat{\varphi}^{s}\left(  z\right)
\rho\left(  s\right)  \right)  \text{\ }+\hat{\theta}_{z}\left(  s,x,z\right)
\left(  r_{0}\left(  s\right)  z+\hat{\varphi}^{s}\left(  z\right)
\rho\left(  s\right)  \right) \\
\text{ \ \ }+\frac{1}{2}\left\vert \hat{\varphi}^{s}\left(  z\right)
\sigma\left(  s\right)  \right\vert ^{2}\left(  \hat{\theta}_{xx}\left(
s,x,z\right)  +\hat{\theta}_{zz}\left(  s,x,z\right)  +2\hat{\theta}%
_{xz}\left(  s,x,z\right)  \right) \\
\text{ \ \ }+%
{\displaystyle\int\limits_{E}}
\left\{  \hat{\theta}\left(  s,x+\hat{\varphi}^{s}\left(  z\right)
\phi\left(  s,e\right)  ,z+\hat{\varphi}^{s}\left(  z\right)  \phi\left(
s,e\right)  \right)  -\hat{\theta}\left(  s,x,z\right)  \right. \\
\text{ \ \ \ }\left.  -\left(  \hat{\theta}_{x}\left(  s,x,z\right)
+\hat{\theta}_{z}\left(  s,x,z\right)  \right)  \hat{\varphi}^{s}\left(
z\right)  \phi\left(  s,e\right)  \right\}  \vartheta\left(  de\right)
\text{, }\\
0=\hat{g}_{s}\left(  s,x,z\right)  +\hat{g}_{x}\left(  s,x,z\right)  \left(
r_{0}\left(  s\right)  x+\hat{\varphi}^{s}\left(  z\right)  \rho\left(
s\right)  \right)  +\hat{g}_{z}\left(  s,x,z\right)  \left(  r_{0}\left(
s\right)  z+\hat{\varphi}^{s}\left(  z\right)  \rho\left(  s\right)  \right)
\\
\text{ \ \ }+\frac{1}{2}\left\vert \hat{\varphi}^{s}\left(  z\right)
\sigma\left(  s\right)  \right\vert ^{2}\left(  \hat{g}_{xx}\left(
s,x,z\right)  +\hat{g}_{zz}\left(  s,x,z\right)  +2\hat{g}_{xz}\left(
s,x,z\right)  \right) \\
\text{ \ \ }+%
{\displaystyle\int\limits_{E}}
\left\{  \hat{g}\left(  s,x+\hat{\varphi}^{s}\left(  z\right)  \phi\left(
s,e\right)  ,z+\hat{\varphi}^{s}\left(  z\right)  \phi\left(  s,e\right)
\right)  -\hat{g}\left(  s,x,z\right)  \right. \\
\text{ \ \ }%
\begin{array}
[c]{r}%
\text{\ }\left.  -\left(  \hat{g}_{x}\left(  s,x,z\right)  +\hat{g}_{z}\left(
s,x,z\right)  \right)  \hat{\varphi}^{s}\left(  z\right)  \phi\left(
s,e\right)  \right\}  \vartheta\left(  de\right)  \text{, }\\
\text{for }\left(  s,x,z\right)  \in\left[  0,T\right]  \times%
\mathbb{R}
\times%
\mathbb{R}
\text{,}%
\end{array}
\\
\hat{\theta}\left(  T,x,z\right)  =\frac{1}{2}x^{2}\text{, }\hat{g}\left(
T,x,z\right)  =x\text{, for }\left(  x,z\right)  \in%
\mathbb{R}
\times%
\mathbb{R}
\text{,}%
\end{array}
\right.  \label{Eq42}%
\end{equation}
with the condition: for all $\left(  s,z\right)  \in\left[  0,T\right]
\times\mathbb{R}_{+}$,%
\begin{equation}
\mathcal{H}\left(  s,z,s,z,z,\hat{\varphi}^{s}\left(  z\right)  \right)
=\min_{u\in%
\mathbb{R}
}\mathcal{H}\left(  s,z,s,z,z,u\right)  \text{.} \label{Eq43}%
\end{equation}

To solve the above system of IPDEs, we suggest the following ansatz: For all
$\left(  s,x,z\right)  \in\left[  0,T\right]  \times%
\mathbb{R}
\times%
\mathbb{R}
$,%
\begin{equation}
\hat{\theta}\left(  s,x,z\right)  =M_{1}\left(  s\right)  \frac{x^{2}}%
{2}+M_{2}\left(  s\right)  \frac{z^{2}}{2}+M_{3}\left(  s\right)  zx\nonumber
\end{equation}
and%
\[
\hat{g}\left(  s,x,z\right)  =N_{1}\left(  s\right)  x+N_{2}\left(  s\right)
z\text{,}%
\]
where $M_{1}\left(  \cdot\right)  ,$ $M_{2}\left(  \cdot\right)  ,$
$M_{3}\left(  \cdot\right)  ,$ $N_{1}\left(  \cdot\right)  $ and $N_{2}\left(
\cdot\right)  $ are deterministic continuously differentiable functions of
time such that%
\[
M_{1}\left(  T\right)  =N_{1}\left(  T\right)  =1\text{ and }M_{2}\left(
T\right)  =M_{3}\left(  T\right)  =N_{2}\left(  T\right)  =0\text{.}%
\]

In this case, the partial derivatives of $\hat{\theta}\left(  s,x,z\right)  $
and $\hat{g}\left(  s,x,z\right)  $ are%
\begin{align*}
\hat{\theta}_{s}\left(  s,x,z\right)   &  =\frac{dM_{1}\left(  s\right)  }%
{ds}\frac{x^{2}}{2}+\frac{dM_{2}\left(  s\right)  }{ds}\frac{z^{2}}{2}%
+\frac{dM_{3}\left(  s\right)  }{ds}zx\text{,}\\
\hat{g}_{s}\left(  s,x,z\right)   &  =\frac{dN_{1}\left(  s\right)  }%
{ds}x+\frac{dN_{2}\left(  s\right)  }{ds}z\text{,}\\
\hat{\theta}_{x}\left(  s,x,z\right)   &  =M_{1}\left(  s\right)
x+M_{3}\left(  s\right)  z\text{, }\hat{\theta}_{z}\left(  s,x,z\right)
=M_{2}\left(  s\right)  z+M_{3}\left(  s\right)  x\text{,}\\
\hat{\theta}_{xx}\left(  s,x,z\right)   &  =M_{1}\left(  s\right)  \text{,
}\hat{\theta}_{zz}\left(  s,x\right)  =M_{2}\left(  s\right)  \text{,}\\
\hat{\theta}_{xz}\left(  s,x,z\right)   &  =M_{3}\left(  s\right)  \text{,
}\hat{g}_{x}\left(  s,x,z\right)  =N_{1}\left(  s\right)  \text{, }\hat{g}%
_{z}\left(  s,x,z\right)  =N_{2}\left(  s\right)  \text{ and}\\
\hat{g}_{xx}\left(  s,x,z\right)   &  =\hat{g}_{zz}\left(  s,x,z\right)
=\hat{g}_{xz}\left(  s,x,z\right)  =0\text{.}%
\end{align*}
Thus, by substituting $\hat{\theta}\left(  s,x,z\right)  $ and $\hat{g}\left(
s,x,z\right)  $ together with the above derivatives into (\ref{Eq42}), this
leads to%
\begin{align}
0  &  =\frac{dM_{1}}{ds}\left(  s\right)  \frac{x^{2}}{2}+\frac{dM_{2}}%
{ds}\left(  s\right)  \frac{z^{2}}{2}+\frac{dM_{3}}{ds}zx\nonumber\\
&  +\left(  M_{1}\left(  s\right)  x+M_{3}\left(  s\right)  z\right)  \left(
r_{0}\left(  s\right)  x+\hat{\varphi}^{s}\left(  z\right)  \rho\left(
s\right)  \right) \nonumber\\
&  +\left(  M_{2}\left(  s\right)  z+M_{3}\left(  s\right)  x\right)  \left(
r_{0}\left(  s\right)  z+\hat{\varphi}^{s}\left(  z\right)  \rho\left(
s\right)  \right) \nonumber\\
&  +\frac{1}{2}\left\vert \hat{\varphi}^{s}\left(  z\right)  \sigma\left(
s\right)  \right\vert ^{2}\left(  M_{1}\left(  s\right)  +M_{2}\left(
s\right)  +2M_{3}\left(  s\right)  \right) \nonumber\\
&  +%
{\displaystyle\int\limits_{E}}
\left\{  M_{1}\left(  s\right)  \frac{\left(  x+\hat{\varphi}^{s}\left(
z\right)  \phi\left(  s,e\right)  \right)  ^{2}}{2}+M_{2}\left(  s\right)
\frac{\left(  z+\hat{\varphi}^{s}\left(  z\right)  \phi\left(  s,e\right)
\right)  ^{2}}{2}\right. \nonumber\\
&  +M_{3}\left(  s\right)  \left(  x+\hat{\varphi}^{s}\left(  z\right)
\phi\left(  s,e\right)  \right)  \left(  z+\hat{\varphi}^{s}\left(  z\right)
\phi\left(  s,e\right)  \right) \nonumber\\
&  -M_{1}\left(  s\right)  \frac{x^{2}}{2}-M_{2}\left(  s\right)  \frac{z^{2}%
}{2}-M_{3}\left(  s\right)  zx\nonumber\\
&  \left.  -\phi\left(  s,e\right)  \hat{\varphi}^{s}\left(  z\right)  \left(
\left(  M_{1}\left(  s\right)  +M_{3}\left(  s\right)  \right)  x+\left(
M_{3}\left(  s\right)  +M_{2}\left(  s\right)  \right)  z\right)  \right\}
\vartheta\left(  de\right)  \label{Eq44}%
\end{align}
and%
\begin{align}
0  &  =\frac{dN_{1}}{ds}x+\frac{dN_{2}}{ds}z\nonumber\\
&  +N_{1}\left(  s\right)  \left(  r_{0}\left(  s\right)  x+\hat{\varphi}%
^{s}\left(  z\right)  \rho\left(  s\right)  \right)  +N_{2}\left(  s\right)
\left(  r_{0}\left(  s\right)  z+\hat{\varphi}^{s}\left(  z\right)
\rho\left(  s\right)  \right) \nonumber\\
&  +%
{\displaystyle\int\limits_{E}}
\left\{  N_{1}\left(  s\right)  \left(  x+\hat{\varphi}^{s}\left(  z\right)
\phi\left(  s,e\right)  \right)  +N_{2}\left(  s\right)  \left(
z+\hat{\varphi}^{s}\left(  z\right)  \phi\left(  s,e\right)  \right)  \right.
\nonumber\\
&  \left.  -N_{1}\left(  s\right)  x-N_{2}\left(  s\right)  z-\left(
N_{1}\left(  s\right)  +N_{2}\left(  s\right)  \right)  \hat{\varphi}%
^{s}\left(  z\right)  \phi\left(  s,e\right)  \right\}  \vartheta\left(
de\right)  . \label{Eq45}%
\end{align}

On the other hand, the minimum condition in (\ref{Eq43}) suggests that: for
all $\left(  s,z\right)  \in\left[  0,T\right]  \times%
\mathbb{R}
$,%
\begin{align*}
&  0=\mathcal{H}_{u}\left(  s,z,s,z,z,\hat{\varphi}^{s}\left(  z\right)
\right) \\
&  =\left(  \hat{\theta}_{xx}\left(  s,z,z\right)  -\left(  \mu z+\hat
{g}\left(  s,z,z\right)  \right)  \hat{g}_{xx}\left(  s,z,z\right)  \right)
\sigma\left(  s\right)  ^{2}\hat{\varphi}^{s}\left(  z\right) \\
&  +\left(  \hat{\theta}_{x}\left(  s,z,z\right)  -\left(  \mu z+\hat
{g}\left(  s,z,z\right)  \right)  \hat{g}_{x}\left(  s,z,z\right)  \right)
\rho\left(  s\right) \\
&  +\left(  \hat{\theta}_{xz}\left(  s,z,z\right)  -\left(  \mu z+\hat
{g}\left(  s,z,z\right)  \right)  \hat{g}_{xz}\left(  s,z,z\right)  \right)
\hat{\varphi}^{s}\left(  z\right)  \sigma\left(  s\right)  ^{2}\\
&  +%
{\displaystyle\int\limits_{E}}
\left\{  \hat{\theta}_{x}\left(  s,z+\hat{\varphi}^{s}\left(  z\right)
\phi\left(  s,e\right)  ,z+\hat{\varphi}^{s}\left(  z\right)  \phi\left(
s,e\right)  \right)  \phi\left(  s,e\right)  -\hat{\theta}_{x}\left(
s,z,z\right)  \phi\left(  s,e\right)  \right\}  \vartheta\left(  de\right) \\
&  -\left(  \mu z+\hat{g}\left(  s,z,z\right)  \right)
{\displaystyle\int\limits_{E}}
\left\{  \hat{g}_{x}\left(  s,z+\hat{\varphi}^{s}\left(  z\right)  \phi\left(
s,e\right)  ,z+\hat{\varphi}^{s}\left(  z\right)  \phi\left(  s,e\right)
\right)  \phi\left(  s,e\right)  \right. \\
&  \left.  -\hat{g}_{x}\left(  s,z,z\right)  \phi\left(  s,e\right)  \right\}
\vartheta\left(  de\right)
\end{align*}%
\begin{align*}
&  =\left(  M_{1}\left(  s\right)  +M_{3}\left(  s\right)  \right)  \left(
\sigma\left(  s\right)  ^{2}+%
{\displaystyle\int_{E}}
\phi\left(  s,e\right)  ^{2}\vartheta\left(  de\right)  \right)  \hat{\varphi
}^{s}\left(  z\right) \\
&  +\left(  M_{1}\left(  s\right)  +M_{3}\left(  s\right)  -\left(  \mu
+N_{1}\left(  s\right)  +N_{2}\left(  s\right)  \right)  N_{1}\left(
s\right)  \right)  \rho\left(  s\right)  z.
\end{align*}
Accordingly, we obtain that $\hat{\varphi}^{s}\left(  z\right)  $ admits the
following explicit representation:%
\begin{equation}
\hat{\varphi}^{s}\left(  z\right)  =-\frac{M_{1}\left(  s\right)
+M_{3}\left(  s\right)  -N_{1}\left(  s\right)  \mu-N_{1}^{2}\left(  s\right)
-N_{1}\left(  s\right)  N_{2}\left(  s\right)  }{\left(  M_{1}\left(
s\right)  +M_{3}\left(  s\right)  \right)  \left(  \sigma\left(  s\right)
^{2}+%
{\displaystyle\int_{E}}
\phi\left(  s,e\right)  ^{2}\vartheta\left(  de\right)  \right)  }\rho\left(
s\right)  z\text{,} \label{Eq46}%
\end{equation}
Denoting%
\begin{equation}
\alpha\left(  s\right)  :=-\frac{M_{1}\left(  s\right)  +M_{3}\left(
s\right)  -N_{1}\left(  s\right)  \mu-N_{1}^{2}\left(  s\right)  -N_{1}\left(
s\right)  N_{2}\left(  s\right)  }{\left(  M_{1}\left(  s\right)
+M_{3}\left(  s\right)  \right)  \left(  \sigma\left(  s\right)  ^{2}+%
{\displaystyle\int_{E}}
\phi\left(  s,e\right)  ^{2}\vartheta\left(  de\right)  \right)  }\rho\left(
s\right)  \label{Eq47}%
\end{equation}
and coming back to (\ref{Eq44})-(\ref{Eq45}), we get: For all $\left(
s,x,z\right)  \in\left[  0,T\right]  \times%
\mathbb{R}
\times%
\mathbb{R}
$%
\begin{align*}
0  &  =\frac{x^{2}}{2}\left(  \frac{dM_{1}}{ds}\left(  s\right)  +M_{1}\left(
s\right)  2r_{0}\left(  s\right)  \right) \\
&  +\frac{z^{2}}{2}\left(  \frac{dM_{2}}{ds}\left(  s\right)  +2r_{0}\left(
s\right)  M_{2}\left(  s\right)  +2\left(  M_{2}\left(  s\right)
+M_{3}\left(  s\right)  \right)  \alpha\left(  s\right)  \rho\left(  s\right)
\right. \\
&  +\left(  M_{1}\left(  s\right)  +M_{2}\left(  s\right)  +2M_{3}\left(
s\right)  \right)  \alpha\left(  s\right)  ^{2}\left(  \left\vert
\sigma\left(  s\right)  \right\vert ^{2}+\int\limits_{E}\phi\left(
s,e\right)  ^{2}\vartheta\left(  de\right)  \right) \\
&  +zx\left(  \frac{dM_{3}}{ds}+2r_{0}\left(  s\right)  M_{3}\left(  s\right)
+\left(  M_{1}\left(  s\right)  +M_{3}\left(  s\right)  \right)  \alpha\left(
s\right)  \rho\left(  s\right)  \right)
\end{align*}
and%
\begin{align*}
0  &  =\left(  \frac{dN_{1}}{ds}\left(  s\right)  +N_{1}\left(  s\right)
r_{0}\left(  s\right)  \right)  x\\
&  +\left(  \frac{dN_{2}}{ds}\left(  s\right)  +r_{0}\left(  s\right)
N_{2}\left(  s\right)  +\left(  N_{1}\left(  s\right)  +N_{2}\left(  s\right)
\right)  \alpha\left(  s\right)  \rho\left(  s\right)  \right)  z,
\end{align*}
which leads to the following systems of ODEs (suppressing $\left(  s\right)
$):%
\begin{equation}
\left\{
\begin{array}
[c]{l}%
\frac{dM_{1}}{ds}=-2r_{0}M_{1},\\
\text{ }\\
\frac{dM_{2}}{ds}=-2r_{0}M_{2}-2\left(  M_{2}+M_{3}\right)  \alpha\rho\\
\text{ \ \ }-\left(  M_{1}+M_{2}+2M_{3}\right)  \left\vert \alpha\right\vert
^{2}\left(  \left\vert \sigma\right\vert ^{2}+%
{\displaystyle\int\limits_{E}}
\phi\left(  e\right)  ^{2}\vartheta\left(  de\right)  \right)  ,\\
\text{ }\\
\frac{dM_{3}}{ds}=-2r_{0}M_{3}-\left(  M_{3}+M_{1}\right)  \alpha\rho,\\
\text{ }\\
\frac{dN_{1}}{ds}=-r_{0}N_{1},\\
\text{ }\\
\frac{dN_{2}}{ds}=-r_{0}N_{2}-\left(  N_{1}+N_{2}\right)  \alpha\rho,\\
\text{ }\\
M_{1}\left(  T\right)  =N_{1}\left(  T\right)  =1,\text{ }\\
N_{2}\left(  T\right)  =M_{2}\left(  T\right)  =M_{3}\left(  T\right)  =0.
\end{array}
\right.  \label{Eq47*}%
\end{equation}
Using the above system of ODEs, it is not difficult to verify that: For all
$s\in\left[  0,T\right]  $,%
\begin{equation}
M_{1}\left(  s\right)  -N_{1}\left(  s\right)  ^{2}=0\nonumber
\end{equation}
and%
\[
M_{3}\left(  s\right)  -N_{1}\left(  s\right)  N_{2}\left(  s\right)  =0.
\]
Consequently, it follows from (\ref{Eq47}) that for all $s\in\left[
0,T\right]  ,$%
\[
\alpha\left(  s\right)  =\frac{\mu\left(  N_{1}\left(  s\right)  +N_{2}\left(
s\right)  \right)  ^{-1}\rho\left(  s\right)  }{\left(  \sigma\left(
s\right)  ^{2}+%
{\displaystyle\int\limits_{E}}
\phi\left(  s,e\right)  ^{2}\vartheta\left(  de\right)  \right)  }=\mu\left(
N_{1}\left(  s\right)  +N_{2}\left(  s\right)  \right)  ^{-1}\kappa\left(
s\right)  ,
\]
where $\kappa\left(  s\right)  $ is as introduced in (\ref{Eq40*}).

Using the above expressions for $\alpha\left(  \cdot\right)  $, we can easily
solve the system of ODEs (\ref{Eq47*}), whose solutions are%
\begin{equation}
N_{1}\left(  s\right)  =e^{\int_{s}^{T}r_{0}\left(  \tau\right)  d\tau
}\text{,} \label{Eq60}%
\end{equation}%
\begin{equation}
N_{2}\left(  s\right)  =\mu e^{\int_{s}^{T}r_{0}\left(  \tau\right)  d\tau
}\int_{s}^{T}e^{-\int_{\tau}^{T}r_{0}\left(  \kappa\right)  d\kappa}%
\kappa\left(  \tau\right)  \rho\left(  \tau\right)  d\tau\text{,} \label{Eq61}%
\end{equation}%
\begin{equation}
M_{1}\left(  s\right)  =e^{2\int_{s}^{T}r_{0}\left(  \tau\right)  d\tau
}\text{,} \label{Eq62}%
\end{equation}%
\begin{equation}
M_{3}\left(  s\right)  =\mu e^{2\int_{s}^{T}r_{0}\left(  \tau\right)  d\tau
}\int_{s}^{T}e^{-\int_{\tau}^{T}r_{0}\left(  \kappa\right)  d\kappa}%
\kappa\left(  \tau\right)  \rho\left(  \tau\right)  d\tau\label{Eq63}%
\end{equation}
and%
\begin{align}
M_{2}\left(  s\right)   &  =\mu e^{\int_{s}^{T}\chi\left(  \tau\right)  d\tau
}\int_{s}^{T}e^{-\int_{\tau}^{T}\chi\left(  \kappa\right)  d\kappa}%
\kappa\left(  \tau\right)  \rho\left(  \tau\right)  \left\{  2\frac
{N_{1}\left(  \tau\right)  N_{2}\left(  \tau\right)  }{N_{1}\left(
\tau\right)  +N_{2}\left(  \tau\right)  }\right.  \text{\ }\nonumber\\
&  +\left.  \frac{N_{1}^{2}\left(  \tau\right)  +2N_{1}\left(  \tau\right)
N_{2}\left(  \tau\right)  }{\left\vert N_{1}\left(  \tau\right)  +N_{2}\left(
\tau\right)  \right\vert ^{2}}\mu\right\}  d\tau, \label{Eq64}%
\end{align}
where, $\chi\left(  \tau\right)  \equiv2r_{0}\left(  \tau\right)  +\left(
\frac{N_{1}\left(  \tau\right)  +N_{2}\left(  \tau\right)  +\mu}{N_{1}\left(
\tau\right)  +N_{2}\left(  \tau\right)  }\right)  ^{2}\kappa\left(
\tau\right)  \rho\left(  \tau\right)  -\kappa\left(  \tau\right)  \rho\left(
\tau\right)  $.

By (\ref{Eq46}) we get%
\begin{equation}
\hat{\varphi}^{s}\left(  z\right)  =\mu\left(  N_{1}\left(  s\right)
+N_{2}\left(  s\right)  \right)  ^{-1}\kappa\left(  s\right)  z\text{, for all
}\left(  s,z\right)  \in\left[  0,T\right]  \times%
\mathbb{R}
\text{.} \label{Eq48**}%
\end{equation}

Substituting this into the wealth equation results%
\[
\left\{
\begin{array}
[c]{l}%
d\hat{X}^{t,y}\left(  s\right)  =\left\{  r_{0}\left(  s\right)  +\mu\left(
N_{1}\left(  s\right)  +N_{2}\left(  s\right)  \right)  ^{-1}\kappa\left(
s\right)  \rho\left(  s\right)  \right\}  \hat{X}^{t,y}\left(  s\right)  ds\\
\text{ \ \ \ \ \ \ \ \ \ \ }+\mu\left(  N_{1}\left(  s\right)  +N_{2}\left(
s\right)  \right)  ^{-1}\kappa\left(  s\right)  \hat{X}^{t,y}\left(  s\right)
\sigma\left(  s\right)  dW\left(  s\right) \\
\text{ \ \ \ \ \ \ \ \ \ \ }%
\begin{array}
[c]{r}%
+%
{\displaystyle\int_{E}}
\mu\left(  N_{1}\left(  s\right)  +N_{2}\left(  s\right)  \right)  ^{-1}%
\kappa\left(  s\right)  \phi\left(  s,e\right)  \hat{X}_{-}^{t,y}\left(
s\right)  \tilde{N}\left(  ds,de\right)  \text{,}\\
\text{for }s\in\left[  t,T\right]  \text{,}%
\end{array}
\\
\hat{X}^{t,y}\left(  t\right)  =y\text{.}%
\end{array}
\right.
\]
The above SDE has a unique solution $\hat{X}^{t,y}\left(  \cdot\right)
\in\mathcal{S}_{%
\mathcal{F}%
}^{4}\left(  t,T;%
\mathbb{R}
\right)  $, for any $\left(  t,y\right)  \in\left[  0,T\right]  \times%
\mathbb{R}
$. So the closed-loop strategy $\hat{\varphi}\left(  \cdot,\cdot\right)  $
defined by (\ref{Eq48**}) is admissible. Moreover, it is not difficult to
check that Assumptions (i), (ii) and (iii) in Theorem \ref{result2} are
satisfied. Hence, $\hat{\varphi}\left(  \cdot,\cdot\right)  $ is an open-loop
equilibrium strategy.\eop

\begin{remark}
Let $N_{1}\left(  \cdot\right)  ,$ $N_{2}\left(  \cdot\right)  ,$
$M_{1}\left(  \cdot\right)  ,$ $M_{2}\left(  \cdot\right)  ,$ $M_{3}\left(
\cdot\right)  $ be the deterministic functions given by (\ref{Eq60}%
)-(\ref{Eq64}),$\ $respectively. Then the open-loop equilibrium control can be
represented by%
\[%
\begin{array}
[c]{r}%
\hat{u}^{x_{0}}\left(  s\right)  =\mu e^{-\int_{s}^{T}r_{0}\left(
\tau\right)  d\tau}\left(  1+\mu\int_{s}^{T}e^{-\int_{\tau}^{T}r_{0}\left(
v\right)  dv}\kappa\left(  \tau\right)  \rho\left(  \tau\right)  d\tau\right)
^{-1}\kappa\left(  s\right)  \hat{X}_{-}^{x_{0}}\left(  s\right)  \text{,}\\
\text{for }s\in\left[  0,T\right]  \text{,}%
\end{array}
\]
where the equilibrium wealth process is explicitly given by%
\begin{align*}
\hat{X}^{x_{0}}\left(  s\right)   &  =x_{0}\exp\left\{  \int_{0}^{s}\left(
r_{0}\left(  \tau\right)  +\mu L\left(  \tau\right)  \kappa\left(
\tau\right)  \rho\left(  \tau\right)  -\frac{1}{2}\left(  \mu L\left(
\tau\right)  \kappa\left(  \tau\right)  \sigma\left(  \tau\right)  \right)
^{2}\right)  d\tau\right. \\
&  +\int_{0}^{s}%
{\displaystyle\int_{E}}
\left\{  \ln\left(  1+\mu L\left(  \tau\right)  \kappa\left(  \tau\right)
\phi\left(  \tau,e\right)  \right)  -\mu L\left(  \tau\right)  \kappa\left(
\tau\right)  \phi\left(  \tau,e\right)  \right\}  \vartheta\left(  de\right)
d\tau\\
&
\begin{array}
[c]{r}%
+\int\limits_{0}^{s}\mu L\left(  \tau\right)  \kappa\left(  \tau\right)
\sigma\left(  \tau\right)  dW\left(  \tau\right)  \left.  +\int\limits_{t}^{s}%
{\displaystyle\int_{E}}
\left\{  \ln\left(  1+\mu L\left(  \tau\right)  \kappa\left(  \tau\right)
\phi\left(  \tau,e\right)  \right)  \right\}  \widetilde{N}\left(
d\tau,de\right)  \right\}  \text{,}\\
\text{for }s\in\left[  0,T\right]  \text{,}%
\end{array}
\end{align*}
with $L\left(  s\right)  \equiv e^{-\int_{s}^{T}r_{0}\left(  \tau\right)
d\tau}\left(  1+\mu\int_{s}^{T}e^{-\int_{\tau}^{T}r_{0}\left(  v\right)
dv}\kappa\left(  \tau\right)  \rho\left(  \tau\right)  d\tau\right)  ^{-1}$.
Moreover, simple calculations show that%
\begin{equation}
\mathbf{J}\left(  t,\hat{X}^{x_{0}}\left(  t\right)  ;\hat{\varphi}\left(
\cdot,\cdot\right)  \right)  =\left(  M_{2}\left(  t\right)  -N_{2}\left(
t\right)  ^{2}-2\mu_{1}\left(  N_{2}\left(  t\right)  +N_{1}\left(  t\right)
\right)  \right)  \frac{\hat{X}^{x_{0}}\left(  t\right)  ^{2}}{2}%
\text{,}\nonumber
\end{equation}
for any $t\in\left[  0,T\right]  $.
\end{remark}

\begin{remark}
If the modeling framework is without jumps, then the function $\kappa\left(
s\right)  $ introduced in (\ref{Eq40*}) reduces to%
\[
\kappa\left(  s\right)  :=\frac{\rho\left(  s\right)  }{\sigma\left(
s\right)  ^{2}}\text{.}%
\]
Accordingly, the open-loop equilibrium control $\hat{u}^{x_{0}}\left(
\cdot\right)  $ takes the form%
\[%
\begin{array}
[c]{r}%
\hat{u}^{x_{0}}\left(  s\right)  =\mu e^{-\int_{s}^{T}r_{0}\left(
\tau\right)  d\tau}\left(  1+\mu\int_{s}^{T}e^{-\int_{\tau}^{T}r_{0}\left(
v\right)  dv}\dfrac{\rho\left(  \tau\right)  ^{2}}{\sigma\left(  \tau\right)
^{2}}d\tau\right)  ^{-1}\frac{\rho\left(  s\right)  }{\sigma\left(  s\right)
^{2}}\hat{X}^{x_{0}}\left(  s\right)  \text{,}\\
\text{for }s\in\left[  0,T\right]  \text{,}%
\end{array}
\]
and the equilibrium wealth process becomes%
\begin{align*}
\hat{X}^{x_{0}}\left(  s\right)   &  =x_{0}\exp\left\{  \int_{0}^{s}\left(
r_{0}\left(  \tau\right)  +\mu L^{0}\left(  \tau\right)  \kappa\left(
\tau\right)  \rho\left(  \tau\right)  -\frac{1}{2}\left(  \mu L^{0}\left(
\tau\right)  \kappa\left(  \tau\right)  \sigma\left(  s\right)  \right)
^{2}\right)  d\tau\right. \\
&
\begin{array}
[c]{c}%
\left.  +\int_{0}^{s}\mu L^{0}\left(  \tau\right)  \kappa\left(  \tau\right)
\sigma\left(  s\right)  dW\left(  \tau\right)  \right\}  \text{, for }%
s\in\left[  0,T\right]  \text{,}%
\end{array}
\end{align*}
with $L^{0}\left(  s\right)  \equiv e^{-\int_{s}^{T}r_{0}\left(  \tau\right)
d\tau}\left(  1+\mu\int_{s}^{T}e^{-\int_{\tau}^{T}r_{0}\left(  v\right)
dv}\dfrac{\rho\left(  \tau\right)  ^{2}}{\sigma\left(  \tau\right)  ^{2}}%
d\tau\right)  ^{-1}$. This essentially coincide with the result derived by Hu
et al. [\cite{Huetal2011}, Subsection 5.4.1] by solving a flow of FBSDEs.
\end{remark}

\section{Conclusion and future work}

In this paper we have investigated open-loop equilibrium controls to a general
class of time-inconsistent stochastic optimal control problems. Compared with
the existing literature, the novelty of this paper is that we propose a new
method that enables us to construct open-loop equilibrium controls by solving
a deterministic coupled system of IPDEs. In Remark \ref{Rem1}, we have
listed\ several advantages of this new approach. Moreover, we illustrated our
main result by solving a mean-variance portfolio selection problem and the
solution we obtained coincides with that obtained in \cite{Sun} and
\cite{Huetal2011} by solving a flow of FBSDEs. Some obvious open research
problems are the following.

\begin{description}
\item[(i)] The task of proving existence and/or uniqueness of solutions to the
system of IPDEs (\ref{Eq38}) to be technically extremely difficult. We have no
idea about how to proceed so we leave it for future research.

\item[(ii)] We characterized open-loop equilibriums via a sufficient condition
only, it remains to derive the necessary condition and investigate the
uniqueness of the equilibrium solution.

\item[(iii)] Assumptions (i) and (ii) in Theorem \ref{result2} may seem
restrictive conditions. We hope that in our future publications, these
conditions can be made weaker.

\item[(iv)] Finally, the present work can be extended in several ways. For
example it is interesting to generalize our results to the case where the
running cost also depend on conditional expectations for the state process and
control process.
\end{description}

\end{document}